\documentclass[12pt,a4paper]{article}
\usepackage{amssymb,latexsym,times,standard,graphicx}

\textheight 
            675pt
\textwidth 
           460pt

\newcommand{\dk}{\, \mathrm{d}k}
\newcommand{\dr}{\, \mathrm{d}r}

\newcommand{\dx}{\, \mathrm{d}x}
\newcommand{\dtildex}{\, \mathrm{d}\tilde{x}}
\newcommand{\ds}{\, \mathrm{d}s}
\newcommand{\dy}{\, \mathrm{d}y}
\newcommand{\dtildey}{\, \mathrm{d}\tilde{y}}
\newcommand{\dz}{\, \mathrm{d}z}
\newcommand{\dtildez}{\, \mathrm{d}\tilde{z}}

\newcommand{\dtheta}{\, \mathrm{d}\theta}

\newcommand{\sm}{\!\setminus\!}

\newcommand{\EE}{\mathcal E}
\newcommand{\FF}{\mathcal F}
\newcommand{\GG}{\mathcal G}
\newcommand{\II}{\mathcal I}
\newcommand{\JJ}{\mathcal J}
\newcommand{\KK}{\mathcal K}
\newcommand{\LL}{\mathcal L}
\newcommand{\MM}{\mathcal M}

\newcommand{\RR}{\mathcal R}
\renewcommand{\SS}{\mathcal S}

\newcommand{\e}{\mathrm{e}}
\renewcommand{\i}{\mathrm{i}}

\newcommand{\supp}{\mathrm{supp}\,}
\newcommand{\dist}{\mathrm{dist}}

\newcommand{\qed}{\hfill$\Box$\bigskip}

\newcommand{\eqn}[1]{(\ref{#1})}

\newcommand{\nn}{|{\mskip-2mu}|{\mskip-2mu}|}
\newcommand{\uunderbrace}[1]{\underbrace{\vphantom{\bigg(}#1}}

\newtheorem{theorem}{Theorem}[section]
\newtheorem{lemma}[theorem]{Lemma}
\newtheorem{definition}[theorem]{Definition}
\newtheorem{proposition}[theorem]{Proposition}
\newtheorem{corollary}[theorem]{Corollary}
\newtheorem{remark}[theorem]{Remark}

\begin{document}
\newcounter{count}

\title{Existence and conditional energetic stability of three-dimensional fully localised solitary gravity-capillary water waves}

\author{B. Buffoni\thanks{Section de math\'{e}matiques (IACS), \'{E}cole polytechnique f\'{e}d\'{e}rale, 1015 Lausanne, Switzerland},\quad
M. D. Groves\thanks{FR 6.1 - Mathematik, Universit\"{a}t des Saarlandes, Postfach 151150, 66041 Saarbr\"{u}cken, Germany;
Department of Mathematical Sciences, Loughborough University, Loughborough, Leics, LE11 3TU, UK},
\quad S. M. Sun\thanks{Department of Mathematics, Virginia Polytechnic Institute and State University,
Blacksburg, VA 24061, USA},
\quad E. Wahl\'{e}n\thanks{Centre for Mathematical Sciences, Lund University, P.O. Box 118, 22100 Lund, Sweden}}
\date{}
\maketitle{}

\begin{abstract}
In this paper we show that the hydrodynamic problem for three-dimensional water waves with strong
surface-tension effects admits a \emph{fully localised solitary wave} which decays to the undisturbed
state of the water in every horizontal direction. The proof is based upon the classical variational
principle that a solitary wave of this type is a critical point of the energy, which is given in
dimensionless coordinates by
$$\EE(\eta,\phi) = \int_{{\mathbb R}^2}\left\{ \frac{1}{2}\int_0^{1+\eta} (\phi_x^2+\phi_y^2+\phi_z^2) \dy
+ \frac{1}{2}\eta^2 + \beta[\sqrt{1+\eta_x^2+\eta_z^2}-1]\right\}\dx\dz,$$
subject to the constraint that the momentum
$$\II(\eta,\phi) = \int_{{\mathbb R}^2} \eta_x \phi|_{y=1+\eta} \dx\dz$$
is fixed; here $\{(x,y,z): x,z \in {\mathbb R},\ y \in (0,1+\eta(x,z))\}$ is the fluid domain, $\phi$ is the
velocity potential and $\beta>1/3$ is the Bond number. These functionals are studied locally for $\eta$
in a neighbourhood of the origin in $H^3({\mathbb R}^2)$.

We prove the existence of a minimiser of $\EE$ subject to the constraint $\II=2\mu$, where $0<\mu \ll 1$.
The existence of a small-amplitude solitary wave is thus assured, and since $\EE$ and $\II$ are both conserved quantities
a standard argument may be used to establish the stability of the set $D_\mu$ of minimisers as a whole.
`Stability' is however understood in a qualified sense due to the lack of a global well-posedness theory for
three-dimensional water waves. We show that solutions to the evolutionary problem
starting near $D_\mu$ remain close to $D_\mu$ in a suitably defined energy space over their
interval of existence; they may however explode in finite time due to higher-order derivatives becoming
unbounded.
\end{abstract}

\newpage

\section{Introduction}

\subsection{The hydrodynamic problem}

The classical \emph{three-dimensional gravity-capillary water wave
problem} concerns the irrotational flow of a perfect fluid of unit
density subject to the forces of gravity and surface tension. The
fluid motion is described by the Euler equations in a domain bounded
below by a rigid horizontal bottom $\{y=0\}$ and above by a free
surface $\{y=h+\eta(x,z,t)\}$, where $h$
denotes the depth of the water in its undisturbed state and the
function $\eta$ depends upon the two horizontal spatial directions
$x$, $z$ and time $t$. In terms of an Eulerian velocity potential $\phi$,
the mathematical problem is to solve Laplace's equation
$$
\phi_{xx} + \phi_{yy} + \phi_{zz} = 0, \qquad\qquad 0<y<h+\eta
$$
with boundary conditions
\begin{eqnarray*}
\phi_{y} & = & \parbox{100mm}{$0$,} y=0, \\
\eta_t & = & \parbox{100mm}{$\phi_{y} - \eta_x\phi_x - \eta_z\phi_z$,}  y=h+\eta, \\
\phi_t & = & -\frac{1}{2}(\phi_x^2+\phi_{y}^2+\phi_z^2) -g\eta \nonumber \\
& &
\parbox{100mm}{$\displaystyle\qquad\mbox{} + \sigma\left[\frac{\eta_x}{\sqrt{1+\eta_x^2+\eta_z^2}}\right]_x
+ \sigma\left[\frac{\eta_z}{\sqrt{1+\eta_x^2+\eta_z^2}}\right]_z$,}
 y=h+\eta,
\end{eqnarray*}
in which $g$ is the acceleration due to gravity and
$\sigma>0$ is the coefficient of surface tension (see, for example, Stoker
\cite[\S\S 1, 2.1]{Stoker}). Introducing the dimensionless variables
$$(x^\prime,y^\prime,z^\prime)=\frac{1}{h}(x,y,z),\qquad t^\prime=\left(\frac{g}{h}\right)^{1/2},$$
$$\eta^\prime(x^\prime,z^\prime,t^\prime) = \frac{1}{h}\eta(x,z,t),\qquad
\phi^\prime(x^\prime,y^\prime,z^\prime,t^\prime) = \frac{1}{(gh)^{3/2}}\phi(x,y,z,t),$$
one obtains the equations
\begin{equation}
\phi_{xx} + \phi_{yy} + \phi_{zz} = 0, \qquad\qquad 0<y<1+\eta
\label{Laplace's equation}
\end{equation}
with boundary conditions
\begin{eqnarray}
\phi_{y} & = & \parbox{100mm}{$0$,} y=0, \label{BC 1} \\
\eta_t & = & \parbox{100mm}{$\phi_{y} - \eta_x\phi_x - \eta_z\phi_z$,}  y=1+\eta,
\label{BC 2} \\
\phi_t & = & -\frac{1}{2}(\phi_x^2+\phi_{y}^2+\phi_z^2) -\eta \nonumber \\
& &
\parbox{100mm}{$\displaystyle\qquad\mbox{} + \beta\left[\frac{\eta_x}{\sqrt{1+\eta_x^2+\eta_z^2}}\right]_x
+ \beta\left[\frac{\eta_z}{\sqrt{1+\eta_x^2+\eta_z^2}}\right]_z$,} y=1+\eta,
\label{BC 3}
\end{eqnarray}
where $\beta=\sigma/gh^2$ and the primes have been dropped for notational simplicity.

\emph{Steady waves} are water waves which travel in a distinguished horizontal
direction with constant speed and without change of shape; without loss of generality we
assume that the waves propagate from right to left in the $x$-direction with speed $\nu$, so that
$\eta(x,z,t)=\eta(x+\nu t,z)$ and $\phi(x,y,z,t)=\phi(x+\nu t,y,z)$. In this paper we study \emph{fully
localised solitary waves,} that is steady waves with the property that $\eta(x+\nu t,z) \rightarrow 0$
as $|(x+\nu t,z)| \rightarrow \infty$; in particular we consider the parameter regime $\beta>1/3$
corresponding to strong surface tension. Interest in this parameter regime stems from the
celebrated Kadomtsev \& Petviashvili (KP-I) equation (a model for long water waves with
strong surface tension and a preferred direction of propagation), which admits a fully localised
solitary-wave solution given by the explicit formula
\begin{equation}
u(x,z)=-8\frac{3-x^2+z^2}{(3+x^2+z^2)^2}
\label{Explicit KP formula}
\end{equation}
in a frame of reference moving with the wave
(see Ablowitz \& Segur \cite{AblowitzSegur79}); the variable $u$ is supposed to approximate
the free surface of the water via the relationship
\begin{equation}
\eta(x,z)=\mu^2 u\Bigg(\frac{\mu x}{2(\beta-1/3)^{1/2}},\mu^2 z\Bigg) + O(\mu^3),
\label{AS explicit formula}
\end{equation}
where $\mu$ is a small parameter associated with the weakly nonlinear scaling limit (see below).

\begin{figure}[h]
\hspace{4cm}\includegraphics[width=7cm]{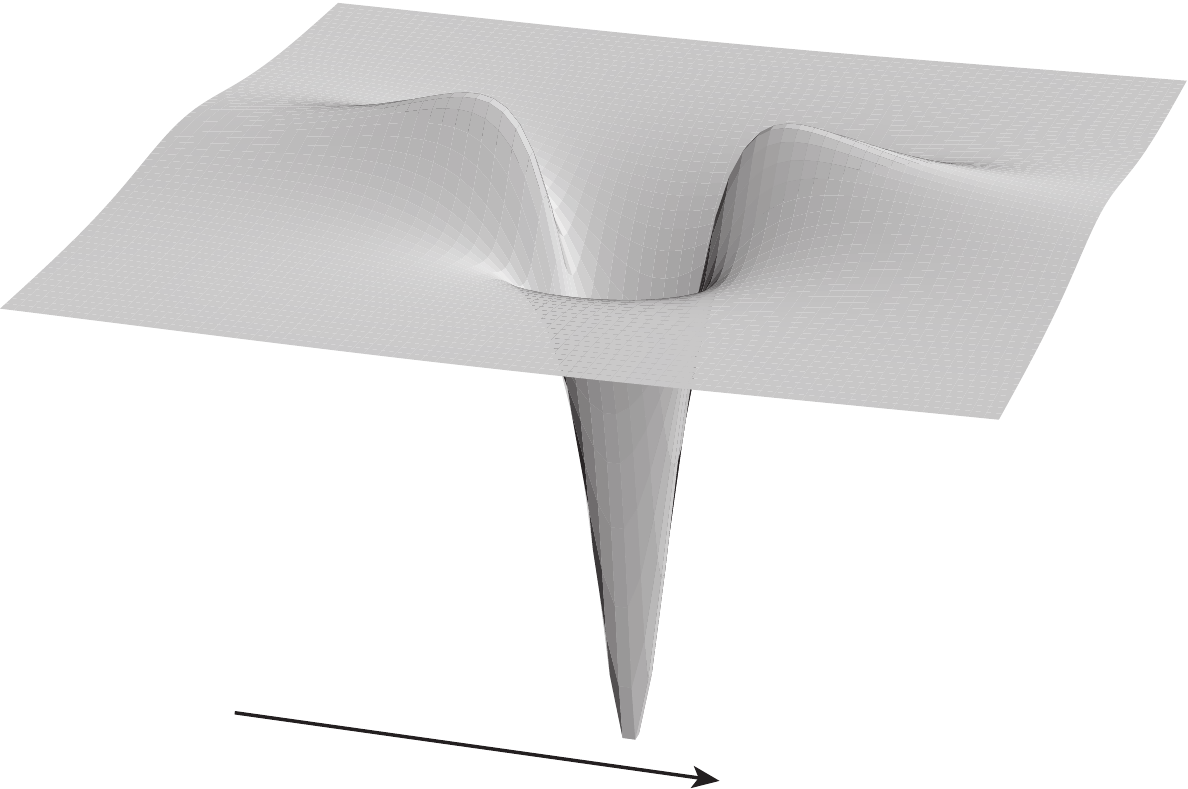}
{\it
\caption{A fully localised solitary wave; the
arrow shows the direction of wave propagation. \label{Fully localised solitary wave}}}
\end{figure}

Fully localised solitary-wave solutions to other models for three-dimensional
surface-tension dominated flows have also been studied. Mathematical existence theories
for generalised KP-I equations were given by
Wang \& Willem \cite{WangWillem96},
de Bouard \& Saut \cite{deBouardSaut97b}
and Pankov \& Pfl\"{u}ger \cite{PankovPflueger00},
and for the Benny-Luke equation (an isotropic version of the KP-I equation)
by Pego \& Quintero \cite{PegoQuintero99} (see Berger \& Milewski
\cite{BergerMilewski00} for numerical computations). The existence of a fully localised
solitary-wave solution to \eqn{Laplace's equation}--\eqn{BC 3} in this parameter regime was recently
established rigorously by Groves \& Sun \cite{GrovesSun08} and
computed numerically by Parau, Vanden-Broeck \& Cooker \cite{ParauVandenBroeckCooker05a}
(see Figure \ref{Fully localised solitary wave}). In the present paper we give an alternative,
more natural existence theory. Groves \& Sun work with Fourier-multiplier operators in $L^p$-based
function spaces for $p>1$ (which require detailed analysis), and use a local reduction technique
to reduce the problem to a single semilinear equation. On the other hand the present
theory employs $L^2$-based function spaces and tackles the original hydrodynamic problem directly;
furthermore, it yields information concerning the stability of the fully localised solitary wave.

The conserved quantities
$$\EE(\eta,\phi) = \int_{{\mathbb R}^2}\left\{ \frac{1}{2}\int_0^{1+\eta} (\phi_x^2+\phi_y^2+\phi_z^2) \dy
+ \frac{1}{2}\eta^2 + \beta[\sqrt{1+\eta_x^2+\eta_z^2}-1]\right\}\dx\dz,$$
$$\II(\eta,\phi) = \int_{{\mathbb R}^2} \eta_x \phi|_{y=1+\eta} \dx\dz$$
of \eqn{Laplace's equation}--\eqn{BC 3},
which represent respectively the energy and momentum of a wave in the $x$-direction
(see Benjamin \& Olver \cite{BenjaminOlver82}), are the key to our existence theory. A
fully localised solitary wave is characterised as a critical point of the energy subject to
the constraint of fixed momentum; it is therefore a critical point of the functional $\EE-\nu\II$,
where the Lagrange multiplier $\nu$ gives the speed of the wave.
Furthermore, Benjamin \cite{Benjamin74} noted that constrained \emph{minimisers}
should be stable, and in this paper we construct an existence and stability theory for
fully localised solitary waves using Benjamin's principle.

\subsection{Constrained minimisation and conditional energetic stability}

The above formulation of the hydrodynamic problem has the disadvantage that it is
posed in the \emph{a priori} unknown domain $\{0<y<1+\eta(x,z,t)\}$. We overcome this
difficulty using the \emph{Dirichlet-Neumann operator} $G(\eta)$ introduced by Craig
\cite{Craig91} and formally defined as follows.
For fixed $\Phi=\Phi(x,z)$ solve the boundary-value problem
\begin{eqnarray*}
& & \parbox{6cm}{$\phi_{xx}+\phi_{yy}+\phi_{zz}=0,$}0<y<1+\eta, \\
& & \parbox{6cm}{$\phi = \Phi,$}y=1+\eta, \\
& & \parbox{6cm}{$\phi_y =0,$}y=0
\end{eqnarray*}
and define
\begin{eqnarray*}
G(\eta)\Phi & = & \sqrt{1+\eta_x^2+\eta_z^2}\left.\frac{\partial \phi}{\partial n}\right|_{y=1+\eta} \\
& = & \phi_y - \eta_x\phi_x - \eta_z\phi_z\Big|_{y=1+\eta}.
\end{eqnarray*}
In terms of the variables $\eta$ and $\Phi$ the energy and momentum functionals are
given by the convenient formulae
\begin{equation}
\EE(\eta,\Phi) = \int_{{\mathbb R}^2}\left\{ \frac{1}{2}\Phi G(\eta)\Phi
+ \frac{1}{2}\eta^2 + \beta[\sqrt{1+\eta_x^2+\eta_z^2}-1]\right\}\dx\dz,
\label{Definition of E}
\end{equation}
\begin{equation}
\II(\eta,\Phi) = \int_{{\mathbb R}^2} \eta_x \Phi \dx\dz, \label{Definition of I}
\end{equation}
and in this paper we study the problem of finding minimisers of $\EE$ subject to the
constraint
$$\II(\eta,\Phi) = 2\mu,$$
where $\mu$ is a small positive number.

We begin our study in Section \ref{FA setting} by identifying function spaces
in which equations \eqn{Definition of E}, \eqn{Definition of I} define
analytic functionals, using in particular the theory of analytic operators
reported by Buffoni \& Toland \cite{BuffoniToland}. We take
$\eta$ in a neighbourhood $U=B_M(0)$ of the origin in $H^3({\mathbb R}^2)$, while the
elements of the function space $H_\star^{1/2}({\mathbb R}^2)$ for $\Phi$ are traces of potential
functions in the fluid domain (see Section \ref{ND operator} for a precise definition).
It follows from the following result, which is proved in Section \ref{Analyticity of NDO},
that $\EE$ and $\II$ are smooth functionals $U \times H_\star^{1/2}({\mathbb R}^2) \mapsto {\mathbb R}$.

\begin{lemma}
The mapping $W^{1,\infty}({\mathbb R}^2) \rightarrow
\LL(H_\star^{1/2}({\mathbb R}^2), (H_\star^{1/2}({\mathbb R}^2))^\prime)$ given by
$\eta \mapsto (\Phi \mapsto G(\eta)\Phi)$
is analytic at the origin and $\Phi \rightarrow G(\eta)\Phi$ is an isomorphism
$H_\star^{1/2}({\mathbb R}^2) \rightarrow (H_\star^{1/2}({\mathbb R}^2))^\prime$
for each $\eta$ in a neighbourhood of the origin.
\end{lemma}

The above variational principle has been used by several authors in existence theories
for three-dimensional steady water waves. Groves \& Sun \cite{GrovesSun08} obtained fully localised
solitary waves as critical points of the functional $\EE-\nu\II$, while Groves \& Haragus \cite{GrovesHaragus03}
interpreted $\EE-\nu\II$ as an action functional to derive a formulation
of the hydrodynamic problem as an infinite-dimensional spatial Hamiltonian system with
a rich solution set. Finally, Craig \& Nicholls \cite{CraigNicholls00} developed an existence
theory for doubly periodic steady waves as critical points of $\EE$ on level sets of $\II$
(here the integrals over ${\mathbb R}^2$ are replaced with integrals over the periodic
domain). These papers all use reduction methods to reduce the quasilinear
problem under investigation to simpler semilinear or finite-dimensional problems.

In this paper we find constrained minimisers of $\EE$ using a genuinely infinite-dimensional
method. Our main existence result is stated in the following theorem (see Section \ref{Outline of J theory}
for an overview of the proof).

\begin{theorem} \hspace{1cm} \label{Result for constrained minimisation - Intro}
\begin{list}{(\roman{count})}{\usecounter{count}}
\item
The set $D_\mu$ of minimisers of $\EE$ over the set
$$S_\mu=\{(\eta,\Phi) \in U \times H_\star^{1/2}({\mathbb R}^2): \II(\eta,\Phi)=2\mu\}$$
is non-empty. The corresponding solitary waves are subcritical, that is their
dimensionless speed is less than unity.
\item
Suppose that $\{(\eta_n,\Phi_n)\} \subset S_\mu$ is a minimising sequence for $\EE$ with
the property that
$$
\sup_{n\in{\mathbb N}} \|\eta_n\|_3 < M.
$$
There exists a sequence $\{(x_n,z_n)\} \subset {\mathbb R}^2$ with the property that
a subsequence of\linebreak $\{\eta_n(x_n+\cdot,z_n+\cdot),\Phi_n(x_n+\cdot,z_n+\cdot)\}$
converges  in
$H^r({\mathbb R}^2) \times H_\star^{1/2}({\mathbb R}^2)$, $0 \leq r < 3$ to a function in $D_\mu$.
\end{list}
\end{theorem}

We discuss the stability of the set $D_\mu$ in Section \ref{Stability}.
The usual informal interpretation of the statement that a set $X$ of solutions to an initial-value problem
is `stable' is that a solution which begins close to a solution in $X$
remains close to a solution in $X$ at all subsequent times. Implicit in this statement
is the assumption that the initial-value problem is globally well-posed, that is every pair
$(\eta_0,\Phi_0)$ in an appropriately chosen set is indeed the initial datum
of a unique solution $t \mapsto (\eta(t),\Phi(t))$, $t \in [0,\infty)$.
At present there is no global well-posedness theory for three-dimensional
water waves, and we work instead under the following assumption (see
Alazard, Burq \& Zuly \cite{AlazardBurqZuily11} for results of this kind).
\begin{quote}
({\bf Well-posedness assumption}) There exists a subset  $S$
of $U \times H_\star^{1/2}({\mathbb R}^2)$ with the following properties.
\begin{list}{(\roman{count})}{\usecounter{count}}
\item
The closure of $S \sm D_\mu$ in $L^2({\mathbb R}^2)$ has a non-empty intersection with
$D_\mu$.
\item
For each $(\eta_0,\Phi_0) \in S$
there exists $T>0$ and a continuous function $t \mapsto (\eta(t), \Phi(t)) \in 
U \times H_\star^{1/2}({\mathbb R}^2)$, $t \in [0, T]$ such that $(\eta(0),\Phi(0)) = (\eta_0,\Phi_0)$,
$$\EE(\eta(t),\Phi(t)) = \EE(\eta_0,\Phi_0),\ \II(\eta(t),\Phi(t))=\II(\eta_0,\Phi_0), \qquad t \in [0,T]$$
and 
$$\sup_{t \in [0,T]} \|\eta(t)\|_3 < M.$$
\end{list}
\end{quote}
It is a general principle that the solution set of a constrained minimisation
problem constitutes a stable set of solutions of the corresponding
initial-value problem (e.g.\ see Cazenave \& Lions \cite{CazenaveLions82} and
de Bouard \& Saut \cite{deBouardSaut96}, Liu \& Wang \cite{LiuWang97} for applications
to generalised KP equations). Combining this general principle
with our well-posedness assumption, we obtain the following stability result.

\begin{theorem} \label{CES - Intro}
Choose $r \in [0,3)$. For each $\varepsilon>0$  there exists $\delta>0$ such that
$$(\eta_0,\Phi_0) \in S,\ \dist((\eta_0,\Phi_0), D_\mu) < \delta \quad \Rightarrow \quad
\dist((\eta(t),\Phi(t)), D_\mu)<\varepsilon,$$
for $t\in[0,T]$,
where `$\dist$' denotes the distance in $H^r({\mathbb R}^2) \times H_\star^{1/2}({\mathbb R}^2)$.
\end{theorem}

This result is a statement of
the \emph{conditional, energetic stability of the set $D_\mu$}. Here
\emph{energetic} refers to the fact that the
distance in the statement of stability is measured in the `energy space'
$H^r({\mathbb R}^2) \times H_\star^{1/2}({\mathbb R}^2)$, while \emph{conditional}
alludes to the well-posedness assumption. Note that the solution $t \mapsto (\eta(t),\Phi(t))$
may exist in a smaller space over the interval $[0,T]$, at each instant of which it remains close
(in energy space) to a solution in $D_\mu$. Furthermore, Theorem \ref{CES - Intro}
is a statement of the stability of the \emph{set} of constrained minimisers $D_\mu$;
establishing the uniqueness of the constrained minimiser would imply
that $D_\mu$ consists of translations of a single solution, so that the statement
that $D_\mu$ is stable is equivalent to classical orbital stability of this unique
solution (Benjamin \cite{Benjamin74}).
The phrase `conditional, energetic stability' was introduced by
Mielke \cite{Mielke02} in his study of the stability of two-dimensional
solitary water waves with strong surface tension using dynamical-systems
methods and further developed in the context of variational methods
for two-dimensional solitary waves by Buffoni \cite{Buffoni04b,Buffoni05}.

\subsection{The minimisation problem} \label{Outline of J theory}
We tackle the problem of finding minimisers of $\EE(\eta,\Phi)$ subject to the constraint
$\II(\eta,\Phi)=2\mu$ in two steps.
\begin{enumerate}
\item \emph{Fix $\eta \neq 0$ and minimise $\EE(\eta,\cdot)$ over
$T_\mu=\{\Phi \in H^{1/2}_\star({\mathbb R}^2): \II(\eta,\Phi)=2\mu\}$.}
This problem (of minimising a quadratic functional over a linear manifold)
admits a unique global minimiser $\Phi_\eta$.
\item \emph{Minimise $\JJ_\mu(\eta):=\EE(\eta,\Phi_\eta)$ over $\eta \in U\sm\{0\}$.}
Because $\Phi_\eta$ minimises
$\EE(\eta,\cdot)$ over $T_\mu$ there exists a Lagrange multiplier $\lambda_\eta$ such that
$$G(\eta)\Phi_\eta = \lambda_\eta \eta_x,$$
and straightforward calculations show that
$\Phi_\eta = \lambda_\eta G(\eta)^{-1}\eta_x$, $\lambda_\eta = \mu/\LL(\eta)$
and
\begin{equation}
\JJ_\mu(\eta)=\KK(\eta)+\frac{\mu^2}{\LL(\eta)},
\label{Definition of J}
\end{equation}
where
\begin{equation}
\KK(\eta) = \int_{{\mathbb R}^2} \left\{\frac{1}{2}\eta^2 + \beta\sqrt{1+\eta_x^2+\eta_z^2}-\beta
\right\} \dx\dz, \label{Definition of K}
\end{equation}
\begin{equation}
\LL(\eta) = \frac{1}{2}\int_{{\mathbb R}^2} \eta_x G(\eta)^{-1} \eta_x \dx\dz.
\label{Definition of L}
\end{equation}
This computation also shows that the dimensionless speed of the solitary wave corresponding to
a constrained minimiser of $\EE(\eta,\Phi)$ is $\mu/\LL(\eta)$.
\end{enumerate}
The above two-step approach was introduced by Buffoni \cite{Buffoni04a} in a corresponding theory
for two-dimen\-sional solitary waves. Buffoni used a conformal mapping to transform $\KK$
and $\LL$ into simpler functionals and hence greatly simplified the analysis necessary to
show that $\JJ_\mu$ has a minimiser. Here we extend his method to our three-dimensional
problem, working directly with the functionals as given above. In this respect we note that
$\KK$ is analytic at the origin in the Sobolev space $H^r({\mathbb R}^2)$ for $r>2$, and it
follows from the following result, which is proved in Section \ref{Analyticity of NDO},
that  ${\mathcal L}$ is analytic at the origin in $H^r({\mathbb R}^2)$ for $r>5/2$.

\begin{lemma}
Suppose that $s>1$. The operator $K(\cdot): H^{s+3/2}({\mathbb R}^2) \rightarrow
\LL(H^{s+1}({\mathbb R}^2), H^s({\mathbb R}^2))$ given by the formula
$$K(\eta) = -\partial_x (G(\eta)^{-1}\partial_x)$$
is analytic at the origin.
\end{lemma}

The above comments show that $\JJ_\mu$ is a smooth functional in a punctured neighbourhood
of the origin in $H^r({\mathbb R}^2)$ for $r>5/2$; we seek minimisers of $\JJ_\mu$ in
the smaller space $H^3({\mathbb R}^2)$, taking advantage of the fact that $H^3({\mathbb R}^2)$
is locally compactly embedded in $H^r({\mathbb R}^2)$ for $r \in [0,3)$. Our main result for $\JJ_\mu$,
from which Theorem \ref{Result for constrained minimisation - Intro} is deduced in Section \ref{Stability},
is stated in the following theorem.

\begin{theorem} \label{Main result for J}
There exists a neighbourhood $U=B_M(0)$ of the origin in $H^3({\mathbb R}^2)$ with the following
properties.
\begin{list}{(\roman{count})}{\usecounter{count}}
\item
The set of minimisers of $\JJ_\mu$ on $U \sm \{0\}$ is non-empty.
\item
Suppose that $\{\eta_n\}$ is a minimising sequence for $\JJ_\mu$ over $U\sm\{0\}$ which satisfies
$$
\sup_{n\in{\mathbb N}} \|\eta_n\|_3 < M.
$$
There exists a sequence $\{(x_n,z_n)\} \subset {\mathbb R}^2$ with the property that
a subsequence of\linebreak $\{\eta_n(x_n+\cdot,z_n+\cdot)\}$ converges
in $H^r({\mathbb R}^2)$ for $r \in [0,3)$ to a critical point $\eta$
which minimises $\JJ_\mu$ on $U\sm\{0\}$.
\end{list}
\end{theorem}

The theorem is proved by first reducing the first assertion to a special case of the second; in so doing
one is immediately confronted by the unfavourable properties of the quadratic parts $\KK_2$ and $\LL_2$
of $\KK$ and $\LL$. The functionals are not coercive, in the sense that $(\KK_2)^{\frac{1}{2}}$ and $(\LL_2)^{\frac{1}{2}}$
are not bounded below by any constant multiple of the $H^3({\mathbb R}^2)$-norm, and furthermore
the functional
$$\LL_2(\eta) = \frac{1}{2} \int_{{\mathbb R}^2}\frac{k_1^2}{|k|^2}|k|\coth|k||\hat{\eta}|^2 \dk$$
is anisotropic and involves a Fourier multiplier which is
not smooth at the origin.
We proceed by introducing the penalised functional
$\JJ_{\rho,\mu}: H^3({\mathbb R}^2) \rightarrow {\mathbb R} \cup \{\infty\}$ defined by
$$\JJ_{\rho,\mu}(\eta) = \left\{\begin{array}{lll}\displaystyle \KK(\eta)+\frac{\mu^2}{\LL(\eta)}
+ \rho(\|\eta\|_3^2), & & u \in U\sm\{0\}, \\
\\
\infty, & & \eta \not\in U\sm\{0\},
\end{array}\right.$$
where $\rho: [0,M^2) \rightarrow {\mathbb R}$ is a smooth, increasing `penalisation' function which
explodes to infinity as $t \uparrow M^2$ and vanishes for $0 \leq t \leq \tilde{M}^2$; the number
$\tilde{M}$ is chosen very close to $M$. This functional enjoys a degree of coercivity and has the advantage that
a minimising sequence over $U\sm\{0\}$ does not approach the boundary of $U$.

Minimising sequences $\{\eta_n\}$ for $\JJ_{\rho,\mu}$,
which clearly satisfy $\sup \|\eta_n\|_3 < M$, are studied in detail in Section \ref{MS} with the help
of the concentration-compactness principle (Lions \cite{Lions84a,Lions84b}). The main difficulty
here lies in discussing the consequences of `dichotomy'. On the one hand the functional  ${\mathcal L}$
is nonlocal, so that a careful argument is required to show that $\eta_n$ splits into two parts
in the usual fashion (see Lemma \ref{Splitting properties 1}(iii) and Appendix D).
On the other hand no \emph{a priori} estimate is available to rule out `dichotomy'
at this stage; proceeding iteratively we find that minimising sequences can theoretically have profiles
with infinitely many `lumps'. In particular we show that $\{\eta_n\}$
asymptotically lies in the region unaffected by the penalisation (Corollary \ref{Asymptotically unpenalised})
and construct a special minimising sequence $\{\tilde{\eta}_n\}$ for $\JJ_{\rho,\mu}$
which lies in a neighbourhood of the origin with radius $O(\mu^{\frac{1}{2}})$ in $H^3({\mathbb R}^2)$
and satisfies $\|\JJ^\prime_\mu(\tilde{\eta}_n)\|_1 \rightarrow 0$ as $n \rightarrow \infty$ (Section
\ref{Special minimising sequence}). The fact that
the construction is independent of the choice of $\tilde{M}$ allows us to conclude that $\{\tilde{\eta}_n\}$
is also a minimising sequence for $\JJ_\mu$ over $U\sm\{0\}$.

The special minimising sequence $\{\tilde{\eta}_n\}$ is used in Section \ref{SSA} to establish the
\emph{strict sub-additivity} property
$$
c_{\mu_1+\mu_2} < c_{\mu_1} + c_{\mu_2}, \qquad \mu_1, \mu_2>0
$$
of the infimum $c_\mu$ of $\JJ_\mu$ over $U \sm \{0\}$. An argument given by
Buffoni \cite{Buffoni04a} shows how this property follows from the fact that the function
\begin{equation}
a \mapsto a^{-\frac{5}{2}}\MM_{a^2\mu}(a\eta_n), \qquad a \in [1,2], \label{Decreasing function of a - Intro}
\end{equation}
is decreasing and strictly negative, where $\{\eta_n\}$ is a minimising sequence
for $\JJ_\mu$ over $U\sm\{0\}$ and
$$\MM_\mu(\eta) := \JJ_\mu(\eta) - \KK_2(\eta) - \frac{\mu^2}{\LL_2(\eta)}$$
is the `nonlinear' part of $\JJ_\mu(\eta)$. The function
\eqn{Decreasing function of a - Intro} would clearly have the required property if $\MM_\mu(\eta_n)$ were
homogeneous cubic and negative; we therefore proceed by approximating $\MM_\mu(\eta_n)$ with an expression
of this kind.

Using the fact that every minimising
sequence $\{\eta_n\}$ satisfies $\LL_2(\eta_n)$, $\LL(\eta_n) \geq c\mu$,\linebreak
$\MM_\mu(\eta_n) \leq -c \mu^3$,
one finds that  the `cubic' part of $\MM_\mu(\eta_n)$ is $-(\mu/\LL_2(\eta_n))^2\LL_3(\eta_n)$ and that
$\MM_\mu(\eta_n)$ can be approximated by a cubic, negative expression provided that 
all other terms in $\MM_\mu(\eta_n)$ are $o(\mu^3)$.
The straightforward estimate $\|\eta_n\|_3^2=O(\mu)$ does not suffice for this purpose (the `quartic' part of $\MM_\mu(\eta_n)$
would for example be merely $O(\mu^2)$). Motivated by the expectation
that a critical point of $\JJ_\mu$, and hence the minimising sequence $\{\tilde{\eta}_n\}$,
should have the KP-I length scales, we prove that $\nn \tilde{\eta}_n \nn_\alpha^2
=O(\mu)$ for each $\alpha<1$, where
$$
\nn \eta \nn_\alpha^2 := \int_{{\mathbb R}^2} \left(1+\mu^{-6\alpha}|k|^6+\mu^{-4\alpha}\frac{k_2^4}{|k|^4}\right)|\hat{\eta}|^2\dk;
$$
the third term in this expression takes account of the anisotropy and the non-smooth Fourier multiplier in the
functional $\LL_2$.
The $L^2({\mathbb R}^2)$-norm of each derivative of $\tilde{\eta}_n$ thus gains a factor of
$\mu^\alpha$, and this fact allows one to obtain better estimates for the parts $\KK_k$ and $\LL_k$
of $\KK$ and $\LL$ which are homogeneous of degree $k$
and hence confirm that $\MM_\mu(\tilde{\eta}_n) = -(\mu/\LL_2(\tilde{\eta}_n))^2\LL_3(\tilde{\eta}_n) + o(\mu^3)$.

Theorem \ref{Main result for J}(ii) is established in Section \ref{Stability}. Its proof, which relies upon the strict-subadditivity
of $c_\mu$ and estimates for general minimising sequences derived in Section \ref{MS}, is now a straightforward application
of the concentration-compactness principle.

\section{The functional-analytic setting} \label{FA setting}

\subsection{The Neumann-Dirichlet operator} \label{ND operator}

Our first task is to find suitable function spaces for the functionals $\EE$ and $\II$ defined in
equations \eqn{Definition of E}, \eqn{Definition of I} and introduce
rigorous definitions of the Dirichlet-Neumann operator $G(\eta)$ and its inverse. Since the functional
$\JJ_\mu$ to be minimised involves $G(\eta)^{-1}$  (see equation \eqn{Definition of J})
we begin with the formal definition of this \emph{Neumann-Dirichlet operator} $N(\eta)$: for fixed
$\xi=\xi(x,z)$ solve the boundary-value problem 
\begin{eqnarray}
& & \parbox{6cm}{$\phi_{xx}+\phi_{yy}+\phi_{zz}=0,$}0<y<1+\eta, \label{BC for NDO 1} \\
& & \parbox{6cm}{$\phi_y - \eta_x\phi_x - \eta_z\phi_z = \xi,$}y=1+\eta, \label{BC for NDO 2}\\
& & \parbox{6cm}{$\phi_y =0,$}y=0 \label{BC for NDO 3}
\end{eqnarray}
and define
$$N(\eta)\xi=\phi|_{y=1+\eta}.$$
We study this boundary-value problem by transforming it to an equivalent problem in a
fixed domain (cf.\ Nicholls \& Reitich \cite{NichollsReitich01a}). The change of variable
$$y^\prime = \frac{y}{1+\eta}, \qquad u(x,y^\prime,z)=\phi(x,y,z)$$
transforms the variable domain $\{0<y<1+\eta(x,z)\}$ into the
slab $\Sigma=\{(x,y^\prime,z) \in {\mathbb R}
\times (0,1) \times {\mathbb R}\}$ and the boundary-value problem \eqn{BC for NDO 1}--\eqn{BC for NDO 3} into
\begin{eqnarray}
& & \parbox{90mm}{$u_{xx}+u_{yy}+u_{zz} = \partial_x F_1 + \partial_z F_2 + \partial_y F_3,$}0 < y <1,
\label{BC for u 1} \\
& & \parbox{90mm}{$u_y = F_3 + \xi,$}y=1, \label{BC for u 2} \\
& & \parbox{90mm}{$u_y=0,$}y=0, \label{BC for u 3}
\end{eqnarray}
where
\begin{eqnarray*}
F_1 & = & -\eta u_x + y \eta_x u_y, \\
F_2 & = & -\eta u_z + y \eta_z u_y, \\
F_3 & = & \frac{\eta u_y}{1+\eta} + y\eta_xu_x + y\eta_zu_z - \frac{y^2}{1+\eta}\eta_x^2u_y
- \frac{y^2}{1+\eta}\eta_z^2u_y
\end{eqnarray*}
and we have again
dropped the primes for notational simplicity; the Neumann-Dirichlet operator is given by
$$N(\eta)\xi = u|_{y=1}.$$

The next step is to develop a convenient theory for weak solutions of the boundary-value
problem \eqn{BC for u 1}--\eqn{BC for u 3}. The observation that solutions of this problem
are unique only up to additive constants leads us to introduce the completion
$H_\star^1(\Sigma)$ of
$$\SS(\Sigma,{\mathbb R})=\{u \in C^\infty(\bar{\Sigma}):
|(x,z)|^m|\partial_x^{\alpha_1}\partial_z^{\alpha_2}u|\mbox{ is bounded for all }m,\alpha_1,\alpha_2 \in {\mathbb N}_0 \}$$
with respect to the norm
$$\|u\|_\star^2:= \int_\Sigma (u_x^2+u_y^2+u_z^2)\dy\dx\dz$$
as an appropriate function space for $u$. The corresponding space for the trace $u|_{y=1}$ is the completion
$H_\star^{1/2}({\mathbb R}^2)$ of the inner product space $X_\star^{1/2}({\mathbb R}^2)$
constructed by equipping the Schwartz class $\SS({\mathbb R}^2,{\mathbb R})$ with the norm
$$\|\eta\|_{\star,1/2}^2:= \int_{{\mathbb R}^2} (1+|k|^2)^{-\frac{1}{2}}|k|^2|\hat{\eta}|^2\dk,$$
where $\hat{\eta}$ denotes the Fourier transform of $\eta$;
its dual $(H_\star^{1/2}({\mathbb R}^2))^\prime$ is the space
$$(X_\star^{1/2}({\mathbb R}^2))^\prime=\Big\{u \in \SS^\prime({\mathbb R}^2,{\mathbb R}):
\sup \{|(u,\eta)|: \eta \in X_\star^{1/2}({\mathbb R}^2),\, \|\eta\|_{\star,1/2}<1\}<\infty\Big\},$$
where $\SS^\prime({\mathbb R}^2,{\mathbb R})$ is the class of two-dimensional, real-valued,
tempered distributions. A more convenient description of $(H_\star^{1/2}({\mathbb R}^2))^\prime$ is
however available.

\begin{proposition}
Let $H_\star^{-1/2}({\mathbb R}^2)$ be the completion  of the inner product space
$X_\star^{-1/2}({\mathbb R}^2)$ constructed by equipping
$\bar{\SS}({\mathbb R}^2,{\mathbb R})$ with the norm
$$\|\eta\|_{\star,-1/2}^2:= \int_{{\mathbb R}^2} (1+|k|^2)^{\frac{1}{2}}|k|^{-2}|\hat{\eta}|^2\dk,$$
where $\bar{\SS}({\mathbb R}^2,{\mathbb R})$ is the subclass of
$\SS({\mathbb R}^2,{\mathbb R})$ consisting
of functions with zero mean. The space $H_\star^{-1/2}({\mathbb R}^2)$
 can be identified with $(H_\star^{1/2}({\mathbb R}^2))^\prime$.
\end{proposition}
{\bf Proof.} In the usual manner we identify $u \in X_\star^{-1/2}({\mathbb R}^2)$ with the distribution
$$(u,\eta)=\int_{{\mathbb R}^2} u\eta \dx\dz,$$
which belongs to $(H_\star^{1/2}({\mathbb R}^2))^\prime$ and satisfies
$\|u\|_{(H_\star^{1/2}({\mathbb R}^2))^\prime} = \|u\|_{\star,-1/2}$; it follows that
$H_\star^{-1/2}({\mathbb R}^2)$ is a subspace of $(H_\star^{1/2}({\mathbb R}^2))^\prime$.

We now demonstrate that $\eta=0$ is the only function $\eta \in X_\star^{1/2}({\mathbb R}^2)$
with the property that $(u,\eta)=0$ for all $u \in X_\star^{-1/2}({\mathbb R}^2)$; this fact implies that
$X_\star^{-1/2}({\mathbb R}^2)$ is dense in $(X_\star^{1/2}({\mathbb R}^2))^\prime=
(H_\star^{1/2}({\mathbb R}^2))^\prime$ and yields the
required result. To this end, we note that the stated property
of $\eta$ asserts in particular that $(\eta_0,\eta)=0$, where $\eta_0 \in X_\star^{-1/2}({\mathbb R}^2)$ is given
by the formula $\hat{\eta}_0 = (1+|k|^2)^{-\frac{1}{2}}|k|^2\hat{\eta}$, and the only solution of the equation
$$0=(\eta_0,\eta) = \int_{{\mathbb R}^2} (1+|k|^2)^{-\frac{1}{2}}|k|^2|\hat{\eta}|^2\dk$$
is indeed $\eta=0$. (The Fourier transform maps $\bar{\SS}({\mathbb R}^2,{\mathbb R})$
bijectively onto the subclass
$$\SS_0({\mathbb R}^2,{\mathbb C})=\{\eta \in \SS({\mathbb R}^2,{\mathbb C}):
\eta(0)=0,\ \eta(-k)=\overline{\eta (k)}\mbox{ for all }k\in {\mathbb R}^2\}$$
of $\SS({\mathbb R}^2,{\mathbb C})$.)
\qed

The following proposition, which is proved by elementary estimates, confirms the required
relationship between $H_\star^1(\Sigma)$ and $H_\star^{1/2}({\mathbb R}^2)$.
\begin{proposition}
The trace map $u \mapsto u|_{y=1}$ defines a continuous operator $H_\star^1(\Sigma) \rightarrow
H_\star^{1/2}({\mathbb R}^2)$ and has a continuous right inverse $H_\star^{1/2}({\mathbb R}^2)
\rightarrow H_\star^1(\Sigma)$.
\end{proposition}

Finally, let us take $\eta \in B_{1/2}(0) \subset W^{1,\infty}({\mathbb R}^2)$, so that the estimates
\begin{equation}
\|F_j\|_0 \leq c\|\eta\|_{1,\infty}^2 \|u\|_\star, \qquad j=1,2,3
\label{L2 estimate for Fj}
\end{equation}
imply $F_1$, $F_2$, $F_3 \in L^2(\Sigma)$. It is then a straightforward matter to
define and prove the existence of a unique weak solution to \eqn{BC for u 1}--\eqn{BC for u 3}.

\begin{definition} \label{Weak soln}
Suppose that
$\xi \in H_\star^{-1/2}({\mathbb R}^2)$ and $\eta \in B_M(0) \subset W^{1,\infty}({\mathbb R}^2)$.
A \underline{weak solution} of \eqn{BC for u 1}--\eqn{BC for u 3} is a
function $u \in H_\star^1(\Sigma)$ which satisfies
\begin{eqnarray*}
\lefteqn{\int_\Sigma (u_xw_x+u_yw_y+u_zw_z)\dx\dy\dz} \quad\\
& & = \int_\Sigma (F_1 w_x + F_2 w_z + F_3 w_y)\dx\dy\dz + \int_{{\mathbb R}^2} \xi w|_{y=1} \dx\dz
\end{eqnarray*}
for all $w \in H_\star^1(\Sigma)$.
\end{definition}

\begin{lemma}
For each $\xi \in H_\star^{-1/2}({\mathbb R}^2)$ and $\eta \in B_{1/2}(0) \subset W^{1,\infty}({\mathbb R}^2)$
the boundary-value problem \eqn{BC for u 1}--\eqn{BC for u 3} has a unique weak solution
$u \in H_\star^1(\Sigma)$.
\end{lemma}
{\bf Proof.} The existence of a unique weak solution $u \in H_\star^1(\Sigma)$ of \eqn{BC for u 1}--\eqn{BC for u 3} follows from the estimates \eqn{L2 estimate for Fj},
$$\int_{{\mathbb R}^2} \xi u|_{y=1}\dx\dz\ \leq\ \|\xi\|_{\star,-1/2}\|u|_{y=1}\|_{\star,1/2}\ \leq\ c\|\xi\|_{\star,-1/2}\|u\|_\star$$
and the Lax-Milgram lemma.\qed

We conclude with a rigorous definition of the Neumann-Dirichlet operator.

\begin{definition}
The \underline{Neumann-Dirichlet operator} for the boundary-value problem
\eqn{BC for u 1}--\eqn{BC for u 3}
is the bounded
linear operator $N(\eta): H_\star^{-1/2}({\mathbb R}^2) \rightarrow H_\star^{1/2}({\mathbb R}^2)$
defined by
$$N(\eta)\xi = u|_{y=1},$$
where $u \in H_\star^1(\Sigma)$ is the unique weak solution of \eqn{BC for u 1}--\eqn{BC for u 3}.
\end{definition}

\begin{remark} \label{Quadratic form}
Observe that
\begin{eqnarray*}
\lefteqn{\int_{{\mathbb R}^2} \xi N(\eta) \xi \dx\dz} \qquad\\
& = & \int_{{\mathbb R}^2} (\phi_y - \eta_x\phi_x - \eta_z \phi_z) \phi|_{y=\eta} \dx\dz \\
& = & \int_{\{y=1+\eta\}}  \phi \frac{\partial \phi}{\partial n} \dx\dz \\
& = & \int_{\{0<y<1+\eta\}} \Big( \phi_x^2+\phi_y^2 + \phi_z^2\Big) \dx\dy\dz\\
& = & \int_\Sigma \left(\!\left(u_x - \frac{y \eta_x u_y}{1+\eta}\right)^{\!\!2}
+\frac{u_y^2}{(1+\eta)^2}
+\left(u_z - \frac{y \eta_z u_y}{1+\eta}\right)^{\!\!2}\right)(1+\eta)\dx\dy\dz.
\end{eqnarray*}
\end{remark}

\subsection{Analyticity of the Neumann-Dirichlet operator} \label{Analyticity of NDO}

In this section we establish that $N(\eta)$ is an analytic function of $\eta$ in the above function
spaces and examine some consequences of this fact. Let us begin by recording the definition
of analyticity given by Buffoni \& Toland \cite[Definition 4.3.1]{BuffoniToland} together with
a useful fact concerning multiplication of multilinear operators.

\begin{definition} \label{BT defn of analyticity}
Let $X$ and $Y$ be Banach spaces,  $U$ be a non-empty, open subset of $X$ and
$\LL_\mathrm{s}^k(X,Y)$ be the space of bounded, $k$-linear
symmetric operators $X^k \rightarrow Y$ with norm
$$\nn m\nn:=\inf \{c: \|m(\{f\}^{(k)})\|_Y \leq c \|f\|_X^k \mbox{ \rm{for all} } f \in X\}.$$
A function $F:U \rightarrow Y$ is \underline{analytic} at a point $x_0 \in U$ if there exist
real numbers $\delta, r>0$ and a sequence $\{m_k\}$, where $m_k \in \LL_\mathrm{s}^k(X,Y)$,
$k=0,1,2,\ldots$, with the properties that
$$F(x) = \sum_{k=0}^\infty m_k(\{x-x_0\}^{(k)}), \qquad x \in B_\delta(x_0)$$
and
$$\sup_{k \geq 0} r^k\nn m_k\nn < \infty.$$
\end{definition}
\begin{remark} \label{Multiply m's}
Let $X$, $Y_1$ and $Y_2$ be Banach spaces.
Suppose that $m_1 \in \LL_\mathrm{s}^{k_1}(X,Y_1)$, $m_2 \in \LL_\mathrm{s}^{k_2}(X,Y_1)$ and that
the operation of pointwise multiplication defines
a bounded bilinear operator $Y_1 \times Y_2 \rightarrow Y$.
There exists a unique $m_3 \in \LL_\mathrm{s}^{k_1+k_2}(X,Y)$ with the
property that
$$m_1(\{f\}^{(k_1)})m_2(\{f\}^{(k_2)})=m_3(\{f\}^{(k_1+k_2)}).$$
\end{remark}

Our first task is to establish the following theorem.

\begin{theorem} \label{NDO is analytic}
The mapping $W^{1,\infty}({\mathbb R}^2) \rightarrow
\LL(H_\star^{-1/2}({\mathbb R}^2), H_\star^{1/2}({\mathbb R}^2))$ given by
$\eta \mapsto (\xi \mapsto u|_{y=1})$, where $u \in H_\star^1(\Sigma)$
is the unique weak solution of \eqn{BC for u 1}--\eqn{BC for u 3},
is analytic at the origin.
\end{theorem}

We prove Theorem \ref{NDO is analytic} using a modification of a method due to Nicholls
\& Reitich \cite{NichollsReitich01a}. Let us seek a solution of \eqn{BC for u 1}--\eqn{BC for u 3}
of the form
\begin{equation}
u(x,y,z) = \sum_{n=0}^\infty u^n(x,y,z), \label{Formal expansion}
\end{equation}
where $u^n$ is a function of $\eta$ and $\xi$ which is homogeneous of degree $n$ in $\eta$
and linear in $\xi$. Substituting this \emph{Ansatz} into the equations, one finds that
\begin{eqnarray}
& & \parbox{60mm}{$u_{xx}^0+u_{yy}^0+u_{zz}^0 = 0,$}0 < y <1, \label{BC for u0 1}\\
& & \parbox{60mm}{$u_y^0 = \xi,$}y=1, \label{BC for u0 2}\\
& & \parbox{60mm}{$u_y^0=0,$}y=0 \label{BC for u0 3}
\end{eqnarray}
and
\begin{eqnarray}
& & \parbox{90mm}{$u_{xx}^n+u_{yy}^n+u_{zz}^n = \partial_x F_1^n + \partial_z F_2^n + \partial_y F_3^n,$}0 < y <1, \label{BC for un 1} \\
& & \parbox{90mm}{$u_y^n = F_3^n,$}y=1, \label{BC for un 2} \\
& & \parbox{90mm}{$u_y^n=0,$}y=0 \label{BC for un 3}
\end{eqnarray}
for $n=1,2,3,\ldots$, where
\begin{eqnarray}
F_1^n & = & -\eta u_x^{n-1} +y \eta_x u_y^{n-1}, \label{Defn of F1n} \\
F_2^n & = & -\eta u_z^{n-1} +y \eta_z u_y^{n-1}, \label{Defn of F2n} \\
F_3^n & = & \eta \sum_{\ell=0}^{n-1} (-\eta)^\ell u_y^{n-1-\ell} + y \eta_x u_x^{n-1}
 + y \eta_z u_z^{n-1}-y^2(\eta_x^2 +\eta_z^2)\sum_{\ell=0}^{n-2}(-\eta)^\ell u_y^{n-2-\ell}.
 \label{Defn of F3n}
\end{eqnarray}

\begin{definition} \hspace{1cm}
\begin{list}{(\roman{count})}{\usecounter{count}}
\item
Suppose that $\xi \in H_\star^{-1/2}({\mathbb R}^2)$.
A \underline{weak solution} of \eqn{BC for u0 1}--\eqn{BC for u0 3} is a
function $u^0 \in H_\star^1(\Sigma)$ which satisfies
\begin{equation}
\int_\Sigma (u_x^0w_x+u_y^0w_y+u_z^0w_z)\dx\dy\dz= \int_S \xi w|_{y=1} \dx\dz
\label{u0 is a WS}
\end{equation}
for all $w \in H_\star^1(\Sigma)$. (The existence and uniqueness of a weak solution $u^0$,
which satisfies the estimate
\begin{equation}
\|u^0\|_\star \leq C_1 \|\xi\|_{\star,-1/2}, \label{Estimate for u0}
\end{equation}
follows from the Lax-Milgram lemma.)
\item Suppose that $F_1^n$, $F_2^n$, $F_3^n \in L^2(\Sigma)$.
A \underline{weak solution} of \eqn{BC for un 1}--\eqn{BC for un 3} is a
function $u^n \in H_\star^1(\Sigma)$ which satisfies
\begin{equation}
\int_\Sigma (u_x^nw_x+u_y^nw_y+u_z^nw_z)\dx\dy\dz = \int_\Sigma (F_1^n w_x + F_2^n w_z + F_3^n w_y)\dx\dy\dz
\label{un is a WS}
\end{equation}
for all $w \in H_\star^1(\Sigma)$. (The existence and uniqueness of a weak solution $u^n$,
which satisfies the estimate
\begin{equation}
\|u^n\|_\star \leq C_2 ( \|F_1^n\|_0 + \|F_2^n\|_0 + \|F_3^n\|_0),  \label{Estimate for un}
\end{equation}
follows from the Lax-Milgram lemma.)
\end{list}
\end{definition}

The next step is to compute the weak solutions $u^0$ and $u^n$, $n=1,2,\ldots$ of
the boundary-value problems \eqn{BC for u0 1}--\eqn{BC for u0 3} and
\eqn{BC for un 1}--\eqn{BC for un 3} inductively. In this fashion we obtain a sequence
$\{m^n\}_{n=0}^\infty$ of $n$-linear symmetric operators such that
$$u^n=m^n(\{\eta\}^{(n)}).$$

\begin{lemma} \label{Inductive lemma}
Suppose there exist $m^k  \in
\LL_\mathrm{s}^k(W^{1,\infty}({\mathbb R}^2), H_\star^1(\Sigma))$
and constants $C_1>0$, $B_1>2$ such that
$$u^k = m^k(\{\eta\}^{(k)}), \qquad \nn m_j^k\nn \leq C_1B_1^k\|\xi\|_{\star,-1/2}$$
for $k=0, \ldots, n-1$.

There exist $\tilde{m}_1^n$, $\tilde{m}_2^n$, $\tilde{m}_3^n \in \LL_\mathrm{s}^n(W^{1,\infty}({\mathbb R}^2),L^2(\Sigma))$
and a constant $C_3>0$ such that
$$F_j^n =\tilde{m}_j^n(\{\eta\}^{(n)}),\quad \nn \tilde{m}_j^n\nn \leq C_1C_3B_1^{n-1}\|\xi\|_{\star,-1/2}, \qquad j=1,2,3.$$
\end{lemma}
{\bf Proof.} The existence of $\tilde{m}_1^n$, $\tilde{m}_2^n$, $\tilde{m}_3^n$ follows from
the formulae \eqn{Defn of F1n}--\eqn{Defn of F3n} defining $F_1^n$, $F_2^n$, $F_3^n$
and Remark \ref{Multiply m's}.

Observe that
$$
\|F_1^n\|_0\ \leq\ \|\eta\|_{1,\infty}(\|u_x^{n-1}\|_0+\|u_y^{n-1}\|_0)\ \leq\ 2C_1B_1^{n-1}\|\xi\|_{\star,-1/2} \|\eta\|_{1,\infty}^n,
$$
$$
\|F_2^n\|_0\ \leq\ \|\eta\|_{1,\infty}(\|u_z^{n-1}\|_0+\|u_y^{n-1}\|_0)\ \leq\ 2C_1B_1^{n-1}\|\xi\|_{\star,-1/2} \|\eta\|_{1,\infty}^n.
$$
A similar calculation yields
\begin{eqnarray*}
\|F_3^n\|_0 & \leq & 2 C_1B_1^{n-1}\|\xi\|_{\star,-1/2} \|\eta\|_{1,\infty}^n \\
& & \qquad\mbox{}+ C_1\|\xi\|_{\star,-1/2}\|\eta\|_{1,\infty}^n \sum_{\ell=0}^{n-1} B_1^{n-1-\ell}
+2C_1\|\xi\|_{\star,-1/2} \|\eta\|_{1,\infty}^n \sum_{\ell=0}^{n-2} B_1^{n-2-\ell},
\end{eqnarray*}
and estimating
$$\sum_{\ell=0}^{n-1} B_1^{n-1-\ell}\ =\ B_1^{n-1}\sum_{\ell=0}^{n-1} B_1^{-\ell}
\ <\ B_1^{n-1}\sum_{\ell=0}^{n-1}\left(\frac{1}{2}\right)^{\!\!\ell}
\ <\ B_1^{n-1}\sum_{\ell=0}^\infty\left(\frac{1}{2}\right)^{\!\!\ell}
\ =\ 2B_1^{n-1},$$
$$\sum_{\ell=0}^{n-2} B_1^{n-2-\ell}\ <\ 2B_1^{n-2}\ <\ B_1^{n-1},$$
we find that
$$
\|F_3^n\|_0 \leq 6C_1B_1^{n-1}\|\xi\|_{\star,-1/2}\|\eta\|_{1,\infty}^n.\eqno{\Box}$$

\begin{theorem} \label{Analyticity theorem 1}
There exist $m^n \in \LL_\mathrm{s}^n(W^{1,\infty}({\mathbb R}^2),H_\star^1(\Sigma))$
and constants $C_1>0$, $B>1$ with the properties that
$$u^n = m^n(\{\eta\}^{(n)}), \qquad \nn m^n\nn \leq C_1B_1^n \|\xi\|_{\star,-1/2}$$
for $n=0,1,2,\ldots$.
\end{theorem}
{\bf Proof.} This result is obtained by mathematical induction. The base step ($n=0$) follows
from estimate \eqn{Estimate for u0}. Suppose the result is true for all $k < n$. The existence of
$m^n$ follows from Lemma \ref{Inductive lemma} and the fact that
$u^n \in H_\star^1(\Sigma)$ is a bounded linear function
of $F_1^n$, $F_2^n$, $F_3^n \in L^2(\Sigma)$;
using the estimate \eqn{Estimate for un}, one finds that
$$\|u\|_\star\ \leq\ 3C_2C_1C_3B_1^{n-1}\|\xi\|_{\star,-1/2}  \|\eta\|_{1,\infty}^n
\ \leq\ C_1 B_1^n\|\xi\|_{\star,-1/2} \|\eta\|_{1,\infty}^n,
$$
upon choosing $B_1>3C_2C_3$.\qed

\begin{corollary}
The mapping $W^{1,\infty}({\mathbb R}^2) \rightarrow
\LL(H_\star^{-1/2}({\mathbb R}^2), H_\star^1(\Sigma))$ given by
$\eta \mapsto (\xi \mapsto u)$, where $u \in H_\star^1(\Sigma)$
is the unique weak solution of \eqn{BC for u 1}--\eqn{BC for u 3},
is analytic at the origin.
\end{corollary}
{\bf Proof.} According to Definition \ref{BT defn of analyticity} the mapping
$W^{1,\infty}({\mathbb R}^2) \rightarrow
\LL(H_\star^{-1/2}({\mathbb R}^2), H_\star^1(\Sigma))$ given by
$\eta \mapsto (\xi \mapsto u)$, where $u$ is defined by  \eqn{Formal expansion},
\eqn{BC for u0 1}--\eqn{BC for u0 3} and \eqn{BC for un 1}--\eqn{BC for un 3}, is analytic at the origin.

It remains to verify that $u$ is a weak solution of \eqn{BC for u 1}--\eqn{BC for u 3}. The facts that
$$\sum_{n=0}^N u_x^n \rightarrow u_x, \quad
\sum_{n=0}^N u_y^n \rightarrow u_y, \quad
\sum_{n=0}^N u_z^n \rightarrow u_z$$
in $L^2(\Sigma)$ and hence that
$$\sum_{n=1}^N F_j^n \rightarrow F_j, \quad j=1,2,3$$
in $L^2(\Sigma)$ imply that
\begin{eqnarray*}
\lefteqn{\int_\Sigma \left\{ \left(\textstyle\sum_{n=0}^N u_x^n \right)w_x
+\left(\textstyle\sum_{n=0}^N u_y^n \right)w_y
+\left(\textstyle\sum_{n=0}^N u_z^n \right)w_z\right\}\dx\dy\dz} \qquad\\
& & \rightarrow \int_\Sigma (u_xw_x+u_yw_y+u_zw_2)\dx\dy\dz\hspace{2.5cm}
\end{eqnarray*}
and
\begin{eqnarray*}
\lefteqn{\int_\Sigma \left\{ \left(\textstyle\sum_{n=1}^N F_1^n \right)w_x
+\left(\textstyle\sum_{n=1}^N F_2^n \right)w_y
+\left(\textstyle\sum_{n=1}^N F_3^n \right)w_z\right\}\dx\dy\dz} \qquad\\
& & \rightarrow \int_\Sigma (F_1w_x+F_2w_y+F_3w_2)\dx\dy\dz\hspace{2.5cm}
\end{eqnarray*}
for each $w \in H_\star^1(\Sigma)$. It follows from these results, equations \eqn{u0 is a WS},
\eqn{un is a WS} and the uniqueness of limits that $u$ is a weak solution of \eqn{BC for u 1}--\eqn{BC for u 3}.\qed

Theorem \ref{NDO is analytic} is a direct consequence of the above corollary. Using this result
and the continuity of the trace operator
$H_\star^1(\Sigma) \rightarrow H_\star^{1/2}({\mathbb R}^2)$, we find that the
Neumann-Dirichlet operator
$W^{1,\infty}({\mathbb R}^2) \rightarrow
\LL(H_\star^{-1/2}({\mathbb R}^2), H_\star^{1/2}({\mathbb R}^2))$ given by
$\eta \mapsto (\xi \mapsto u|_{y=1})$ is analytic at the origin; its series representation is given by
$$
N(\eta)=\sum_{n=0}^\infty N^n(\eta),
$$
where $N^n(\eta)\xi = u^n|_{y=1}$.

The next step is to show that the Neumann-Dirichlet operator is invertible and that its inverse
is also an analytic function of $\eta$ at the origin in $W^{1,\infty}({\mathbb R}^2)$.

\begin{proposition}
The operator $N(0): H_\star^{-1/2}({\mathbb R}^2) \rightarrow H_\star^{1/2}({\mathbb R}^2)$
is an isomorphism and the norm
$$\xi \mapsto \left(\int_{{\mathbb R}^2} \xi N(0) \xi \dx\dz\right)^{\!\!\frac{1}{2}}$$
is equivalent to the usual norm on $H_\star^{-1/2}({\mathbb R}^2)$.
\end{proposition}
{\bf Proof.} Observe that $N(0)$ admits the representation
$$N(0)\xi = \FF^{-1}\left[\frac{\coth |k|}{|k|} \hat{\xi}\right]$$
as a Fourier-multiplier operator. Using the estimate
$$
c (1+|k|^2)^{\frac{1}{2}}\ \leq\ |k|\coth|k|\ \leq\ c (1+|k|^2)^{\frac{1}{2}}
$$
one finds that
\begin{equation}
c \|\xi\|_{\star,-1/2}^2\ \leq\ \int_{{\mathbb R}^2} \xi N(0) \xi \dx\dz \leq c \|\xi\|_{\star,-1/2}^2,
\label{Equivalence of norms at 0}
\end{equation}
which establishes the second assertion; the first assertion follows from the first inequality in the above estimate.\qed

\begin{corollary} \label{Equivalence of norms}
The estimate
$$c\|\xi\|_{\star,-1/2}^2\ \leq\ \int_{{\mathbb R}^2} \xi N(\eta) \xi \dx\dz \leq c\|\xi\|_{\star,-1/2}^2$$
holds for each $\eta \in B_M(0) \subset W^{1,\infty}({\mathbb R}^2)$.
In particular, the operator $N(\eta): H_\star^{-1/2}({\mathbb R}^2) \rightarrow H_\star^{1/2}({\mathbb R}^2)$
is an isomorphism and the norm
$$\xi \mapsto \left(\int_{{\mathbb R}^2} \xi N(\eta) \xi \dx\dz\right)^{\!\!\frac{1}{2}}$$
is equivalent to the usual norm on $H_\star^{-1/2}({\mathbb R}^2)$.

\end{corollary}
{\bf Proof.}  It follows from the analyticity of $N(\cdot): W^{1,\infty}({\mathbb R}^2) \rightarrow
\LL(H_\star^{-1/2}({\mathbb R}^2), H_\star^{1/2}({\mathbb R}^2))$ at the origin that
$$\|N(\eta)-N(0)\|_{\LL(H_\star^{-1/2}({\mathbb R}^2), H_\star^{1/2}({\mathbb R}^2))}
\ \leq\ c\|\eta\|_{1,\infty}\ \leq\ cM$$
and hence that
$$\left| \int_{{\mathbb R}^2} \xi N(\eta) \xi \dx\dz - \int_{{\mathbb R}^2} \xi N(0) \xi \dx\dz \right|
\leq cM\|\xi\|_{\star,-1/2}^2.$$
The result is obtained by choosing $M$ sufficiently small and combining the above estimate with
\eqn{Equivalence of norms at 0}.\qed

\begin{definition}
The \underline{Dirichlet-Neumann operator} for the boundary-value problem
\eqn{BC for u 1}--\eqn{BC for u 3}
is the bounded linear operator $G(\eta): H_\star^{1/2}({\mathbb R}^2) \rightarrow H_\star^{-1/2}({\mathbb R}^2)$
defined by $G(\eta)\Phi = N(\eta)^{-1}\Phi$.
\end{definition}

\begin{lemma} \label{Properties of DNO}
The Dirichlet-Neumann operator $G(\cdot): W^{1,\infty}({\mathbb R}^2) \rightarrow
\LL(H_\star^{1/2}({\mathbb R}^2), H_\star^{-1/2}({\mathbb R}^2))$ is analytic at the origin
and the estimate
$$c\|\Phi\|_{\star,1/2}^2\ \leq\ \int_{{\mathbb R}^2} \Phi G(\eta) \Phi \dx\dz \leq c\|\Phi\|_{\star,1/2}^2$$
holds for each $\eta \in B_M(0) \subset W^{1,\infty}({\mathbb R}^2)$.
\end{lemma}
{\bf Proof.} Define $F_1: \LL(H_\star^{1/2}({\mathbb R}^2), H_\star^{-1/2}({\mathbb R}^2)) \times
B_M(0) \rightarrow \LL(H_\star^{1/2}({\mathbb R}^2), H_\star^{1/2}({\mathbb R}^2))$ and
$F_2: \LL(H_\star^{1/2}({\mathbb R}^2), H_\star^{-1/2}({\mathbb R}^2)) \times
B_M(0) \rightarrow \LL(H_\star^{-1/2}({\mathbb R}^2), H_\star^{-1/2}({\mathbb R}^2))$
by the formulae
$$F_1(A,\eta)=N(\eta)A-I_1, \qquad F_2(A,\eta)=AN(\eta)-I_2,$$
where $I_1$ and $I_2$ are the identity operators on respectively $H_\star^{1/2}({\mathbb R}^2)$
and $H_\star^{-1/2}({\mathbb R}^2)$. It follows from the implicit function theorem that the equations
$$F_1(A,\eta)=0, \qquad F_2(A,\eta)=0$$
have unique solutions $A_1=A_1(\eta)$, $A_2=A_2(\eta)$ which are analytic at the origin;
by uniqueness we deduce that $A_1(\eta)=A_2(\eta)=G(\eta)$.

The inequality is obtained by writing $\xi = G(\eta)\Phi$ in the inequality given in Corollary
\ref{Equivalence of norms}.\qed

\subsection{The operator $K$}

Observe that the formula \eqn{Definition of L} defining $\LL$ may be written as
\begin{equation}
\LL(\eta) = \frac{1}{2}\int_{{\mathbb R}^2} \eta K(\eta) \eta \dx\dz,
\label{Definition of LL via K}
\end{equation}
where
$$K(\eta) = -\partial_x (N(\eta)\partial_x),$$
and we now study this operator in detail. Our first result
is obtained from the material presented in Section \ref{Analyticity of NDO}.
\begin{corollary} \label{K is W1infinity analytic}
The operator $K(\cdot): W^{1,\infty}({\mathbb R}^2) \rightarrow
\LL(H^{1/2}({\mathbb R}^2), H^{-1/2}({\mathbb R}^2))$ is analytic at the origin.
\end{corollary}
{\bf Proof.} This result follows from the definition of $K$
and the continuity of the operators $\partial_x: H^{1/2}({\mathbb R}^2) \rightarrow H_\star^{-1/2}({\mathbb R}^2)$ and $\partial_x: H_\star^{1/2}({\mathbb R}^2) \rightarrow H^{-1/2}({\mathbb R}^2)$.\qed

In the remainder of this section we establish the following result concerning the
analyticity of $K$ in the Sobolev spaces
$$H^r({\mathbb R}^2) = \left\{\eta \in (\SS({\mathbb R}^2,{\mathbb R}))^\prime: \|\eta\|_r^2:= \int_{{\mathbb R}^2} (1+|k|^2)^r |\hat{\eta}(k)|^2 \dk < \infty \right\}.$$

\begin{theorem} \label{K is analytic}
Suppose that $s>1$. The operator $K(\cdot): H^{s+3/2}({\mathbb R}^2) \rightarrow
\LL(H^{s+1}({\mathbb R}^2), H^s({\mathbb R}^2))$ is analytic at the origin.
\end{theorem}

The first step in the proof of this theorem is to establish additional regularity of the
weak solutions $u^0$ and $u^n$, $n=1,2,\ldots$ of the boundary-value problems
\eqn{BC for u0 1}--\eqn{BC for u0 3} and \eqn{BC for un 1}--\eqn{BC for un 3}.
This task is accomplished in the Propositions \ref{Regularity proposition 1} and
\ref{Regularity proposition 2} below; the proof of the latter is given in Appendix A.
We work in the function spaces $(H^r(\Sigma),\|\cdot\|_r)$; for 
$r \notin {\mathbb N}_0$ the space is defined by interpolation in the sense of Lions \& Magenes
\cite{LionsMagenes61} (see also Adams \& Fournier
\cite[\S7.57]{AdamsFournier}). 

\begin{proposition} \label{Regularity proposition 1}
For each $r \geq 0$ the weak solution to the boundary-value problem
\eqn{BC for u0 1}--\eqn{BC for u0 3} with $\xi=\zeta_x$ satisfies
$$\|u_x^0\|_r\, \|u_y^0\|_r,\, \|u_z^0\|_r \leq C_4 \|\zeta\|_{r+1/2}.$$
\end{proposition}
{\bf Proof.} The Fourier transform  of the weak solution $u^0$ of \eqn{BC for u0 1}--\eqn{BC for u0 3}
is given by
$$
\hat{u}^0 = \frac{\cosh |k|y}{|k|\sinh|k|} \hat{\xi}.
$$
Using this formula with $\xi=\zeta_x$, we find that
\begin{equation}
\hat{u}_y^0 = \frac{\i k_1\sinh|k|y}{\sinh |k|}\hat{\zeta}, \quad
\hat{u}_x^0 = \frac{-k_1^2\cosh|k|y}{|k|\sinh |k|}\hat{\zeta}, \quad
\hat{u}_z^0 = \frac{-k_1k_2\cosh|k|y}{|k|\sinh |k|}\hat{\zeta},
\label{Formulae for derivs of u0}
\end{equation}
whereby
$$\Big|\FF[\partial_x^{\alpha_1}\partial_y^{\alpha_2}\partial_z^{\alpha_3}\left\{\begin{array}{c}u_y^0 \\
u_x^0 \\ u_z^0 \end{array}\right\} \Big| \leq \frac{|k|^{\alpha_1+\alpha_2+\alpha_3+1}}{\sinh |k|}
\left\{\begin{array}{c} \sinh|k|y \\ \cosh|k|y\end{array}\right\}|\hat{\zeta}|,$$
where we have estimated $|k_1|, |k_2| \leq |k|$. It follows that
\begin{eqnarray*}
\left\|\FF[\partial_x^{\alpha_1}\partial_y^{\alpha_2}\partial_z^{\alpha_3}\left\{\begin{array}{c}u_y^0 \\
u_x^0 \\ u_z^0 \end{array}\right\} \right\|_0^2
& = & \int_\Sigma \Big|\FF[\partial_x^{\alpha_1}\partial_y^{\alpha_2}\partial_z^{\alpha_3}\left\{\begin{array}{c}u_y^0 \\
u_x^0 \\ u_z^0 \end{array}\right\} \Big|^2 \dy\dx\dz \\
& = & \int_{{\mathbb R}^2} \frac{|k|^{2(\alpha_1+\alpha_2+\alpha_3+1)}}{\sinh^2|k|} \left\{
\pm \frac{1}{2} + \frac{\sinh 2|k|}{4|k|}\right\} |\hat{\zeta}|^2 \dx\dz \\
& \leq & \int_{{\mathbb R}^2} (1+|k|^2)^{\alpha_1+\alpha_2+\alpha_3+1/2}
|\hat{\zeta}|^2 \dx\dz \\
& = & c \|\zeta\|_{\alpha_1+\alpha_2+\alpha_3+1/2}^2,
\end{eqnarray*}
in which the estimates
$$
\frac{|k|^2}{\sinh^2|k|}\left(-\frac{1}{2}+\frac{\sinh 2|k|}{4|k|}\right) \leq c|k|, \qquad
\frac{|k|^2}{\sinh^2|k|}\left(\frac{1}{2}+\frac{\sinh 2|k|}{4|k|}\right) \leq c(1+|k|^2)^{\frac{1}{2}}$$
have been used.

The above calculations show that
$$\left\|\left\{\begin{array}{c} u_y^0 \\ u_x^0 \\ u_z^0 \end{array}\right\}\right\|_n^2
\ =\!\!\!\!\sum_{0 \leq \alpha_1+\alpha_2+\alpha_3 \leq n} 
\left\|\partial_x^{\alpha_1}\partial_y^{\alpha_2}\partial_z^{\alpha_3}\left\{\begin{array}{c} u_y^0 \\ u_x^0 \\ u_z^0 \end{array}\right\}\right\|_n^2
\ \leq\ c\!\!\!\!\sum_{0 \leq \alpha_1+\alpha_2+\alpha_3 \leq n} \|\zeta\|_{n+1/2}^2$$
for $n \in {\mathbb N}_0$, and it follows by interpolation that
$$\left\|\left\{\begin{array}{c} u_y^0 \\ u_x^0 \\ u_z^0 \end{array}\right\}\right\|_r
\leq c\|\zeta\|_{r+1/2}, \qquad r \geq 0.\eqno{\Box}$$

\begin{proposition} \label{Regularity proposition 2}
Suppose that $F_1^n$, $F_2^n$, $F_3^n \in H^r(\Sigma)$ for $r \geq 0$.
The weak solution to the boundary-value problem \eqn{BC for un 1}--\eqn{BC for un 3}
satisfies
$$\|u^n_x\|_r,\,\|u^n_y\|_r,\,\|u^n_z\|_r \leq C_5(\|F_1^n\|_r + \|F_2^n\|_r + \|F_3^n\|_r).$$
\end{proposition}

\begin{lemma} \label{Inductive step}
Suppose that $s>1$ and there exist $m_1^k$, $m_2^k$, $m_3^k \in
\LL_\mathrm{s}^k(H^{s+3/2}({\mathbb R}^2),H^{s+1/2}(\Sigma))$
and constants $C_4>0$, $B_2>2$ such that
$$u_x^k = m_1^k(\{\eta\}^{(k)}), \quad u_y^k = m_2^k(\{\eta\}^{(k)}), \quad u_z^k = m_3^k(\{\eta\}^{(k)})$$
and
$$\nn m_j^k\nn \leq C_4 B_2^k \|\zeta\|_{s+1}, \qquad j=1,2,3$$
for $k=0, \ldots, n-1$.

There exist $m_1^n$, $m_2^n$, $m_3^n \in
\LL_\mathrm{s}^n(H^{s+3/2}({\mathbb R}^2),H^{s+1/2}(\Sigma))$
and a constant $C_6>0$ such that
$$F_j^n = m_j^n(\{\eta\}^{(n)}),\quad \nn m_j^n\nn \leq C_4C_6B_2^{n-1}\|\zeta\|_{s+1}, \qquad j=1,2,3.$$
\end{lemma}
{\bf Proof.} This result is proved inductively in the same way as Theorem \ref{Analyticity theorem 1}. The base step follows from Proposition \ref{Regularity proposition 1}, while the inductive step is treated according to the strategy of 
Lemma \ref{Inductive lemma}. Using the inequalities
$$\|w_1w_2\|_{s+1/2} \leq c_s\|w_1\|_{s+1/2}\|w_2\|_{s+1/2}, \qquad s>1$$
and
$$
\|w\|_{H^r(\Sigma)} \leq \|w\|_{H^r({\mathbb R}^2)}, \qquad w=w(x,z), \qquad r \geq0,
$$
one finds that
\begin{eqnarray*}
\|F_1^n\|_{s+1/2}
& \leq & c_s\|\eta\|_{s+1/2}\|u_x^{n-1}\|_{s+1/2} + c_s^2\|y\|_{s+1/2} \|\eta_x\|_{s+1/2}\|u_y^{n-1}\|_{s+1/2} \\
& \leq & C_4\|\zeta\|_{s+1}(c_s+c_s^2\|y\|_{s+1/2})B_2^{n-1}\|\eta\|_{s+3/2}^n
\end{eqnarray*}
and similarly
\begin{eqnarray*}
\|F_2^n\|_{s+1/2} & \leq & C_4\|\zeta\|_{s+1}(c_s+c_s^2\|y\|_{s+1/2})B_2^{n-1}\|\eta\|_{s+3/2}^n, \\
\|F_3^n\|_{s+1/2} & \leq & 2 C_4c_s^2\|y\|_{s+1/2}\|\zeta\|_{s+1} B_2^{n-1}\|\eta\|_{s+3/2}^n
+ C_4c_s\|\zeta\|_{s+1}\|\eta\|_{s+3/2}^n \sum_{\ell=0}^{n-1} c_s^\ell B_2^{n-1-\ell} \\[-2mm]
& & \qquad\mbox{} +2C_4c_s^4\|y\|_{s+1/2}^2\|\zeta\|_{s+1} \|\eta\|_{s+3/2}^n \sum_{\ell=0}^{n-2} c_s^\ell B_2^{n-2-\ell};
\end{eqnarray*}
the inductive step is completed by choosing
$B_2>\max\{2c_s,6C_5(c_s+c_s^2\|y\|_{s+1/2}+c_s^4\|y\|_{s+1/2}^2)\}$
and using Proposition \ref{Regularity proposition 2}.\qed

\begin{corollary}
For each $s>1$ the mappings $H^{s+3/2}({\mathbb R}^2) \rightarrow
\LL(H^{s+1}({\mathbb R}^2), H^{s+1/2}(\Sigma))$ given by
$\eta \mapsto (\zeta \mapsto u_x)$, $\eta \mapsto (\zeta \mapsto u_y)$
and $\eta \mapsto (\zeta \mapsto u_z)$, where $u$
is the weak solution of \eqn{BC for u 1}--\eqn{BC for u 3} with $\xi=\zeta_x$,
are analytic at the origin.
\end{corollary}

Theorem \ref{K is analytic} follows from the above corollary, the definition
$K: \eta \mapsto (\zeta \mapsto -u_x|_{y=1})$
and the continuity of the trace operator $H^{s+1/2}(\Sigma) \rightarrow H^s({\mathbb R}^2)$ for $s > 1$.
We write 
$$
K(\eta)=\sum_{n=0}^\infty K^n(\eta),
$$
where $K^n(\eta)\zeta = -u_x^n|_{y=1}$ and observe that
$$K^0 \zeta = \FF^{-1}\left[\frac{k_1^2}{|k|^2} |k|\coth |k| \hat{\zeta}\right]$$
(see \eqn{Formulae for derivs of u0}).

\begin{remark} \label{Other analytic operators}
A straightforward modification of the above analysis yields analogous
results (with the same function spaces) for the operators
$$L(\eta):=-\partial_z(N(\eta)\partial_x), \qquad M(\eta):=-\partial_z(N(\eta)\partial_z);$$
we write
$$
L(\eta)=\sum_{n=0}^\infty L^n(\eta), \qquad M(\eta)=\sum_{n=0}^\infty M^n(\eta)
$$
and note that
$$L^0 \zeta = \FF^{-1}\left[\frac{k_1k_2}{|k|^2} |k|\coth |k| \hat{\zeta}\right], \qquad
M^0 \zeta = \FF^{-1}\left[\frac{k_2^2}{|k|^2} |k|\coth |k| \hat{\zeta}\right].$$
\end{remark}

\subsection{The functionals $\KK$, $\LL$ and $\JJ_\mu$} \label{Properties of J, K, L}

The following lemma formally states the analyticity property of $\KK$ (examine
the explicit formula for $\KK$) and $\LL$ (see Theorem \ref{K is analytic}).
In particular this result implies that $\KK, \LL$ belong to the class $C^\infty(U,{\mathbb R})$ and that
equation \eqn{Definition of J} defines an operator $\JJ_\mu \in C^\infty(U\sm\{0\},{\mathbb R})$, where
$U=B_M(0) \subset H^3({\mathbb R}^2)$ and $M$ is chosen sufficiently small.

\begin{lemma} \label{KK and LL are analytic}
Equations \eqn{Definition of K}, \eqn{Definition of L} define functionals
$\KK: H^{s+1}({\mathbb R}^2) \rightarrow {\mathbb R}$,
$\LL: H^{s+3/2}({\mathbb R}^2) \rightarrow {\mathbb R}$ for $s>1$
which are analytic at the origin and satisfy $\KK(0)=\LL(0)=0$.
\end{lemma}

The following results state further useful properties of the operators $\KK$ and $\LL$.

\begin{proposition} \label{Quadratic estimates}
The functionals $\KK$ and $\LL$ satisfy
$$\KK(\eta) \geq c\|\eta\|_1^2, \qquad \LL(\eta) \geq c\|\eta_x\|_{\star,-1/2}^2$$
for each $\eta \in U$.
\end{proposition}
{\bf Proof.} Using the estimates
$$\eta_x^2, \eta_z^2\ \leq\ \|\eta\|_{1,\infty}^2\ \leq\ c\|\eta\|_3^2\ \leq\ cM^2, \qquad \eta \in U$$
and choosing $M$ small enough so that $(1+\eta_x^2+\eta_z^2)^{1/2} \leq 2$, we find that
\begin{eqnarray*}
\KK(\eta) & = & \int_{{\mathbb R}^2} \left\{\frac{\beta(\eta_x^2+\eta_z^2)}{1+(1+\eta_x^2
+\eta_z^2)^{1/2}}+\frac{\eta^2}{2}\right\}\dx\dz \\
& \geq & \frac{1}{3} \int_{{\mathbb R}^2} \left\{\beta \eta_x^2+\beta \eta_z^2 + \eta^2 \right\}\dx\dz \\
& \geq & c \|\eta\|_1^2,
\end{eqnarray*}
and furthermore
$$
\LL(\eta)\ = \ \frac{1}{2}\int_{{\mathbb R}^2} \eta_x G(\eta)^{-1} \eta_x \dx\dz\ \geq\ c\|\eta_x\|_{\star,-1/2}^2, \qquad \eta \in U
$$
(see Corollary \ref{Equivalence of norms}).\qed

\begin{lemma} \label{Formulae for the gradients}
The gradients $\KK^\prime(\eta)$ and $\LL^\prime(\eta)$ in $L^2({\mathbb R}^2)$
exist for each $\eta \in U$ and are given by the formulae
\begin{eqnarray*}
\KK^\prime(\eta) & = & -\left(\frac{\beta\eta_x}{\sqrt{1+\eta_x^2+\eta_z^2}}\right)_{\!\!x}
-\left(\frac{\beta\eta_z}{\sqrt{1+\eta_x^2+\eta_z^2}}\right)_{\!\!z}+\eta, \\
\LL^\prime(\eta) & = & -\frac{1}{2}(u_x^2+u_z^2)+\frac{u_y^2}{2(1+\eta)^2}(\eta_x^2+\eta_z^2) + \frac{u_y^2}{2(1+\eta)^2}\Bigg|_{y=1}+K(\eta)\eta,
\end{eqnarray*}
where $u$ is the weak solution of the boundary-value problem \eqn{BC for u 1}--\eqn{BC for u 3} with $\xi=\eta_x$. These
formulae define functions $\KK^\prime: H^3({\mathbb R}^2) \rightarrow H^1({\mathbb R}^2)$,
$\LL^\prime: H^{s+3/2}({\mathbb R}^2) \rightarrow H^s({\mathbb R}^2)$ for $s>1$ which are analytic at the
origin and satisfy $\KK^\prime(0)=\LL^\prime(0)=0$.
\end{lemma}
{\bf Proof.}  Differentiating the formulae \eqn{Definition of K} and
$$\LL(\eta)=\frac{1}{2}\int_\Sigma \left(\!\left(u_x - \frac{y \eta_x u_y}{1+\eta}\right)^{\!\!2}
+\frac{u_y^2}{(1+\eta)^2}
+\left(u_z - \frac{y \eta_z u_y}{1+\eta}\right)^{\!\!2}\right)(1+\eta)\dx\dy\dz$$
(see Remark \ref{Quadratic form}), one finds that
\begin{equation}
\mathrm{d}\KK[\eta](\omega)=\int_{{\mathbb R}^2} \left\{
\frac{\beta \eta_x \omega_x}{\sqrt{1+\eta_x^2+\eta_z^2}}+
\frac{\beta \eta_z \omega_z}{\sqrt{1+\eta_x^2+\eta_z^2}}+\eta\omega\right\}\dx\dz \label{Basic formula for dK}
\end{equation}
and
\begin{eqnarray}
\lefteqn{\mathrm{d}\LL[\eta](\omega)} \nonumber \\
& & = \int_\Sigma \Bigg\{ (1+\eta)(w_xu_x+w_zu_z) - y \eta_x w_x u_y -  y\eta_xu_x w_y -y\eta_zw_zu_y-y\eta_zu_zw_y \nonumber \\
& & \qquad\quad\mbox{}+ \frac{y^2u_yw_y}{1+\eta}(\eta_x^2+\eta_z^2)
+\frac{u_yw_y}{1+\eta} +\frac{\omega}{2}(u_x^2+u_z^2) - y\omega_xu_xu_y -y\omega_z u_zu_y \nonumber \\
& & \qquad\quad\mbox{}
+ \frac{y^2u_y^2}{1+\eta}(\eta_x\omega_x +\eta_z\omega_z)- \frac{y^2u_y^2}{2(1+\eta)^2}(\eta_x^2+\eta_z^2)\omega- \frac{\omega u_y^2}{2(1+\eta)^2} \Bigg\}\dy\dx\dz, \label{Basic formula for dL}
\end{eqnarray}
where $w=\mathrm{d}u[\eta](\omega)$; the expression for $\KK^\prime(\eta)$ follows directly from \eqn{Basic formula for dK}
and the expression for $\LL^\prime(\eta)$ is obtained by eliminating $w$ from \eqn{Basic formula for dL} using
the following argument.

Recall that $u$ satisfies
\begin{eqnarray*}
\lefteqn{\int_\Sigma \Bigg\{ (1+\eta)(u_xv_x+u_zv_z) - y\eta_xv_x u_y - y\eta_xu_xv_y
-y\eta_zv_z u_y - y\eta_zu_zv_y} \\
& & \quad\mbox{}+\frac{y^2u_yv_y}{1+\eta}(\eta_x^2+\eta_z^2)
+\frac{u_yv_y}{1+\eta}\Bigg\}\dy\dx\dz\hspace{2in}\\
\lefteqn{=  \int_{{\mathbb R}^2} \eta_xv|_{y=1}\dx\dz}
\end{eqnarray*}
for every $v \in H_\star^1(\Sigma)$ (Definition \ref{Weak soln} with $\xi=\eta_x$).
Differentiating this equation with respect to $\eta$,
we find that
\begin{eqnarray*}
\lefteqn{\int_\Sigma \Bigg\{ 
(1+\eta)(w_xv_x+w_zv_z)
- y\eta_xw_x v_y - y \eta_x v_x w_y - y\eta_zw_zv_y-y\eta_zv_zw_y} \\
& & \quad\mbox{}+ \frac{y^2v_yw_y}{1+\eta}(\eta_x^2+\eta_z^2)
+\frac{v_yw_y}{1+\eta}+\omega(u_xv_x+u_zv_z)- y\omega_xv_x u_y - y\omega_xu_xv_y \\
& & \quad\mbox{}-y\omega_zv_z u_y - y\omega_zu_zv_y
+2\frac{y^2u_yv_y}{1+\eta}(\eta_x\omega_x+\eta_z\omega_z) \\
& & \quad\mbox{}
-\frac{y^2u_yv_y}{(1+\eta)^2}(\eta_x^2+\eta_z^2)\omega-\frac{u_yv_y}{(1+\eta)^2}\omega
\Bigg\}\dy\dx\dz \nonumber \\
\lefteqn{= \int_{{\mathbb R}^2} \omega_xv|_{y=1}\dx\dz}
\end{eqnarray*}
for every $v \in H_\star^1(\Sigma)$; subtracting this equation with $v=u$ from
\eqn{Basic formula for dL} yields
\begin{eqnarray}
\lefteqn{\mathrm{d}\LL[\eta](\omega)} \nonumber\\
& & = \int_\Sigma \Bigg\{ -\frac{\omega}{2}(u_x^2+u_z^2) + y\omega_xu_xu_y + y\omega_zu_zu_y
-\frac{y^2u_y^2}{1+\eta}(\eta_x\omega_x +\eta_z\omega_z) \nonumber \\
& & \qquad\qquad\mbox{}
+\frac{y^2u_y^2}{2(1+\eta)^2}(\eta_x^2+\eta_z^2)\omega+ \frac{\omega u_y^2}{2(1+\eta)^2} \Bigg\}\dy\dx\dz +
\int_{{\mathbb R}^2} \omega_x u|_{y=1} \dx\dz. \label{Intermediate formula for dL}
\end{eqnarray}

Finally, observe that
\begin{eqnarray*}
\frac{1}{2}\lefteqn{\int_{{\mathbb R}^2}\Bigg\{-y\left(u_x - \frac{\eta_x y u_y}{1+\eta}\right)^{\!\!2}
-y\left(u_z - \frac{\eta_z y u_y}{1+\eta}\right)^{\!\!2} + \frac{yu_y^2}{(1+\eta)^2}
\Bigg\}\Bigg|_{y=1}\omega \dx\dz}\qquad \\
& = & \frac{1}{2}\int_\Sigma\frac{\mathrm{d}}{\mathrm{d}y}\Bigg\{
-y\omega \left(u_x - \frac{\eta_x y u_y}{1+\eta}\right)^{\!\!2}
-y\omega \left(u_z - \frac{\eta_z y u_y}{1+\eta}\right)^{\!\!2} + \frac{y\omega u_y^2}{(1+\eta)^2}
\Bigg\}\dy\dx\dz \nonumber \\
& = & \int_\Sigma \Bigg\{ -\frac{\omega}{2}(u_x^2+u_z^2) + y\omega_xu_xu_y + y\omega_zu_zu_y
-\frac{y^2u_y^2}{1+\eta}(\eta_x\omega_x +\eta_z\omega_z) \\
& & \qquad\qquad\mbox{}
+\frac{y^2u_y^2}{2(1+\eta)^2}(\eta_x^2+\eta_z^2)\omega+ \frac{\omega u_y^2}{2(1+\eta)^2} \Bigg\}\dy\dx\dz \\
& & +\int_{{\mathbb R}^2}\Bigg\{y\left(u_x - \frac{\eta_x y u_y}{1+\eta}\right)\frac{y\eta_x u_y}{1+\eta}
+y\left(u_z - \frac{\eta_z y u_y}{1+\eta}\right)\frac{y\eta_z u_y}{1+\eta}
\Bigg\}\Bigg|_{y=1}\omega\dx\dz
\end{eqnarray*}
in which the third line follows from the second by differentiating the term in braces with respect to $y$,
integrating by parts and using the fact that $u$ satisfies \eqn{BC for u 1}--\eqn{BC for u 3}. Subtracting this
formula from \eqn{Intermediate formula for dL} and multiplying out the remaining brackets yields
\begin{eqnarray*}
& & \mathrm{d}\LL[\eta](\omega) = 
\int_{{\mathbb R}^2}\Bigg\{
-\frac{1}{2}(u_x^2+u_z^2)+\frac{u_y^2}{2(1+\eta)^2}(\eta_x^2+\eta_z^2) + \frac{u_y^2}{2(1+\eta)^2}
\Bigg\}\Bigg|_{y=1}\omega \dx\dz \\
& & \hspace{1in}\mbox{}+ \underbrace{\int_{{\mathbb R}^2} \omega_x u|_{y=1} \dx\dz}_{\displaystyle = \langle K(\eta)\eta,\omega\rangle_0}.
\end{eqnarray*}

\vspace{-8.5mm}\qed

\begin{corollary} \hspace{1cm}
\begin{list}{(\roman{count})}{\usecounter{count}}
\item
The gradient $\LL^\prime(\eta)$ in $L^2({\mathbb R}^2)$ exists for each $\eta \in U$ and defines a function
$\LL^\prime: H^3({\mathbb R}^2) \rightarrow H^1({\mathbb R}^2)$ which is analytic at the
origin and satisfies $\LL^\prime(0)=0$.
\item
The gradient $\JJ_\mu^\prime(\eta)$ in $L^2({\mathbb R}^2)$ exists for each $\eta \in U$ and defines a
function $\JJ_\mu^\prime \in C^\infty(H^3({\mathbb R}^2),H^1({\mathbb R}^2))$.
\end{list}
\end{corollary}

Let us now write
$$\KK(\eta)=\KK_2(\eta)+\KK_\mathrm{nl}(\eta),
\qquad
\LL(\eta) = \LL_2(\eta) + \LL_\mathrm{nl}(\eta),$$
where
\begin{eqnarray*}
& & \KK_2(\eta) = \int_{{\mathbb R}^2}\left\{\frac{\beta}{2}\eta_x^2+\frac{\beta}{2}\eta_z^2
+\frac{\eta^2}{2}\right\}\dx\dz, \\
& & \LL_2(\eta) = \frac{1}{2}\int_{{\mathbb R}^2} \eta K^0 \eta \dx \dz = \frac{1}{2} \int_{{\mathbb R}^2}\frac{k_1^2}{|k|^2}|k|\coth|k||\hat{\eta}|^2 \dk,
\end{eqnarray*}
so that $\KK_\mathrm{nl}(\eta)$ is given by the explicit formula
\begin{equation}
\KK_\mathrm{nl}(\eta) = -\int_{{\mathbb R}^2} \frac{\beta(\eta_x^2+\eta_z^2)^2}{2(1+\sqrt{1+\eta_x^2+\eta_z^2})^2}\dx\dz,
\label{Explicit formula for KKnl}
\end{equation}
while
$$\LL_\mathrm{nl}(\eta) = \sum_{k=3}^\infty \LL_k(\eta), \qquad \LL_k(\eta) := \frac{1}{2}\int_{{\mathbb R}^2}\eta K^{k-2}(\eta)\eta\dx\dz.$$
According to Lemmata \ref{KK and LL are analytic} and \ref{Formulae for the gradients} there exist constants $B$, $C_0>0$
such that
$$|\LL_k(\eta)| \leq C_0 B^k\|\eta\|_3^k, \qquad \|\LL_k^\prime (\eta)\|_1 \leq C_0 B^{k-1} \|\eta\|_3^{k-1}$$
for $\eta \in U$ and $k=2,3,\ldots$; the following proposition gives another useful estimate for $\LL_k(\eta)$.

\begin{proposition} \label{Z estimate for LL}
The estimates
$$\|K^j(\eta)\eta\|_1 \leq C_0B^j\|\eta\|_Z^j\|\eta\|_3, \qquad j=1,2,\ldots$$
and
$$\LL_k(\eta) \leq C_0 B^{k-2} \|\eta\|_Z^{k-2} \|\eta\|_3^2, \quad \|\LL_k^\prime(\eta)\|_1 \leq C_0B^{k-3}\|\eta\|_Z^{k-3}\|\eta\|_3^2,
\qquad k=3,4,\ldots$$
hold for each $\eta \in U$, where
$$\|\eta\|_Z = \|\eta\|_{1,\infty} + \|\eta_{xx}\|_1 + \|\eta_{xz}\|_1+ \|\eta_{zz}\|_1.$$
\end{proposition}
{\bf Proof.} We establish the existence of constants $C_7>0$ and $B_2>2$
such that
$$\|u_x^n\|_2,\ \|u_y^n\|_2,\ \|u_z^n\|_2 \leq C_7B_2^n\|\eta\|_Z^n \|\zeta\|_{5/2}, \qquad n=0,1,2,3,\ldots,$$
where $u=\sum_{n=0}^\infty u^n$ is the weak solution $u$ of
\eqn{BC for u 1}--\eqn{BC for u 3} with $\xi=\zeta_x$.
Proceeding inductively, note that
\begin{eqnarray*}
\|F_1^n\|_2
& \leq & c_Z\|\eta\|_Z \|u_x^{n-1}\|_2 + c_Z \|\eta\|_Z \|y\|_2\|u_y^{n-1}\|_2 \\
& \leq & C_7c_Z\|\zeta\|_{5/2}(1+\|y\|_2)B_2^{n-1}\|\eta\|_Z^n,
\end{eqnarray*}
where we have used the elementary inequality
$$\|\eta w\|_{H^2(\Sigma)},\ \|\eta_x w\|_{H^2(\Sigma)},\ \|\eta_z w\|_{H^2(\Sigma)} \leq c_Z \|\eta\|_Z \|w\|_{H^2(\Sigma)},$$
and similarly
\begin{eqnarray*}
\|F_2^n\|_2 & \leq & C_7c_Z\|\zeta\|_{5/2}(1+\|y\|_2)B_2^{n-1}\|\eta\|_Z^n, \\
\|F_3^n\|_2 & \leq & 2C_7c_Z\|\zeta\|_{5/2}\|y\|_2 B_2^{n-1}\|\eta\|_Z^n \\
& & \qquad\mbox{} + C_7c_Z\|\zeta\|_{5/2} \|\eta\|_Z^n \sum_{\ell=0}^{n-1} c_Z^\ell B_2^{n-1-\ell}
+ 2C_7c_Z^2\|\zeta\|_{5/2}\|y\|_2^2 \|\eta\|_Z^n \sum_{\ell=0}^{n-2} c_Z^\ell B_2^{n-2-\ell}.
\end{eqnarray*}
The base step follows from Proposition \ref{Regularity proposition 1}, while the inductive step follows from the
above estimates by using Proposition \ref{Regularity proposition 2} and choosing
$B_2>\max\{2c_Z,6C_5(c_Z+c_Z\|y\|_2+c_Z^2\|y\|_2^2)\}$.
One obtains the stated estimates by setting $\zeta = \eta$ and using equation \eqn{Definition of LL via K},
Lemma \ref{Formulae for the gradients} and the fact that $K(\eta)\zeta=-u_x|_{y=1}$.\qed

A more precise description of $\LL_3^\prime(\eta)$ and $\LL_4^\prime(\eta)$ is afforded by the following semi-explicit formulae.

\begin{lemma} \label{Formulae for L3prime and L4prime}
The operators $\LL_3^\prime, \LL_4^\prime: H^3({\mathbb R}^2) \rightarrow H^1({\mathbb R}^2)$ are given by the formulae
\begin{eqnarray*}
\LL_3^\prime(\eta) & = & \frac{1}{2}\eta_x^2 - \frac{1}{2}(K^0\eta)^2 - \frac{1}{2}(L^0\eta)^2 + K^1(\eta)\eta, \\
\LL_4^\prime(\eta) & = & K^0\eta K^0(\eta K^0 \eta) + K^0\eta L^0(\eta L^0\eta)+ L^0\eta L^0(\eta K^0 \eta) + L^0\eta M_0(\eta L^0\eta)\\
& &  \quad\mbox{} +\eta\eta_{xx}K^0\eta +\eta\eta_{xz}L^0\eta + K^2(\eta)\eta.
\end{eqnarray*}
\end{lemma}
{\bf Proof.} It follows from the expression for $\LL^\prime(\eta)$ given in Lemma \ref{Formulae for the gradients} that
\begin{eqnarray*}
\LL_3^\prime(\eta) & = & -\frac{1}{2}(u_x^0)^2 -\frac{1}{2}(u_z^0)^2 +\frac{1}{2}(u_y^0)^2\Bigg|_{y=1} + K^1(\eta)\eta, \\
\LL_4^\prime(\eta) & = & -u_x^0u_x^1-u_z^0u_z^1 +u_y^0u_y^1-\eta(u_y^0)^2\Bigg|_{y=1} + K^2(\eta)\eta.
\end{eqnarray*}
The stated formulae are obtained from these equations by noting that
\begin{eqnarray*}
u_y^0|_{y=1} & = & \eta_x, \\
u_x^0|_{y=1} & = & -K^0 \eta, \\
u_z^0|_{y=1} & = & -L^0\eta, \\
\\
u_y^1|_{y=1} & = & F_3^1|_{y=1} \\
& = & \eta u_y^0 + \eta_x u_x^0 + \eta_z u_z^0 \Big|_{y=1} \\
& = & \eta \eta_x - \eta_x K^0 \eta - \eta_z L^0 \eta
\end{eqnarray*}
and that $u^1 = y\eta u_y^0+u^{1,1}+u^{1,2}$, where $u^{1,1}$, $u^{1,2}$ are the weak solutions
of the boundary-value problems
\begin{eqnarray*}
\parbox{4.5cm}{$u^{1,1}_{xx} + u^{1,1}_{yy} + u^{1,1}_{zz}=0$,}
\parbox{4.5cm}{$u^{1,2}_{xx} + u^{1,2}_{yy} + u^{1,2}_{zz}=0$,} & & 0 < y < 1, \\
\parbox{4.5cm}{$u^{1,1}_y=0$,}
\parbox{4.5cm}{$u^{1,2}_y=0$,}& & y=0, \\
\parbox{4.5cm}{$u^{1,1}_y = (\eta u_x^0)_x$,}
\parbox{4.5cm}{$u^{1,2}_y = (\eta u_z^0)_z$,}& & y=1,
\end{eqnarray*}
so that
$$u^1_x|_{y=1} = (y\eta u_y^0)_x + u^{1,1}_x + u^{1,2}_x\Big|_{y=1} = (\eta \eta_x)_x + K^0(\eta K^0 \eta) + L^0(\eta L^0 \eta),$$
$$u^1_z|_{y=1} = (y\eta u_y^0)_z + u^{1,1}_z + u^{1,2}_z\Big|_{y=1} = (\eta \eta_x)_z + L^0(\eta K^0 \eta) + M_0(\eta L^0 \eta)$$
(cf.\ Remark \ref{Other analytic operators}).\qed

\begin{corollary} \label{Explicit formulae for LL3}
The operators $\LL_3: H^3({\mathbb R}^2) \rightarrow {\mathbb R}$ and $\LL_3^\prime: H^3({\mathbb R}^2) \rightarrow H^1({\mathbb R})$ are given by the formulae
$$\LL_3(\eta)
= \frac{1}{2}\int_{{\mathbb R}^2} \Big\{ (\eta_x)^2\eta
-\eta (K^0\eta)^2 -\eta (L^0\eta)^2\Big\}\dx\dz
$$
and
$$\LL_3^\prime(\eta) = -\frac{1}{2}\eta_x^2 -\eta\eta_{xx}-\frac{1}{2}(K^0\eta)^2 -\frac{1}{2}(L^0\eta)^2 -K^0(\eta K^0\eta) -L^0(\eta L^0 \eta).$$
\end{corollary}
{\bf Proof.} The formula for $\LL_3(\eta)$ follows from Lemma \ref{Formulae for L3prime and L4prime} together with the relationships\linebreak
$\langle \LL_3^\prime(\eta),\eta \rangle_0 = 3 \LL_3(\eta)$ and $\langle K^1(\eta)\eta,\eta \rangle_0 = 2\LL_3(\eta)$, while the
formula for $\LL_3^\prime(\eta)$ is a direct consequence of that for $\LL_3(\eta)$.\qed

\begin{remark} \label{Estimates for KKnl}
The corresponding estimates $|\KK_\mathrm{nl}(\eta)|=O(\|\eta\|_3^4)$, $\|\KK_\mathrm{nl}^\prime(\eta)\|_1 = O(\|\eta\|_3^3)$
follow from Lemma \ref{Formulae for the gradients}, while the more precise estimate
$$|\KK_\mathrm{nl}(\eta)| \leq c (\|\eta_x\|_\infty + \|\eta_z\|_\infty)^2 \|\eta\|_3^2$$
is a consequence of equation \eqn{Explicit formula for KKnl}. The calculation
$$
\langle \KK_\mathrm{nl}^\prime(\eta),\eta \rangle_0 = -\int_{{\mathbb R}^2} \frac{\beta(\eta_x^2+\eta_z^2)^2}{\sqrt{1+\eta_x^2+\eta_z^2}(1+\sqrt{1+\eta_x^2+\eta_z^2})}\dx\dz
$$
and concomitant estimate
$$|\langle\KK_\mathrm{nl}^\prime(\eta),\eta \rangle_0| \leq c (\|\eta_x\|_\infty + \|\eta_z\|_\infty)^2 \|\eta\|_3^2$$
are also used in the subsequent analysis.
\end{remark}

Our final results are useful \emph{a priori} estimates. Lemma \ref{Test function 1},
whose proof is recorded in Appendix B, shows in particular that
\begin{equation}
\inf_{\eta \in U\sm\{0\}} \JJ_\mu(\eta) < 2\mu; \label{inf less than 2mu}
\end{equation}
on the other hand
\begin{equation}
\KK_2(\eta) + \frac{\mu^2}{\LL_2(\eta)} \geq 2\mu, \qquad \eta \in U\sm\{0\} \label{Lower bound on quadratic part}
\end{equation}
(because
$$\KK_2(\eta)-\LL_2(\eta) = \frac{1}{2}\int_{{\mathbb R}^2} \underbrace{\left(1+\beta|k|^2 - \frac{k_1^2}{|k|^2} |k| \coth |k|\right)}_{\displaystyle \geq 0}|\hat{\eta}|^2\dk$$
for $\beta > 1/3$, so that
$$\KK_2(\eta) + \frac{\mu^2}{\LL_2(\eta)}\ \geq\ 2\mu\sqrt{\frac{\KK_2(\eta)}{\LL_2(\eta)}}\ \geq\ 2\mu, \qquad
\eta \in U\sm\{0\}.)$$

\begin{lemma} \label{Test function 1}
There exists $\eta_{\mu}^\star \in U\sm\{0\}$ with compact support and a positive constant $c^\star$
such that
$\|\eta_\mu^\star\|_3^2 \leq c^\star \mu$ and
$\JJ_\mu(\eta_\mu^\star) < 2\mu - c\mu^3$.
\end{lemma}

Proposition \ref{Critical point estimate 1} and Corollary \ref{Critical point estimate for J}
give estimates on the size of critical points of $\JJ_\mu$ and a class of related functionals.

\begin{proposition} \label{Critical point estimate 1}
Any critical point $\eta$ of the functional $\tilde{\JJ}_{\gamma_1}: U\sm\{0\} \rightarrow {\mathbb R}$
defined by the formula
$$\tilde{\JJ}_{\gamma_1}(\eta)=\KK(\eta) - \gamma_1 \LL(\eta), \qquad \gamma_1 \in (0,4]$$
satisfies the estimate
$$\|\eta\|_3^2 \leq D\KK(\eta),$$
where $D$ is a positive constant which does not depend upon $\gamma_1$.
\end{proposition}
{\bf Proof.} Observe that
\begin{eqnarray*}
\lefteqn{\langle \KK_2^\prime(\eta), \eta \rangle_0-\langle (\KK_2^\prime(\eta))_x, \eta_{xxx}\rangle_0
-\langle (\KK_2^\prime(\eta))_z, \eta_{zzz}\rangle_0}\quad \\
& = & \int_{{\mathbb R}^2} (\beta \eta_{xxx}^2+\beta \eta_{xxz}^2+\beta\eta_{xzz}^2+\beta\eta_{zzz}^2
+\eta_{xx}^2+\eta_{zz}^2+\beta \eta_x^2+\beta\eta_z^2+\eta^2)\dx\dz \\
& \geq & \int_{{\mathbb R}^2} (\beta \eta_{xxx}^2 + \beta \eta_{zzz}^2 + \eta^2)\dx\dz
\end{eqnarray*}
and
$$\langle \KK_\mathrm{nl}^\prime(\eta), \eta \rangle_0
- \langle (\KK_\mathrm{nl}^\prime(\eta))_x, \eta_{xxx} \rangle_0
-\langle (\KK_\mathrm{nl}^\prime(\eta))_z, \eta_{zzz} \rangle_0 = O(\|\eta\|_3^4)$$
because $\|\KK_\mathrm{nl}^\prime(\eta)\|_1 \leq c \|\eta\|_3^3$ (see Remark \ref{Estimates for KKnl}).
One therefore finds that
\begin{eqnarray}
\lefteqn{\langle \KK^\prime(\eta), \eta \rangle_0
-\langle (\KK^\prime(\eta))_x, \eta_{xxx}\rangle_0
-\langle (\KK^\prime(\eta))_z, \eta_{zzz}\rangle_0}\hspace{0.5in} \nonumber \\
& \geq & \int_{{\mathbb R}^2} (\beta \eta_{xxx}^2 + \beta \eta_{zzz}^2 + \eta^2)\dx\dz
+ O(\|\eta\|_3^4) \qquad \nonumber \\
& \geq & D_1 \|\eta\|_3^2. \label{D1 estimate}
\end{eqnarray}

Choose $s \in (1,\frac{3}{2})$ and note that
$$\langle \LL^\prime(\eta), \eta \rangle_0,\,
\langle (\LL^\prime(\eta))_x, \eta_{xxx} \rangle_0,\,
\langle (\LL^\prime(\eta))_z, \eta_{zzz} \rangle_0\ \leq\ \|\LL^\prime(\eta)\|_s\|\eta\|_3
\ \leq\ c\|\eta\|_{s+3/2}\|\eta\|_3
$$
because $\LL^\prime: H^{s+3/2}({\mathbb R}^2) \rightarrow H^s({\mathbb R}^2)$ is
analytic at the origin with $\LL^\prime(0)=0$ (Lemma \ref{Formulae for the gradients});
combining this estimate and the interpolation inequality
$$\|\eta\|_{s+3/2}\ \leq\ \|\eta\|_1^{2q}\|\eta\|_3^{1-2q}
\ \leq\ c(\KK(\eta))^q\|\eta\|_3^{2-2q},$$
where $q=3/8-s/4$ (see Proposition \ref{Quadratic estimates}), one finds that
 \begin{equation}
 \langle \LL^\prime(\eta), \eta \rangle_0
-\langle (\LL^\prime(\eta))_x, \eta_{xxx}\rangle_0
-\langle (\LL^\prime(\eta))_z, \eta_{zzz}\rangle_0
\leq D_2 \KK(\eta)^q\|\eta\|_3^{2(1-q)}. \label{D2 estimate}
\end{equation}

Applying the operator
$$\langle \cdot,\eta \rangle_0 -\eta\langle (\cdot)_x,\eta_{xxx} \rangle_0-\langle (\cdot)_z,\eta_{zzz} \rangle_0$$
to the equation
$$\KK^\prime(\eta)-\gamma_1\LL^\prime(\eta)=0$$
and using the estimates \eqn{D1 estimate}, \eqn{D2 estimate} and $\gamma_1 \in (0,4]$,
we find that
$$D_1\|\eta\|_3^2 \leq 4D_2 \KK(\eta)^q\|\eta\|_3^{2(1-q)}$$
and hence that
$$\|\eta\|_3^2 \leq D \KK(\eta), \qquad D=\left(\frac{4D_2}{D_1}\right)^{\!\!\frac{1}{q}}.\eqno{\Box}$$

\begin{corollary} \label{Critical point estimate for J}
Any critical point $\eta$ of $\JJ_\mu$ with $\JJ_\mu(\eta)<2\mu$ satisfies the estimate
$$\|\eta\|_3^2 \leq D\KK(\eta).$$
\end{corollary}
{\bf Proof.} Observe that $\mu^2/\LL(\eta)\leq\JJ_\mu(\eta)$, so that $\JJ_\mu(\eta)<2\mu$
implies $\LL(\eta)>\mu/2$. Furthermore, any
critical point $\eta$ of $\JJ_\mu$ with  $\LL(\eta)>\mu/2$ is also a critical point of the
functional $\tilde{\JJ}_{\gamma_1}$, where $\gamma_1 = \mu^2/\LL(\eta)^2 \in (0,4)$.\qed

\section{Minimising sequences} \label{MS}

\subsection{The penalised minimisation problem}
In this section we study the functional 
$\JJ_{\rho,\mu}: H^3({\mathbb R}^2) \rightarrow {\mathbb R} \cup \{\infty\}$ defined by
$$\JJ_{\rho,\mu}(\eta) = \left\{\begin{array}{lll}\displaystyle \KK(\eta)+\frac{\mu^2}{\LL(\eta)}
+ \rho(\|\eta\|_3^2), & & u \in U\sm\{0\}, \\
\\
\infty, & & \eta \not\in U\sm\{0\},
\end{array}\right.$$
in which $\rho: [0,M^2) \rightarrow {\mathbb R}$ is a smooth, increasing `penalisation' function such that
$\rho(t)=0$ for $0 \leq t \leq \tilde{M}^2$ and $\rho(t) \rightarrow \infty$ as $t \uparrow M^2$. The number
$\tilde{M} \in (0,M)$ is chosen so that
$$\tilde{M}^2 > (c^\star+2D)\mu$$
(see below), and the following analysis is valid for every such choice of $\tilde{M}$, which
in particular may be chosen arbitrarily close to $M$.
In particular we give a detailed description of the qualitative properties of an arbitrary minimising sequence
$\{\eta_n\}$ for $\JJ_{\rho,\mu}$; the penalisation function ensures that $\{\eta_n\}$ does not
approach the boundary of the set $U \sm \{0\}$ in which $\JJ_\mu$ is defined. This information is used
in Section \ref{Special minimising sequence} to construct a minimising sequence $\{\tilde{\eta}_n\}$ for $\JJ_\mu$
over $U \sm \{0\}$ with $\|\tilde{\eta}_n\|_3^2 = O(\mu)$, the existence of which is a key ingredient
in the proof that the infimum $c_\mu$ of $\JJ_\mu$ over $U \sm \{0\}$ is strictly sub-additive as a function of
$\mu$ (see Section \ref{SSA}). The subadditivity property of $c_\mu$ is in turn used in Section \ref{Stability}
to establish the convergence (up to subsequences and translations) of \emph{any}
minimising sequence for $\JJ_\mu$ over $U \sm \{0\}$ which does not approach the boundary of $U$.

We begin with the following results which explain the choice of $\tilde{M}$ and
confirm that the \emph{a priori} estimates for $\JJ_\mu$
established in Section \ref{Properties of J, K, L} remain valid for $\JJ_{\rho,\mu}$.

\begin{proposition} \label{Test function 2}
The function $\eta_\mu^\star$ satisfies
$$\rho(\|\eta_\mu^\star\|_3^2)=0, \qquad \JJ_{\rho,\mu}(\eta_\mu^\star)<2\mu-c\mu^3,$$
so that
$$c_{\rho,\mu}: = \inf \JJ_{\rho,\mu} < 2\mu-c\mu^3.$$
\end{proposition}
{\bf Proof.} This result follows from the choice of $\tilde{M}$, which implies that
$\JJ_{\rho,\mu}(\eta_\mu^\star) = \JJ_\mu(\eta_\mu^\star)$, and Lemma \ref{Test function 1}.\qed

\begin{proposition} \label{Critical point estimate 2}
Any critical point $\eta$ of the functional $\tilde{\JJ}_{\gamma_1,\gamma_2}: U\sm\{0\} \rightarrow {\mathbb R}$ defined by the formula
$$\tilde{\JJ}_{\gamma_1,\gamma_2}(\eta)=
\KK(\eta) - \gamma_1 \LL(\eta)+\gamma_2 \|\eta\|_3^2, \qquad \gamma_1 \in (0,4],\ \gamma_2 \geq 0$$
satisfies the estimate
$$\|\eta\|_3^2 \leq D\KK(\eta),$$
where $D$ is a positive constant which does not depend upon $\gamma_1$ and $\gamma_2$.
\end{proposition}
{\bf Proof.} Suppose that $\gamma_2>0$ and recall that
$$\langle \tilde{\JJ}_{\gamma_1,\gamma_2}^\prime(\eta), \phi \rangle_0
\ =\ 2\gamma_2 \langle \eta, \phi \rangle_3
+ \langle \KK^\prime(\eta), \phi \rangle_0
- \gamma_1 \langle \LL^\prime(\eta),\phi \rangle_0
\ =\ 0$$
for all $\phi \in H^3({\mathbb R}^2)$. It follows from elliptic regularity theory that
$\eta \in H^6({\mathbb R}^2)$ and that the equation
$$2\gamma_2 (1-\partial_x^2-\partial_z^2)^3\eta
+\KK^\prime(\eta) -\gamma_1 \LL^\prime(\eta)=0$$
holds in $L^2({\mathbb R}^2)$. Taking the $L^2({\mathbb R}^2)$ inner product of this equation
with $\eta+\eta_{xxxx}+\eta_{zzzz}$, integrating by parts and using the estimates
\eqn{D1 estimate}, \eqn{D2 estimate}, $\gamma_1 \in (0,4]$ and
$$\int_{{\mathbb R}^2} (1-\partial_x^2-\partial_z^2)^3\eta
\, (1+\partial_x^4+\partial_z^4)\eta\dx\dz
\ = \int_{{\mathbb R}^2} (1+k_1^2+k_2^2)^3(1+k_1^4+k_2^4)|\hat{\eta}|^2 \dk\ >\ 0,$$
we find that
$$D_1\|\eta\|_3^2 \leq 4D_2 \KK(\eta)^q\|\eta\|_3^{2(1-q)}$$
and hence that $\|\eta\|_3^2 \leq D \KK(\eta)$.

The result for $\gamma_2=0$ follows directly from
Proposition \ref{Critical point estimate 1}.\qed

\begin{corollary} \label{Critical point estimate for weighted J}
Any critical point $\eta$ of $\JJ_{\rho,\mu}$ with $\JJ_{\rho,\mu}(\eta)<2\mu$ satisfies the estimates
$$\|\eta\|_3^2 \leq D\KK(\eta), \qquad  \rho(\|\eta\|_3^2)=0.$$
\end{corollary}
{\bf Proof.} Observe that $\mu^2/\LL(\eta)\leq\JJ_{\rho,\mu}(\eta)$, so that $\JJ_{\rho,\mu}(\eta)<2\mu$
implies $\LL(\eta)>\mu/2$. Furthermore, any
critical point $\eta$ of $\JJ_{\rho,\mu}$ with  $\LL(\eta)>\mu/2$ is also a critical point of the
functional $\tilde{\JJ}_{\gamma_1,\gamma_2}$, where $\gamma_1 = \mu^2/\LL(\eta)^2 \in (0,4)$
and $\gamma_2=\rho^\prime(\|\eta\|^2) \geq 0$. Using the previous proposition, one finds that
$\|\eta\|_3^2 \leq 2D\mu$ and hence $\rho(\|\eta\|_3^2)=0$ because of the choice of $\tilde{M}$.\qed

Let us now establish some basic properties of a minimising sequence $\{\eta_n\}$ for $\JJ_{\rho,\mu}$.
Without loss of generality we may assume that
$$\sup \|\eta_n\|_3 < M$$
($\|\eta_n\|_3 \rightarrow M$ would imply that $\JJ_{\rho,\mu}(\eta_n) \rightarrow \infty$),
and it follows that $\lim_{n \rightarrow \infty}
\|\eta_n\|_3$ exists and is positive ($\eta_n \rightarrow 0$ in $H^3({\mathbb R}^2)$ would also imply
that $\JJ_{\rho,\mu}(\eta_n) \rightarrow \infty$). The following lemma records further
useful properties of $\{\eta_n\}$.

\begin{lemma} \label{General properties}
Every minimising sequence $\{\eta_n\}$ for $\JJ_{\rho,\mu}$ has the properties that
$$
\JJ_{\rho,\mu}(\eta_n) < 2\mu, \quad \LL(\eta_n) > \frac{\mu}{2},
\quad \LL_2(\eta_n) \geq c\mu,
\quad \MM_{\rho,\mu}(\eta_n) \leq - c\mu^3,
\quad \|\eta_n\|_{1,\infty} \geq c \mu^3
$$
for each $n \in {\mathbb N}$, where
$$\MM_{\rho,\mu}(\eta) = \JJ_{\rho,\mu}(\eta) - \KK_2(\eta) - \frac{\mu^2}{\LL_2(\eta)}.$$
\end{lemma}
{\bf Proof.}  The first and second estimates are obtained from Proposition \ref{Test function 2} and the
elementary inequality $\mu^2/\LL(\eta_n) \leq \JJ_{\rho,\mu}(\eta_n)$, while the third and fourth are consequences of the
calculations
$$\LL_2(\eta)\ \geq\ c\|\eta_x\|_{\star,-1/2}^2\ \geq c\LL(\eta), \qquad \eta \in U$$
(see Corollary \ref{Equivalence of norms}) and
$$\MM_{\rho,\mu}(\eta_n)\ \leq\ \JJ_{\rho,\mu}(\eta_n) - 2\mu\ \leq\ -c\mu^3$$
(see inequality \eqn{Lower bound on quadratic part} and Proposition \ref{Test function 2}).

Finally, it follows from the calculation
$$|\MM_{\rho,\mu}(\eta_n)-\rho(\|\eta_n\|_3^2)|\ =\ \left|\KK_\mathrm{nl}(\eta_n) - \frac{\mu^2\LL_\mathrm{nl}(\eta_n)}{\LL(\eta_n)\LL_2(\eta_n)}\right|,$$
and the estimates
$$\MM_{\rho,\mu}(\eta_n)-\rho(\|\eta_n\|_3^2) \leq -c\mu^3, \quad
|\KK_\mathrm{nl}(\eta)|,\ |\LL_\mathrm{nl}(\eta)| \leq c\|\eta\|_{1,\infty}$$
(see Corollary \ref{K is W1infinity analytic} and Remark \ref{Estimates for KKnl}) that
$$c\mu^3\ \leq\ |\MM_{\rho,\mu}(\eta_n)-\rho(\|\eta_n\|_3^2)|\ \leq\ c\|\eta_n\|_{1,\infty}.\eqno{\Box}$$

\begin{remark} \label{General properties without rho}
Replacing $\JJ_{\rho,\mu}(\eta)$ by $\JJ_\mu(\eta)$ and
$\MM_{\rho,\mu}(\eta)$ by
$$\MM_\mu(\eta) := \JJ_\mu(\eta) - \KK_2(\eta) - \frac{\mu^2}{\LL_2(\eta)}$$
in its statement, one finds that the above lemma is also valid for a minimising sequence $\{\eta_n\}$ for $\JJ_\mu$ over $U\sm\{0\}$.
\end{remark}

\subsection{Application of the concentration-compactness principle} \label{Application of cc}

The next step is to perform a more detailed analysis of the behaviour of a minimising sequence $\{\eta_n\}$
for $\JJ_{\rho,\mu}$ by applying the concentration-compactness principle (Lions \cite{Lions84a,Lions84b}); Theorem \ref{Concentration-compactness theorem} below states this result in a form suitable for the present situation.

\begin{theorem} \label{Concentration-compactness theorem}
Any sequence $\{u_n\} \subset L^1({\mathbb R}^2)$ of non-negative functions with the
property that
$$\lim_{n \rightarrow \infty} \int_{{\mathbb R}^2} u_n(x,z) \dx\dz = \ell > 0$$
admits a subsequence for which one of the following phenomena occurs.\\
\\
\underline{Vanishing}: For each $r>0$ one has that
$$\lim_{n \rightarrow \infty}\Bigg(\sup_{(\tilde{x},\tilde{z}) \in {\mathbb R}^2} \int_{B_r(\tilde{x},\tilde{z})} \!\!\!\!\!\!\! u_n(x,z) \dx\dz\Bigg) = 0.
$$
\\
\underline{Concentration}: There is a sequence $\{(x_n,z_n)\} \subset {\mathbb R}^2$ with the property that for each
$\varepsilon>0$ there exists a positive real number $R$ with
$$\mathop{\int}_{B_R(0)} u_n(x+x_n,z+z_n) \dx\dz \geq \ell-\varepsilon$$
for each $n \in {\mathbb N}$.\\
\\
\underline{Dichotomy}: There are sequences $\{(x_n,z_n)\} \subset {\mathbb R}^2$, $\{M_n^{(1)}\}, \{M_n^{(2)}\} \subset {\mathbb R}$ and a real number
$\kappa \in (0,\ell)$ with the properties that $M_n^{(1)}$, $M_n^{(2)} \rightarrow \infty$,
$M_n^{(1)}/M_n^{(2)} \rightarrow 0$,
$$\mathop{\int}_{B_{\!M_n^{(1)}}\!(0)} \!\!\!\!u_n(x+x_n,z+z_n) \dx\dz \rightarrow \kappa,$$
$$\mathop{\int}_{B_{\!M_n^{(2)}}\!(0)} \!\!\!\!u_n(x+x_n,z+z_n) \dx\dz \rightarrow \kappa,$$
as $n \rightarrow \infty$. Furthermore
$$\lim_{n \rightarrow \infty} \Bigg(
\sup_{(\tilde{x},\tilde{z}) \in {\mathbb R}^2} \int_{B_r(\tilde{x},\tilde{z})} \!\!\!\!\!\!\! u_n(x,z) \dx\dz\Bigg) \leq \kappa
$$
for each $r>0$, and for each $\varepsilon>0$ there is a positive,
real number $R$ such that
$$\mathop{\int}_{B_R(0)} u_n(x+x_n,z+z_n) \dx\dz \geq \kappa-\varepsilon$$
for each $n \in {\mathbb N}$.
\end{theorem}

Standard interpolation inequalities show that the norms
$\|\!\cdot\!\|_r$ are metrically equivalent on $U$ for $r \in [0,3)$; we therefore study
the convergence properties of $\{\eta_n\}$ in $H^r({\mathbb R}^2)$ for $r \in [0,3)$ by focussing
on the concrete case $r=2$. One
may assume that $\|\eta_n\|_2 \rightarrow \ell$ as $n \rightarrow \infty$, where $\ell>0$
because $\eta_n \rightarrow 0$ in $H^r({\mathbb R}^2)$ for $r>5/2$ would imply that
$\JJ_{\rho,\mu}(\eta_n) \rightarrow \infty$. This observation suggests applying 
Theorem \ref{Concentration-compactness theorem} to the sequence $\{u_n\}$ defined by 
\begin{equation}
u_n = \eta_{nxx}^2 + 2\eta_{nxz}^2+\eta_{nzz}^2 + 2\eta_{nx}^2 +2\eta_{nz}^2 + \eta_n^2,
\label{Definition of un}
\end{equation}
so that $\|u_n\|_{L^1({\mathbb R}^2)} = \|\eta_n\|_2^2$.

\begin{lemma} \label{Vanishing}
The sequence $\{u_n\}$ does not have the `vanishing' property.
\end{lemma}
{\bf Proof.} This result is proved by contradiction. Suppose that $\{u_n\}$ has the
`vanishing' property. The embedding inequality
$$\|\eta_n\|_{W^{1,p}(B_1(\tilde{x},\tilde{z}))}^p \leq c
\|\eta_n\|_{H^2(B_1(\tilde{x},\tilde{z}))}^p, \qquad p>2$$
shows that
$$\int_{B_1(\tilde{x},\tilde{z})} \!\!\!\!\!\!\! (|\eta_{nx}|^p + |\eta_{nz}|^p
+|\eta_n|^p) \dx\dz \leq c
 \Bigg(\sup_{(\tilde{x},\tilde{z}) \in {\mathbb R}^2}\int_{B_1(\tilde{x},\tilde{z})} \!\!\!\!\!\!\! u_n\dx\dz
\Bigg)^{\!\!\frac{p}{2}-1}\!\!
\int_{B_1(\tilde{x},\tilde{z})} \!\!\!\!\!\!\! u_n \dx\dz$$
for each $(\tilde{x},\tilde{z}) \in {\mathbb R}^2$. Cover ${\mathbb R}^2$ by unit balls in such
a fashion that each point of ${\mathbb R}^2$ is contained in at most three balls. Summing
over all the balls, we find that
\begin{eqnarray*}
\| \eta_n \|_{1,p}^p & \leq & c
\Bigg(\sup_{(\tilde{x},\tilde{z}) \in {\mathbb R}^2}\int_{B_1(\tilde{x},\tilde{z})} \!\!\!\!\!\!\! u_n\dx\dz
\Bigg)^{\!\!\frac{p}{2}-1} \|\eta_n\|_2^2\\
& \leq & c
\Bigg(\sup_{(\tilde{x},\tilde{z}) \in {\mathbb R}^2}\int_{B_1(\tilde{x},\tilde{z})} \!\!\!\!\!\!\! u_n\dx\dz
\Bigg)^{\!\!\frac{p}{2}-1} \\
& \rightarrow & 0
\end{eqnarray*}
as $n \rightarrow \infty$ (Theorem \ref{Concentration-compactness theorem} `vanishing'),
and choosing $\delta \in (2/p,1)$, we conclude that
$$
\|\eta_n\|_{1,\infty}\ \leq\ c\|\eta_n\|_{1+\delta,p}\ \leq\ c\|\eta_n\|_{1,p}^{1-\delta}\|\eta_n\|_{2,p}^\delta
\ \leq\ c\|\eta_n\|_{1,p}^{1-\delta}\|\eta_n\|_3^\delta\ \leq\ c\|\eta_n\|_{1,p}^{1-\delta}\ \rightarrow\ 0
$$
as $n \rightarrow \infty$, which contradicts the fact that $\|\eta_n\|_{1,\infty} \geq c\mu^3$
(see Lemma \ref{General properties}).\qed

Let us now investigate the consequences of `concentration' and `dichotomy',
replacing $\{u_n\}$ by the subsequence identified by the relevant clause in
Theorem \ref{Concentration-compactness theorem} and, with a slight abuse of
notation, abbreviating the
sequences $\{u_n(\cdot+x_n,\cdot+z_n)\}$ and $\{\eta_n(\cdot+x_n,\cdot+z_n)\}$
to respectively $\{u_n\}$ and $\{\eta_n\}$. The fact that $\{\|\eta_n\|_3\}$ is
bounded implies that $\{\eta_n\}$ is weakly
convergent in $H^3({\mathbb R}^2)$; we denote its weak limit by $\eta^{(1)}$.

Lemma \ref{Concentration} deals with the `concentration' case; the following
proposition, which is proved by an argument given by Groves \& Sun \cite[p.\ 53]{GrovesSun08},
is used in its proof.

\begin{proposition} \label{Concentration implies convergence}
Suppose that $\{w_n\} \subset H^2({\mathbb R}^2)$ and $w \in H^2({\mathbb R}^2)$
have the property that for each $\tilde{\varepsilon}>0$ there exists a positive
real number $R$ with $$\|w_n\|_{H^2(|(x,z)| > R)} < \tilde{\varepsilon}$$
for every sufficiently large $n \in {\mathbb N}$ and
$w_n \rightharpoonup w$ in $H^2(|(x,z)|<R)$ as $n \rightarrow \infty$.
The sequence $\{w_n\}$ converges to $w$ in $H^r({\mathbb R}^2)$ for
each $r \in [0,2)$.
\end{proposition}
\begin{lemma} \label{Concentration}
Suppose that $\{u_n\}$ has the `concentration' property. The sequence
$\{\eta_n\}$ admits a subsequence which satisfies
$$\lim_{n \rightarrow \infty} \|\eta_n\|_3 \leq \tilde{M}$$
and converges in $H^r({\mathbb R}^2)$ for $r \in [0,3)$ to
$\eta^{(1)}$. The function $\eta^{(1)}$ satisfies the estimate
$$\|\eta^{(1)}\|_3^2\ \leq\ D\KK(\eta^{(1)})\ <\ 2D\mu,$$
minimises $\JJ_{\rho,\mu}$ and minimises $\JJ_\mu$ over
$\tilde{U}\sm\{0\}$, where $\tilde{U}=\{\eta \in H^3({\mathbb R}^2): \|\eta\|_3 < \tilde{M}\}$.
\end{lemma}
{\bf Proof.} Choose $\varepsilon>0$. The `concentration'
property asserts the existence of $R>0$ such that
$$\|\eta_n\|_{H^2(|(x,z)| > R)} < \varepsilon$$
for each $n \in {\mathbb N}$. Clearly $\eta_n \rightharpoonup \eta^{(1)}$ in $H^2(|(x,z)|<R)$,
and it follows from Proposition \ref{Concentration implies convergence}
that $\eta_n \rightarrow \eta^{(1)}$ in $H^r({\mathbb R}^2)$
for every $r \in [0,2)$ and hence for every $r \in [0,3)$.
Choosing $r > 5/2$, we find that $\KK(\eta_n) \rightarrow \KK(\eta^{(1)})$,
$\LL(\eta_n) \rightarrow \LL(\eta^{(1)})$ (see Lemma \ref{KK and LL are analytic});
furthermore
$$\rho(\|\eta^{(1)}\|_3^2) \leq \lim_{n \rightarrow\infty} \rho(\|\eta_n\|_3^2)$$
(because $\rho(\|\!\cdot\!\|_3^2)$ is weakly lower semicontinuous on $U$).

It follows that
$$\JJ_{\rho,\mu}(\eta^{(1)})\ \leq\ \lim_{n \rightarrow \infty} \JJ_{\rho,\mu}(\eta_n)\ =\ c_{\rho,\mu},$$
so that $\eta^{(1)}$ is a minimiser and hence a critical point of $\JJ_{\rho,\mu}$,
and Corollary \ref{Critical point estimate for weighted J} implies that
$$\|\eta^{(1)}\|_3^2\ \leq\ D\KK(\eta^{(1)})\ \leq\ D\JJ_{\rho,\mu}(\eta^{(1)})\ \leq\ D\JJ_{\rho,\mu}(\eta^\star_\mu)\ <\ 2D\mu,$$
so that $\JJ_\mu(\eta^{(1)}) = \JJ_{\rho,\mu}(\eta^{(1)})$. The function $\eta$ is therefore a minimiser of $\JJ_\mu$
over $\tilde{U}\sm\{0\}$, since the existence of a function $u \in \tilde{U}\sm\{0\}$ with $\JJ_\mu(u)<\JJ_\mu(\eta^{(1)})$ would lead
to the contradiction
$$\JJ_{\rho,\mu}(u)\ =\ \JJ_\mu(u)\ <\ \JJ_\mu(\eta^{(1)})\ =\ \JJ_{\rho,\mu}(\eta^{(1)})\ =\ c_{\rho,\mu}.$$

Finally, notice that
$$\lim_{n \rightarrow \infty} \JJ_{\rho,\mu}(\eta_n)\ =\ \JJ_{\rho,\mu}(\eta^{(1)})\ =\ \JJ_\mu(\eta^{(1)})\ =\ \lim_{n \rightarrow \infty}
\JJ_\mu(\eta_n),$$
whereby
$$\rho\left(\lim_{n \rightarrow \infty}\|\eta_n\|_3^2\right)=0,$$
which implies that $\lim_{n \rightarrow \infty}\|\eta_n\|_3 \leq \tilde{M}$.\qed

We now present the more involved discussion of the remaining case (`dichotomy'), using the
real number $\kappa \in (0,\ell)$ and sequences $\{M_n^{(1)}\}$, $\{M_n^{(2)}\}$ described in
Theorem \ref{Concentration-compactness theorem}.
Let $\chi:[0,\infty) \rightarrow [0,\infty)$ be a smooth, decreasing `cut-off' function such that
$$\chi(t) = \left\{\begin{array}{lll}
1, & & 0 \leq t \leq 1, \\
\\
0, & & t \geq 2,
\end{array}\right.$$
and define sequences $\{\eta_n^{(1)}\}$, $\{\eta_n^{(2)}\}$ by the formulae
$$\eta_n^{(1)}(x,z)=\eta_n(x,z)\chi\left(\frac{|(x,z)|}{M_n^{(1)}}\right), \qquad
\eta_n^{(2)}(x,z) = \eta_n(x,z) \left(1-\chi\left(\frac{|(x,z)|}{M_n^{(2)}}\right)\right),$$
so that
$$\supp \eta_n^{(1)} \subset B_{2M_n^{(1)}}(0), \qquad
\supp \eta_n^{(2)} \subset {\mathbb R}^2\sm B_{\!M_n^{(2)}}(0)$$
and the supports of $\eta_n^{(1)}$ and $\eta_n^{(2)}$ are disjoint. The following lemma,
whose proof is given in Appendix C, shows in particular how the operators
$\KK$ and $\LL$ decompose into separate parts for $\{\eta_n^{(1)}\}$ and $\{\eta_n^{(2)}\}$;
its corollary explains how the construction also induces a decomposition of $\JJ_\mu$ into the
sum of two new operators.

\begin{lemma} \label{Splitting properties 1} \hspace{1cm}
\begin{list}{(\roman{count})}{\usecounter{count}}
\item
The sequences $\{\eta_n\}$, $\{\eta_n^{(1)}\}$ and $\{\eta_n^{(2)}\}$ satisfy
$$\|\eta_n^{(1)}\|_2^2 \rightarrow \kappa, \qquad
\|\eta_n^{(2)}\|_2^2 \rightarrow \ell-\kappa, \qquad
\|\eta_n-\eta_n^{(1)}-\eta_n^{(2)}\|_2 \rightarrow 0$$
as $n \rightarrow \infty$.
\item
The sequences $\{\eta_n\}$, $\{\eta_n^{(1)}\}$ and $\{\eta_n^{(2)}\}$ satisfy the bounds
$$\sup \|\eta_n^{(1)}\|_3<M, \quad \sup \|\eta_n^{(2)}\|_3<M, \quad
\sup \|\eta_n^{(1)}+\eta_n^{(2)}\|_3 < M.$$
\item
The functional $\KK$ satisfies
$$\KK(\eta_n)-\KK(\eta_n^{(1)})-\KK(\eta_n^{(2)}) \rightarrow 0$$
$$\|\KK^\prime(\eta_n)-\KK^\prime(\eta_n^{(1)})-\KK^\prime(\eta_n^{(2)})\|_1 \rightarrow 0$$
as $n \rightarrow \infty$. The functional $\LL$ has the same properties.
\item
The limits $\lim_{n \rightarrow \infty} \LL(\eta_n^{(1)})$ and $\lim_{n \rightarrow \infty} \LL(\eta_n^{(2)})$
are positive.
\end{list}
\end{lemma}

\begin{corollary} \label{J splitting}
The sequences $\{\eta_n\}$, $\{\eta_n^{(1)}\}$ and $\{\eta_n^{(2)}\}$ satisfy
$$\lim_{n \rightarrow \infty} \JJ_\mu(\eta_n) = \lim_{n \rightarrow \infty} \JJ_{\mu^{(1)}}(\eta_n^{(1)})+
\lim_{n \rightarrow \infty} \JJ_{\mu^{(2)}}(\eta_n^{(2)}),$$
$$\lim_{n \rightarrow \infty} \JJ_\mu^\prime(\eta_n) = \lim_{n \rightarrow \infty} \JJ_{\mu^{(1)}}^\prime(\eta_n^{(1)})+
\lim_{n \rightarrow \infty} \JJ_{\mu^{(2)}}^\prime(\eta_n^{(2)}),$$
where the positive numbers $\mu^{(1)}$, $\mu^{(2)}$ are defined by
$$\mu^{(1)}\ =\ \mu\, \frac{\displaystyle\lim_{n \rightarrow \infty} \LL(\eta_n^{(1)})}{\displaystyle\lim_{n \rightarrow \infty} \vphantom{\LL^2}\LL(\eta_n)},
\qquad
\mu^{(2)}\ =\ \mu\, \frac{\displaystyle\lim_{n \rightarrow \infty} \LL(\eta_n^{(2)})}{\displaystyle\lim_{n \rightarrow \infty} \vphantom{\LL^2}\LL(\eta_n)}$$
and the limits in the second equation are taken in $H^1({\mathbb R}^2)$.
\end{corollary}
{\bf Proof.} Observe that
\begin{eqnarray*}
\lim_{n \rightarrow \infty} \JJ_\mu(\eta_n)
& = & \lim_{n \rightarrow \infty} \left\{\KK(\eta_n)+\frac{\mu^2}{\LL(\eta_n)}\right\} \\
& = & \lim_{n \rightarrow \infty} \KK(\eta_n)
+\frac{\mu^2}{\Big(\displaystyle\lim_{n \rightarrow \infty} \LL(\eta_n)\Big)^{\!\!2}}
\lim_{n \rightarrow \infty} \LL(\eta_n) \\
& = & \lim_{n \rightarrow \infty} \left\{\KK(\eta_n^{(1)})+\KK(\eta_n^{(2)})\right\}
+\frac{\mu^2}{\Big(\displaystyle\lim_{n \rightarrow \infty} \LL(\eta_n)\Big)^{\!\!2}}\lim_{n \rightarrow \infty}
\left\{\LL(\eta_n^{(1)})+\LL(\eta_n^{(2)})\right\} \\
& = & \lim_{n \rightarrow \infty} \left\{\KK(\eta_n^{(1)})+\frac{(\mu^{(1)})^2}{\LL(\eta_n^{(1)})}\right\}
+\lim_{n \rightarrow \infty} \left\{\KK(\eta_n^{(2)})+\frac{(\mu^{(2)})^2}{\LL(\eta_n^{(2)})}\right\} \\[2mm]
& = &  \lim_{n \rightarrow \infty} \JJ_{\mu^{(1)}}(\eta_n^{(1)}) + \lim_{n \rightarrow \infty} \JJ_{\mu^{(2)}}(\eta_n^{(2)});
\end{eqnarray*}
the second identity is proved in a similar fashion.\qed

The convergence properties of the sequence $\{\eta_n^{(1)}\}$ are discussed in Lemma
\ref{eta1 converges strongly} and Corollary \ref{Size of eta1} below.

\begin{lemma} \label{eta1 converges strongly}
 \hspace{1cm}
\begin{list}{(\roman{count})}{\usecounter{count}}
\item
The sequence $\{\eta_n^{(1)}\}$ converges to $\eta^{(1)}$
in $H^r({\mathbb R}^2)$ for $r \in [0,3)$.
\item
The function $\eta^{(1)}$ satisfies the estimates
$\|\eta^{(1)}\|_3^2 \leq D\KK(\eta^{(1)})$ and $\|\eta^{(1)}\|_2 \geq c\mu^6$.
\end{list}
\end{lemma}
{\bf Proof.} (i) Choose $\tilde{\varepsilon}>0$. The `dichotomy' property asserts the existence of
$R>0$ such that
$$\|\eta_n\|_{H^2(|(x,z)|<R)}^2 > \kappa - {\textstyle\frac{1}{2}}\tilde{\varepsilon}^2.$$
Taking $n$ large enough so that $M_n^{(1)}>R$, we find that
$$\|\eta_n^{(1)}\|_2^2 - \|\eta_n^{(1)}\|_{H^2(|(x,z)|>R)}^2\ =\ \|\eta_n^{(1)}\|_{H^2(|(x,z)|<R)}^2
\ =\ \|\eta_n\|_{H^2(|(x,z)|<R)}^2\ >\ \kappa-{\textstyle\frac{1}{2}}\tilde{\varepsilon}^2,$$
whereby
$$\|\eta_n^{(1)}\|_{H^2(|(x,z)|>R)}^2\ <\ \|\eta_n^{(1)}\|_2^2
-(\kappa-{\textstyle\frac{1}{2}}\tilde{\varepsilon}^2)\ <\ \tilde{\varepsilon}^2$$
for sufficiently large $n$, since $\|\eta_n^{(1)}\|_3^2 \rightarrow \kappa$ as $n \rightarrow \infty$.
Clearly $\eta_n \rightharpoonup \eta^{(1)}$ in $H^2(|(x,z)|<R)$ as $n \rightarrow \infty$, and this
fact implies that $\eta_n^{(1)} \rightharpoonup \eta^{(1)}$ in $H^2(|(x,z)|<R)$ because
$\eta_n(x,z)=\eta_n^{(1)}(x,z)$ for $(x,z) \in B_R(0)$. The first assertion
now follows from Proposition \ref{Concentration implies convergence}.

(ii) Note that the second derivative of $\JJ_{\rho,\mu}$ is bounded on every subset of $U$
on which $\JJ_{\rho,\mu}$ is bounded. It follows that the minimising
sequence $\{\eta_n\}$ for $\JJ_{\rho,\mu}$ is also a Palais-Smale sequence for this functional
(cf.\ Mawhin \& Willem \cite[Corollary 4.1]{MawhinWillem}), so that
$$\langle \JJ_{\rho,\mu}^\prime(\eta_n), \phi \rangle_0
\ =\ 2\rho^\prime(\|\eta_n\|_3^2) \langle \eta_n, \phi \rangle_3
+ \langle \KK^\prime(\eta_n), \phi \rangle_0
- \frac{\mu^2}{\LL(\eta_n)^2}\langle \LL^\prime(\eta_n),\phi \rangle_0
\ \rightarrow\ 0$$
and hence
$$2\rho^\prime(\|\eta_n\|_3^2) \langle \eta_n, \phi \rangle_3
+ \langle \KK^\prime(\eta_n^{(1)}) + \KK^\prime(\eta_n^{(2)}), \phi \rangle_0
- \frac{\mu^2}{\LL(\eta_n)^2}\langle \LL^\prime(\eta_n^{(1)})+\LL^\prime(\eta_n^{(2)}),\phi \rangle_0
\ \rightarrow\ 0$$
as $n \rightarrow \infty$ for each $\phi \in C_0^\infty({\mathbb R}^2)$ (see Lemma
\ref{Splitting properties 1}(iii)).
Choosing $R$ so that $\supp \phi \subset B_{2R}(0)$, one finds that $\langle \KK^\prime(\eta_n^{(2)}),\phi \rangle_0=0$ and the corresponding result
$$\lim_{n \rightarrow \infty} \langle \LL^\prime(\eta_n^{(2)}), \phi \rangle_0 = 0$$ 
for $\LL$ is proved in in Appendix D (Theorem \ref{Key splitting theorem}). Furthermore
$\eta_n \rightharpoonup \eta^{(1)}$ in $H^3({\mathbb R}^2)$, so that
$\langle \eta_n, \phi \rangle_0 \rightarrow \langle \eta^{(1)}, \phi \rangle_0$, and
$\eta_n^{(1)} \rightarrow \eta^{(1)}$ in $H^t({\mathbb R}^2)$ for $t>5/2$, so that
$\LL^\prime(\eta_n^{(1)}) \rightarrow \LL^\prime(\eta^{(1)})$ (see Lemma
\ref{Formulae for the gradients}) and
\begin{eqnarray*}
\langle \KK^\prime(\eta_n^{(1)}),\phi \rangle_0
& = &
\int_{{\mathbb R}^2} \left(\frac{\beta\eta_{nx}^{(1)}\phi_x}{\sqrt{1+(\eta_{nx}^{(1)})^2+(\eta_{nz}^{(1)})^2}}
+\frac{\beta\eta_{nz}^{(1)}\phi_z}{\sqrt{1+(\eta_{nx}^{(1)})^2+(\eta_{nz}^{(1)})^2}} + \eta^{(1)}\phi\right)\dx\dz \\
& \rightarrow &
\int_{{\mathbb R}^2} \left(\frac{\beta\eta_x^{(1)}\phi_x}{\sqrt{1+(\eta_x^{(1)})^2+(\eta_z^{(1)})^2}}
+\frac{\beta\eta_z^{(1)}\phi_z}{\sqrt{1+(\eta_x^{(1)})^2+(\eta_z^{(1)})^2}} + \eta^{(1)}\phi\right)\dx\dz \\
& = &
\langle \KK^\prime(\eta^{(1)}),\phi \rangle_0.
\end{eqnarray*}
We conclude that
$$2\rho^\prime\left(\lim_{n \rightarrow \infty} \|\eta_n\|_3^2\right)\langle \eta^{(1)},\phi \rangle_3
 + \langle \KK^\prime(\eta^{(1)}),\phi \rangle_0
- \frac{\mu^2}{\displaystyle \lim_{n \rightarrow \infty} \LL(\eta_n)^2}
\langle \LL^\prime(\eta^{(1)}),\phi\rangle_0 = 0$$
for each $\phi \in C_0^\infty({\mathbb R}^2)$ and hence for each $\phi \in H^3({\mathbb R}^2)$.
Because
$$\rho^\prime\left(\lim_{n \rightarrow \infty} \|\eta_n\|_3^2\right) \geq 0, \qquad
\frac{\mu^2}{\displaystyle \lim_{n \rightarrow \infty} \LL(\eta_n)^2} \leq 4$$
Proposition \ref{Critical point estimate 2} asserts that $\|\eta^{(1)}\|_3^2 \leq D\KK(\eta^{(1)})$.

Because $\|\eta_n\|_{1,\infty} \geq c\mu^3$ (see Lemma \ref{General properties})
there exists a sequence $\{(\tilde{x}_n,\tilde{z}_n)\}$ with
the property that
$$|\eta_n(\tilde{x}_n,\tilde{z}_n)| + |\eta_{nx}(\tilde{x}_n,\tilde{z}_n)| + |\eta_{nz}(\tilde{x}_n,\tilde{z}_n)|
\geq c \mu^3,$$
and using the embedding $H^{5/2}(B_1(\tilde{x}_n,\tilde{z}_n)) \subset W^{1,\infty}(B_1(\tilde{x}_n,\tilde{z}_n))$,
we find that
\begin{eqnarray*}
c\mu^3 & \leq& \|\eta_n\|_{H^{5/2}(B_1(\tilde{x}_n,\tilde{z}_n))} \\
& \leq& c\|\eta_n\|_{H^2(B_1(\tilde{x}_n,\tilde{z}_n))}^{\frac{1}{2}}\|\eta_n\|_{H^3(B_1(\tilde{x}_n,\tilde{z}_n))}^{\frac{1}{2}}\\
& \leq& c\|\eta_n\|_{H^2(B_1(\tilde{x}_n,\tilde{z}_n))}^{\frac{1}{2}}\|\eta_n\|_3^{\frac{1}{2}} \\
& \leq& c\|\eta_n\|_{H^2(B_1(\tilde{x}_n,\tilde{z}_n))}^{\frac{1}{2}}.
\end{eqnarray*}
It follows that
$$\sup_{(\tilde{x},\tilde{z})\in{\mathbb R}^2} \int_{B_1(\tilde{x},\tilde{z})} u_n(x,z)\dx\dz \geq c\mu^{12}$$
and hence that
$$\|\eta^{(1)}\|_2^2\ =\ \lim_{n \rightarrow \infty} \|\eta_n^{(1)}\|_2^2\ =\ \kappa
\ \geq\ \lim_{n \rightarrow \infty}
 \Bigg(
\sup_{(\tilde{x},\tilde{z}) \in {\mathbb R}^2} \int_{B_1(\tilde{x},\tilde{z})} \!\!\!\!\!\!\! u_n(x,z) \dx\dz\Bigg)
\ \geq\ c\mu^{12}.\eqno{\Box}$$

\begin{corollary} \label{Size of eta1}
The function $\eta^{(1)}$ satisfies the estimate $\|\eta^{(1)}\|_3^2 \leq 2D\mu$.
\end{corollary}
{\bf Proof.} Using Corollary \ref{J splitting}, we find that
$$\KK(\eta^{(1)})\ \leq\ \JJ_{\mu^{(1)}}(\eta^{(1)})\ = \lim_{n \rightarrow \infty} \JJ_{\mu^{(1)}}(\eta_n^{(1)})
\ \leq\ \lim_{n \rightarrow \infty} \JJ_{\mu}(\eta_n)\ \leq
\ \lim_{n \rightarrow \infty} \JJ_{\rho,\mu}(\eta_n)\ =\ c_{\rho,\mu}\ <\ 2\mu;$$
the assertion follows from this result and the estimate $\|\eta^{(1)}\|_3^2 \leq D\KK(\eta^{(1)})$.\qed

The above results show that $\{\eta_n^{(1)}\}$
essentially `concentrates' and converges. The behaviour of $\{\eta_n^{(2)}\}$ on the other hand is analogous
to that of the original sequence $\{\eta_n\}$: it is a minimising sequence for the
functional
$\JJ_{\rho_2,\mu^{(2)}}: H^3({\mathbb R}^2) \rightarrow {\mathbb R} \cup \{\infty\}$ defined by
$$\JJ_{\rho_2,\mu^{(2)}}(\eta) = \left\{\begin{array}{lll}\displaystyle \KK(\eta)+\frac{(\mu^{(2)})^2}{\LL(\eta)}
+ \rho_2(\|\eta\|_3^2), & & \eta \in U_2\sm\{0\}, \\
\\
\infty, & & \eta \not\in U_2\sm\{0\},
\end{array}\right.$$
where
$$U_2 = \{\eta \in H^3({\mathbb R}^2): \|\eta\|_3^2 \leq M^2 - \|\eta^{(1)}\|_3^2\}, \qquad
\rho_2(\|\eta\|_3^2) =\rho(\|\eta^{(1)}\|_3^2 + \|\eta\|_3^2).$$
This fact is established in Lemma \ref{Iteration results} below; the following result is used in its proof.

\begin{proposition} \label{Splitting properties 2}
For every $\{v_n\} \subset U$ with $\|\eta^{(1)}\|_3^2+\sup \|v_n\|_3^2<M^2$
there exists an increasing, unbounded sequence
$\{S_n\}$ of positive real numbers such that
$$\lim_{n \rightarrow \infty} \|\eta^{(1)} + \tau_{S_n}v_n\|_3^2 = \|\eta^{(1)}\|_3^2 + \lim_{n \rightarrow \infty} \|v_n\|_3^2$$
and
$$\lim_{n \rightarrow \infty} \JJ_{\mu}(\eta^{(1)} + \tau_{S_n}v_n) \leq \JJ_{\mu^{(1)}}(\eta^{(1)})+
\lim_{n \rightarrow \infty} \JJ_{\mu^{(2)}}(v_n),$$
where $(\tau_X v_n)(x,z) := v_n(x+X,z)$.
\end{proposition}
{\bf Proof.} Choose $\varepsilon>0$, take $R>0$ large enough so that
$$\|\eta^{(1)}\|_{H^3(|(x,z)|>R)} < \varepsilon,$$
let $\{R_n\}$ be an increasing, unbounded sequence of positive real numbers
such that
$$\|v_n\|_{H^3(|(x,z)|>R_n)} < n^{-1}$$
and choose $S_n>2R+2R_n$, $n=1,2, \ldots$.
Defining
$$w(x,z)=\eta^{(1)}(x,z)\chi\left(\frac{|(x,z)|}{R}\right), \qquad w_n(x,z)=(\tau_{S_n}v_n)(x,z)\chi\left(\frac{|(x+S_n,z)|}{R_n}\right),$$
note that the supports of $w$ and $w_n$ are mutually disjoint ($\supp w \subset B_{2R}(0)$ while
$\supp w_n \subset B_{2R_n}(-S_n,0) \subset {\mathbb R}^2 \sm B_{2R}(0)$) and
$\|\eta^{(1)}-w\|_3 =O(\varepsilon)$, $\|\tau_{S_n}v_n-w_n\|_3=O(n^{-1})$.

The first assertion follows from the calculation
\begin{eqnarray*}
\|\eta^{(1)}+\tau_{S_n}v_n\|_3-\|\eta^{(1)}\|_3-\|v_n\|_3 & = & \|\eta^{(1)}+\tau_{S_n}v_n\|_3-\|\eta^{(1)}\|_3-\|\tau_{S_n}v_n\|_3 \\
& = & \underbrace{\|w+w_n\|_3 - \|w\|_3 - \|w_n\|_3}_{\displaystyle = 0} + O(\varepsilon)+ O(n^{-1}),
\end{eqnarray*}
and the same method yields the corresponding results for $\KK$ and $\LL$ in place of
$\|\cdot\|_3$ provided that
$\sup \|\eta^{(1)} + \tau_{S_n}v_n\|_3$,  $\sup \|w+w_n\|_3 < M$ (when dealing with $\LL$ the first
three terms on the right-hand side of the above equation are $o(1)$ rather than zero (see
Theorem \ref{Key splitting theorem})).
Clearly
$$\|\eta^{(1)}\|_3^2 + \sup \|\tau_{S_n}v_n\|_3^2 < M^2,$$
whereby
$$\|w+w_n\|_3^2\ =\ \|w\|_3^2 + \|w_n\|_3^2
\ =\ \|\eta^{(1)}\|_3^2 + \|\tau_{S_n}v_n\|_3^2 + O(\varepsilon) + O(n^{-1})$$
and
$$\|\eta^{(1)} + \tau_{S_n}v_n\|_3^2 = \|w+w_n\|_3^2 + O(\varepsilon) + O(n^{-1}).$$
Replacing $\{v_n\}$ by a subsequence if necessary, we conclude that
$\sup \|v + \tau_{S_n}v_n\|_3$ and $\sup \|w+w_n\|_3$ are indeed strictly smaller than
$M$ for sufficiently small values of $\varepsilon$.

Turning to the second assertion, observe that
\begin{eqnarray*}
\lim_{n \rightarrow \infty} \JJ_{\mu}(\eta^{(1)}+\tau_{S_n}v_n)
& = &  \lim_{n \rightarrow \infty} \KK(\eta^{(1)}+\tau_{S_n}v_n) + \frac{\mu^2}{\displaystyle
\lim_{n \rightarrow \infty} \LL(\eta^{(1)}+\tau_{S_n}v_n)} \\
& = & \KK(\eta^{(1)}) + \lim_{n \rightarrow \infty} \KK(v_n)
+ \frac{\mu^2}{\displaystyle \LL(\eta^{(1)})+\lim_{n \rightarrow \infty} \LL(v_n)} \\
& \leq & \KK(\eta^{(1)}) + \lim_{n \rightarrow \infty} \KK(v_n) + \frac{(\mu^{(1)})^2}{\LL(\eta^{(1)})}
+ \frac{(\mu^{(2)})^2}{\displaystyle \lim_{n \rightarrow \infty} \LL(v_n)} \\
& = & \JJ_{\mu^{(1)}}(\eta^{(1)}) + \lim_{n \rightarrow \infty} \JJ_{\rho_2,\mu^{(2)}} (v_n),
\end{eqnarray*}
in which the inequality
$$\frac{(\mu^{(1)}+\mu^{(2)})^2}{\ell_1+\ell_2} \leq \frac{(\mu^{(1)})^2}{\ell_1} + \frac{(\mu^{(2)})^2}{\ell_2},
\qquad \ell_1,\ell_2>0$$
has been used.\qed

\begin{lemma} \hspace{1in} \label{Iteration results}
\begin{list}{(\roman{count})}{\usecounter{count}}
\item
The sequence $\{\eta_n^{(2)}\}$ is a minimising sequence for $\JJ_{\rho_2,\mu^{(2)}}$.
\item
The sequences $\{\eta_n\}$ and $\{\eta_n^{(2)}\}$ satisfy
\begin{equation}
\lim_{n \rightarrow \infty} \rho(\|\eta_n\|_3^2) = \lim_{n \rightarrow \infty} \rho_2(\|\eta_n^{(2)}\|_3^2),
\label{rho does not change}
\end{equation}
\begin{equation}
\lim_{n \rightarrow \infty} \JJ_{\rho,\mu}(\eta_n) = \JJ_{\mu^{(1)}}(\eta^{(1)})+
\lim_{n \rightarrow \infty} \JJ_{\rho_2,\mu^{(2)}}(\eta_n^{(2)}) \label{Penalised splitting}
\end{equation}
and
$$
\|\eta^{(1)}\|_3^2 + \lim_{n \rightarrow \infty} \|\eta_n^{(2)}\|_3^2 \leq \lim_{n \rightarrow \infty} \|\eta_n\|_3^2
$$
with equality if $\lim_{n \rightarrow \infty} \rho(\|\eta_n\|_3^2) > 0$.
\end{list}
\end{lemma}
{\bf Proof.} (i) The existence of a minimising sequence $\{v_n\}$ for $\JJ_{\rho_2,\mu^{(2)}}$ with
$$\lim_{n \rightarrow \infty} \JJ_{\rho_2,\mu^{(2)}}(v_n) < \lim_{n \rightarrow \infty} \JJ_{\rho_2,\mu^{(2)}}(\eta_n^{(2)})$$
implies that $\JJ_{\rho_2,\mu^{(2)}}(v_n) \not\rightarrow \infty$, so that
$\|\eta^{(1)}\|_3^2 + \sup \|v_n\|_3^2 < M^2$. One therefore obtains the contradiction
\begin{eqnarray*}
\lim_{n \rightarrow \infty} \JJ_{\rho,\mu}(\eta^{(1)}+\tau_{S_n}v_n)
& \leq & \rho\left(\|\eta^{(1)}\|_3^2 + \lim_{n \rightarrow \infty} \|v_n\|_3^2\right)+
\JJ_{\mu^{(1)}}(\eta^{(1)}) + \lim_{n \rightarrow \infty} \JJ_{\mu^{(2)}} (v_n) \\
& = & \rho_2\left(\lim_{n \rightarrow \infty} \|v_n\|_3^2\right)+\JJ_{\mu^{(1)}}(\eta^{(1)}) + \lim_{n \rightarrow \infty} \JJ_{\mu^{(2)}} (v_n) \\
& = & \JJ_{\mu^{(1)}}(\eta^{(1)}) + \lim_{n \rightarrow \infty} \JJ_{\rho_2,\mu^{(2)}} (v_n) \\
& < & \JJ_{\mu^{(1)}}(\eta^{(1)}) + \lim_{n \rightarrow \infty} \JJ_{\rho_2,\mu^{(2)}} (\eta_n^{(2)}) \\
& = & c_{\rho,\mu},
\end{eqnarray*}
where the sequence $\{S_n\}$ is constructed in Proposition \ref{Splitting properties 2}).

(ii) Using the inequality \eqn{Sum in Y} and the facts that $\eta_n^{(1)}$, $\eta_n^{(2)}$ have disjoint support
and $\eta_n \rightharpoonup \eta^{(1)}$ in $H^3({\mathbb R}^2)$ as $n \rightarrow \infty$, one finds that
$$
\lim_{n \rightarrow \infty} \|\eta_n\|_3^2 \ \geq \ \lim_{n \rightarrow \infty} \|\eta_n^{(1)}+\eta_n^{(2)}\|_3^2
\ = \ \lim_{n \rightarrow \infty} \|\eta_n^{(1)}\|_3^2 + \lim_{n \rightarrow \infty} \|\eta_n^{(2)}\|_3^2
\ \geq \ \|\eta^{(1)}\|_3^2 + \lim_{n \rightarrow \infty} \|\eta_n^{(2)}\|_3^2.
$$
We now treat the cases $\lim_{n \rightarrow \infty} \rho(\|\eta_n\|_3^2)=0$ and
$\lim_{n \rightarrow \infty} \rho(\|\eta_n\|_3^2)>0$ separately.
\begin{itemize}
\item
The condition $\rho\left(\lim_{n \rightarrow \infty} \|\eta_n\|_3^2\right)=0$ implies that
$\lim_{n \rightarrow \infty} \|\eta_n\|_3 \leq \tilde{M}$, from which it follows that
$$\|\eta^{(1)}\|_3^2 + \lim_{n \rightarrow \infty} \|\eta_n^{(2)}\|_3^2 \leq \tilde{M}^2$$
and hence
$$\lim_{n \rightarrow \infty} \rho(\|\eta^{(1)}\|_3^2 + \|\eta_n^{(2)}\|_3^2) = 0,$$
that is $\rho_2\left(\lim_{n \rightarrow \infty} \|\eta_n^{(2)}\|_3^2\right)=0$.
\item
The condition $\rho\left(\lim_{n \rightarrow \infty} \|\eta_n\|_3^2\right)>0$ is not compatible with
the strict inequality
$$
\|\eta^{(1)}\|_3^2 + \lim_{n \rightarrow \infty} \|\eta_n^{(2)}\|_3^2 < \lim_{n \rightarrow \infty} \|\eta_n\|_3^2,
$$
which would imply that
$$\lim_{n \rightarrow \infty} \rho(\|\eta^{(1)}\|_3^2 + \|\eta_n^{(2)}\|_3^2) < \lim_{n \rightarrow \infty} \rho(\|\eta_n\|_3^2),$$
so that
\begin{eqnarray*}
\lim_{n \rightarrow \infty} \JJ_{\rho,\mu}(\eta^{(1)}+\tilde{\eta}_n^{(2)})
& \leq & \lim_{n \rightarrow \infty} \rho(\|\eta^{(1)}\|_3^2 + \|\eta_n^{(2)}\|_3^2)
+\JJ_{\mu^{(1)}}(\eta^{(1)}) + \lim_{n \rightarrow \infty} \JJ_{\mu^{(2)}}(\eta_n^{(2)}) \\
& < & \lim_{n \rightarrow \infty} \rho(\|\eta_n\|_3^2)
+\JJ_{\mu^{(1)}}(\eta^{(1)}) + \lim_{n \rightarrow \infty} \JJ_{\mu^{(2)}}(\eta_n^{(2)}) \\
& = & \lim_{n \rightarrow \infty} \rho(\|\eta_n\|_3^2) + \lim_{n \rightarrow \infty} \JJ_\mu(\eta_n) \\
& = & \lim_{n \rightarrow \infty} \JJ_{\rho,\mu}(\eta_n)
\end{eqnarray*}
(see Proposition \ref{Splitting properties 2}),
and contradict the fact that $\{\eta_n\}$ is a minimising sequence for $\JJ_{\rho,\mu}$.
We conclude that
$$
\|\eta^{(1)}\|_3^2 + \lim_{n \rightarrow \infty} \|\eta_n^{(2)}\|_3^2 = \lim_{n \rightarrow \infty} \|\eta_n\|_3^2,
$$
whereby
$$\lim_{n \rightarrow \infty} \rho_2(\|\eta_n^{(2)}\|_3^2)
\ =\ \lim_{n \rightarrow \infty} \rho(\|\eta^{(1)}\|_3^2+\|\eta_n^{(2)}\|_3^2)
\ =\ \lim_{n \rightarrow \infty} \rho(\|\eta_n\|_3^2).$$
 \end{itemize}

In both cases the limit \eqn{Penalised splitting} follows from Corollary \ref{J splitting},
equation \eqn{rho does not change} and the fact that
$\JJ_{\mu^{(1)}}(\eta_n^{(1)}) \rightarrow \JJ_{\mu^{(1)}}(\eta^{(1)})$ as $n \rightarrow \infty$.\qed

\subsection{Iteration}

The next step is to apply the concentration-compactness principle to the sequence $\{u_{2,n}\}$
given by
$$
u_{2,n} = (\partial_{xx}\eta_{2,n})^2 + 2 (\partial_x\partial_z \eta_{2,n})^2
+ (\partial_{zz}\eta_{2,n})^2 + 2(\partial_x\eta_{2,n})^2
+ 2(\partial_z\eta_{2,n})^2 + \eta_{2,n}^2,
$$
where $\eta_{2,n}=\eta_n^{(2)}$,
and repeat the above analysis.
We proceed iteratively in this fashion, writing $\{\eta_n\}$, $\mu$ and $U$ in iterative
formulae as respectively $\{\eta_{1,n}\}$, $\mu_1$ and $U_1$. The following lemma describes
the result of one step in this procedure.

\begin{lemma} \label{Results of iteration}
Suppose there exist functions $\eta^{(1)}$, \ldots, $\eta^{(k)} \in H^3({\mathbb R}^2)$
and a sequence $\{\eta_{{k+1},n}\}
\subset H^3({\mathbb R}^2)$ with the following properties.
\begin{list}{(\roman{count})}{\usecounter{count}}
\item
The sequence $\{\eta_{{k+1},n}\}$ is a minimising sequence for
$\JJ_{\rho_{k+1},\mu_{k+1}}: H^3({\mathbb R}^2) \rightarrow {\mathbb R} \cup \{\infty\}$ defined by
$$\JJ_{\rho_{k+1},\mu_{k+1}}(\eta) = \left\{\begin{array}{lll}\displaystyle \KK(\eta)+\frac{\mu_{k+1}^2}{\LL(\eta)} + \rho_{k+1}(\|\eta\|_3^2), & & \eta \in U_{k+1}\sm\{0\}, \\
\\
\infty, & & \eta \not\in U_{k+1}\sm\{0\},
\end{array}\right.$$
where
$$U_{k+1} = \left\{\eta \in H^3({\mathbb R}^2): \|\eta\|_3^2 \leq M^2 - \sum_{j=1}^k\|\eta^{(j)}\|_3^2\right\}$$
and
$$\rho_{k+1}(\|\eta\|_3^2)=\rho\left(\sum_{j=1}^k\|\eta^{(j)}\|_3^2 + \|\eta\|_3^2\right),
\qquad
\mu_{k+1} = \mu \frac{\displaystyle \lim_{n \rightarrow \infty} \LL(\eta_{k+1,n})}{\displaystyle \lim_{n \rightarrow \infty} \LL(\eta_n)}>0.
$$
\item
The functions $\eta^{(1)}$, \ldots, $\eta^{(k)}$ satisfy
\begin{equation}
0<\|\eta^{(j)}\|_3^2 \leq D\KK(\eta^{(j)}), \qquad j=1,\ldots,k \label{TF hypothesis 1}
\end{equation}
and
\begin{equation}
c_{\rho,\mu} = \sum_{j=1}^k \JJ_{\mu_j^{(1)}} (\eta^{(j)}) + c_{\rho_{k+1},\mu_{k+1}},
\label{TF hypothesis 2}
\end{equation}
where
$$\mu_j^{(1)} = \mu \frac{\LL(\eta^{(j)})}{\displaystyle \lim_{n \rightarrow \infty} \LL(\eta_n)}, \qquad
j=1,\ldots,k$$
and $c_{\rho_{k+1},\mu_{k+1}} = \inf \JJ_{\rho_{k+1},\mu_{k+1}}$.
\item
The sequences $\{\eta_n\}$, $\{\eta_{k+1,n}\}$ and functions $\eta^{(1)}$,
\ldots, $\eta^{(k)}$ satisfy
\begin{eqnarray*}
& & \sum_{j=1}^k \KK(\eta^{(j)}) + \lim_{n \rightarrow \infty} \KK(\eta_{k+1,n}) = \lim_{n \rightarrow \infty} \KK(\eta_n),\\
& & \sum_{j=1}^k \LL(\eta^{(j)}) + \lim_{n \rightarrow \infty} \LL(\eta_{k+1,n}) = \lim_{n \rightarrow \infty} \LL(\eta_n),
\end{eqnarray*}
$$
\lim_{n \rightarrow \infty} \rho(\|\eta_n\|_3^2) = \lim_{n \rightarrow \infty} \rho_{k+1}(\|\eta_{k+1,n}\|_3^2)
$$
and
$$
\sum_{j=1}^k \|\eta^{(j)}\|_3^2
+ \lim_{n \rightarrow \infty} \|\eta_{k+1,n}\|_3^2 \leq \lim_{n \rightarrow \infty} \|\eta_n\|_3^2
$$
with equality if $\lim_{n \rightarrow \infty} \rho(\|\eta_n\|_3^2) > 0$.
\end{list}

Under these hypotheses an application of the concentration-compactness principle to the
sequence
$$
u_{k+1,n} = (\partial_{xx}\eta_{k+1,n})^2 + 2 (\partial_x\partial_z \eta_{k+1,n})^2
+ (\partial_{zz}\eta_{k+1,n})^2 + 2(\partial_x\eta_{k+1,n})^2
+ 2(\partial_z\eta_{k+1,n})^2 + \eta_{k+1,n}^2
$$
yields the following results.

\begin{enumerate}
\item
The sequence $\{u_{k+1,n}\}$ does not have the `vanishing property'.
\item
Suppose that $\{u_{k+1,n}\}$ has the `concentration' property. There exists a sequence\linebreak
$\{(x_{k+1,n},z_{k+1,n})\} \subset {\mathbb R}^2$ and a subsequence of
$\{\eta_{k+1,n}(\cdot+x_{k+1,n},\cdot+z_{k+1,n})\}$ which satisfies
$$\lim_{n \rightarrow \infty} \|\eta_{k+1,n}(\cdot+x_{k+1,n},\cdot+z_{k+1,n})\|_3^2 \leq \tilde{M}^2 - \sum_{j=1}^k \|\eta^{(j)}\|_3^2$$
and converges in $H^r({\mathbb R}^2)$ for $r \in [0,3)$. The limiting function $\eta^{(k+1)}$ satisfies
$$\sum_{j=1}^{k+1} \KK(\eta^{(j)}) = \lim_{n \rightarrow \infty} \KK(\eta_n), \qquad
\sum_{j=1}^{k+1} \LL(\eta^{(j)}) = \lim_{n \rightarrow \infty} \LL(\eta_n),$$
$$\|\eta^{(k+1)}\|_3^2 \leq D\KK(\eta^{(k+1)}), \qquad c_{\rho,\mu} = \sum_{j=1}^{k+1} \JJ_{\mu_j^{(1)}}(\eta^{(j)}),$$
with $\mu_{k+1}^{(1)} =\mu_{k+1}$,
minimises $\JJ_{\rho_{k+1},\mu_{k+1}}$ and minimises $\JJ_{\mu_{k+1}^{(1)}}$ over
$\tilde{U}_{k+1}\sm\{0\}$, where
$$\tilde{U}_{k+1} = \left\{\eta \in H^3({\mathbb R}^2): \|\eta\|_3^2 \leq \tilde{M}^2 - \sum_{j=1}^k\|\eta^{(j)}\|_3^2\right\}.$$
The step concludes the iteration.
\item
Suppose that $\{u_{k+1,n}\}$ has the `dichotomy' property. There exist sequences
$\{\eta_{k+1,n}^{(1)}\}$, $\{\eta_{k+1,n}^{(2)}\}$ with the following properties.
\begin{list}{(\roman{count})}{\usecounter{count}}
\item
The sequence $\{\eta_{k+1,n}^{(1)}\}$ converges in $H^r({\mathbb R}^2)$ for $r \in [0,3)$ to a function $\eta^{(k+1)}$
which satisfies the estimates
$$\|\eta^{(k+1)}\|_3^2 \leq D\KK(\eta^{(k+1)}), \qquad \|\eta^{(k+1)}\|_2 \geq c\mu_{k+1}^6.$$
\item
The sequence $\{\eta_{k+1,n}^{(2)}\}$ is a minimising sequence for
$\JJ_{\rho_{k+2},\mu_{k+1}^{(2)}}: H^3({\mathbb R}^2) \rightarrow {\mathbb R} \cup \{\infty\}$ defined by
$$\JJ_{\rho_{k+2},\mu_{k+1}^{(2)}}(\eta) = \left\{\begin{array}{lll}\displaystyle \KK(\eta)+\frac{(\mu_{k+1}^{(2)})^2}{\LL(\eta)} + \rho_{k+2}(\|\eta\|_3^2), & & \eta \in U_{k+2}\sm\{0\}, \\
\\
\infty, & & \eta \not\in U_{k+2}\sm\{0\},
\end{array}\right.$$
where
$$U_{k+2} = \left\{\eta \in H^3({\mathbb R}^2): \|\eta\|_3^2 \leq M^2 - \sum_{j=1}^{k+1}\|\eta^{(j)}\|_3^2\right\}$$
and
$$\rho_{k+2}(\|\eta\|_3^2)=\rho\left(\sum_{j=1}^{k+1}\|\eta^{(j)}\|_3^2 + \|\eta\|_3^2\right),
\qquad
\mu_{k+1}^{(2)} = \mu \frac{\displaystyle \lim_{n \rightarrow \infty} \LL(\eta_{k+1,n}^{(2)})}{\displaystyle \lim_{n \rightarrow \infty} \LL(\eta_n)}>0;
$$
furthermore
$$c_{\rho,\mu} = \sum_{j=1}^{k+1} \JJ_{\mu_j^{(1)}} (\eta^{(j)}) + c_{\rho_{k+2},\mu_{k+1}^{(2)}},$$
where
$$\mu_{k+1}^{(1)} = \mu \frac{\LL(\eta^{(k+1)})}{\displaystyle \lim_{n \rightarrow \infty} \LL(\eta_n)}, \qquad
c_{\rho_{k+2},\mu_{k+1}^{(2)}} = \inf \JJ_{\rho_{k+2},\mu_{k+1}^{(2)}}.$$
\item
The sequences $\{\eta_n\}$, $\{\eta_{k+1,n}^{(2)}\}$ and functions $\eta^{(1)}$,
\ldots, $\eta^{(k+1)}$ satisfy
\begin{eqnarray*}
& & \sum_{j=1}^k \KK(\eta^{(j+1)}) + \lim_{n \rightarrow \infty} \KK(\eta_{k+1,n}^{(2)}) = \lim_{n \rightarrow \infty} \KK(\eta_n),\\
& & \sum_{j=1}^k \LL(\eta^{(j+1)}) + \lim_{n \rightarrow \infty} \LL(\eta_{k+1,n}^{(2)}) = \lim_{n \rightarrow \infty} \LL(\eta_n),
\end{eqnarray*}
$$
\lim_{n \rightarrow \infty} \rho(\|\eta_n\|_3^2) = \lim_{n \rightarrow \infty} \rho_{k+2}(\|\eta_{k+1,n}^{(2)}\|_3^2)
$$
and
\begin{equation}
\sum_{j=1}^{k+1} \|\eta^{(j)}\|_3^2
+ \lim_{n \rightarrow \infty} \|\eta_{k+1,n}^{(2)}\|_3^2 \leq \lim_{n \rightarrow \infty} \|\eta_n\|_3^2
\label{Uniformity of sequences}
\end{equation}
with equality if $\lim_{n \rightarrow \infty} \rho(\|\eta_n\|_3^2) > 0$.
\end{list}
\end{enumerate}
The iteration continues to the next step with $\eta_{k+2,n} = \eta_{k+1}^{(2)}$, $n \in {\mathbb N}$.
\end{lemma}
{\bf Proof.} It follows from \eqn{TF hypothesis 1} and \eqn{TF hypothesis 2} that
$$\sum_{j=1}^k \|\eta^{(j)}\|_3\  \leq\ D\sum_{j=1}^k \KK(\eta^{(j)})\ \leq\ D\sum_{j=1}^k \JJ_{\mu_j^{(1)}} (\eta^{(j)})
\ \leq\ 2D\mu$$
and therefore that
$$\|\eta_{\mu_{k+1}}^\star\|_3^2 + \sum_{j=1}^k\|\eta^{(j)}\|_3^2\ \leq\ (c^\star + 2D)\mu\ <\ \tilde{M}^2.$$
This estimate shows that $\rho_{k+1}(\|\eta_{\mu_{k+1}}^\star\|_3^2)=0$ and hence
$$c_{\rho_{k+1},\mu_{k+1}}\ \leq\ \JJ_{\rho_{k+1},\mu_{k+1}}(\eta_{\mu_{k+1}}^\star)
\ =\ \JJ_{\mu_{k+1}}(\eta^\star_{\mu_{k+1}})\ <\ 2\mu_{k+1}-c(\mu_{k+1})^3.$$
The analysis of the sequence $\{u_{k+1,n}\}$ by means of the concentration-compactness principle
is therefore the same as that given for the sequence $\{u_n\}$ in Section \ref{Application of cc}.\qed

The above construction does not assume that the iteration terminates (that is `concentration'
occurs after a finite number of iterations). If it does not terminate
we let $k \rightarrow \infty$ in Lemma \ref{Results of iteration} and find that
$\|\eta^{(k)}\|_3 \rightarrow 0$ (because $\sum_{j=1}^k \|\eta^{(j)}\|_3^2< 2D\mu$
for each $k \in {\mathbb N}$, so that the series $\sum_{j=1}^\infty \|\eta^{(j)}\|_3^2$ converges),
$\mu_k \rightarrow 0$ (because $\|\eta^{(k)}\|_3^2 \geq c\mu_k^6$), $c_{\rho_k,\mu_k} \rightarrow 0$
(because $c_{\rho_k,\mu_k} < 2\mu_k$) and
$$c_{\rho,\mu} = \sum_{j=1}^\infty \JJ_{\mu_j^{(1)}}(\eta^{(j)}).$$

For completeness we conclude our analysis of minimising sequences by recording the following corollary
of Lemma \ref{Results of iteration} which is not used in the remainder of the paper.

\begin{corollary} \label{Asymptotically unpenalised}
Every minimising sequence $\{\eta_n\}$ for $\JJ_{\rho,\mu}$ satisfies
$\lim_{n \rightarrow \infty} \|\eta_n\|_3 \leq \tilde{M}$.
\end{corollary}
{\bf Proof.} We proceed by contradiction. Apply the iterative scheme described above to
$\{\eta_n\}$ and suppose that $\lim_{n \rightarrow \infty} \|\eta_n\|_3 > \tilde{M}$,
that is $\lim_{n \rightarrow \infty} \rho(\|\eta_n\|_3^2)>0$. Notice that the iteration does not terminate and
equality holds in \eqn{Uniformity of sequences} for all $k \in {\mathbb N}$; passing to the limit
$k \rightarrow \infty$, we therefore find that
\begin{equation}
\sum_{j=1}^\infty \|\eta^{(j)}\|_3^2 + \limsup_{k \rightarrow \infty}\lim_{n \rightarrow \infty} \|\eta_{k,n}^{(2)}\|_3^2 =\lim_{n \rightarrow \infty} \|\eta_n\|_3^2. \label{Rule out 1}
\end{equation}
On the other hand the limit $\lim_{k \rightarrow \infty} c_{\rho_k,\mu_k} = 0$ implies that
$$\lim_{k \rightarrow \infty} \left( \lim_{n \rightarrow \infty} \rho_k (\|\eta_{k,n}\|_3^2) \right) = 0,$$
and hence
$$\rho \left(\sum_{j=1}^\infty \|\eta^{(j)}\|_3^2 + \limsup_{k \rightarrow \infty}\lim_{n \rightarrow \infty}
\|\eta_{k,n}^{(2)}\|_3^2\right)=0,$$
that is
\begin{equation}
\sum_{j=1}^\infty \|\eta^{(j)}\|_3^2 + \limsup_{k \rightarrow \infty}\lim_{n \rightarrow \infty} \|\eta_{k,n}^{(2)}\|_3^2
\leq \tilde{M}^2. \label{Rule out 2}
\end{equation}
Combining \eqn{Rule out 1} and \eqn{Rule out 2}, one obtains the contradiction
$\lim_{n \rightarrow \infty} \|\eta_n\|_3 \leq \tilde{M}$.\qed

\subsection{Construction of a special minimising sequence} \label{Special minimising sequence}

The goal of this section is the proof of the following theorem, the sequence advertised in
which has properties beyond those enjoyed by a general minimising sequence 
(cf.\ Remark \ref{General properties without rho}).

\begin{theorem} \label{Special MS}
There exists a minimising sequence $\{\tilde{\eta}_n\}$ for $\JJ_\mu$ over $U\sm\{0\}$
with the properties that $\|\tilde{\eta}_n\|_3^2 \leq c \mu$ for each $n \in {\mathbb N}$ and
$$\lim_{n \rightarrow \infty} \JJ_\mu(\tilde{\eta}_n) = c_\mu, \qquad
\lim_{n \rightarrow \infty} \|\JJ_\mu^\prime(\tilde{\eta}_n)\|_1=0.$$
\end{theorem}

The sequence $\{\tilde{\eta}_n\}$ is constructed by gluing together
the functions $\eta^{(j)}$ identified in Section \ref{MS}\linebreak above with increasingly large distances between them,
so that the interactions between the `tails' of the individual functions is negligible; the minimal distance is chosen
so that $\|\tilde{\eta}_n\|_3^2$ is approximately $\sum_{j=1}^m \|\eta^{(j)}\|_3^2 = O(\mu)$.
(Here, and in the remainder of the section, the index $j$ is taken between $1$ and $m$, where
$m=k$ if the iteration described in Section \ref{MS} terminates after $k$ steps and $m=\infty$
if it does not terminate.)
Because $\KK$ is a local, translation-invariant operator we find that
$$\lim_{n \rightarrow \infty} \KK(\tilde{\eta}_n) = \sum_{j=1}^m \KK(\eta^{(j)}),$$
and in fact the corresponding result for $\LL$ also holds; it
is obtained by a careful analysis of an integral-operator representation of the functions
$u^n$ defining $\LL$. We deduce from these results that
$$
\lim_{n \rightarrow \infty}\JJ_{\rho,\mu}(\tilde{\eta}_n)\ =\ \sum_{j=1}^m \JJ_{\mu_j^{(1)}}(\eta^{(j)})\ =\ c_{\rho,\mu},
$$
so that $\{\tilde{\eta}_n\}$ is a minimising sequence for $\JJ_{\rho,\mu}$,
and the fact that the construction is independent of the choice of $\tilde{M}$ allows us
to conclude that $\{\tilde{\eta}_n\}$ is also a minimising sequence for $\JJ_\mu$
over $U\sm\{0\}$. Finally, a similar argument yields
$$
\lim_{n \rightarrow \infty}\|\JJ_\mu^\prime(\tilde{\eta}_n)\|_1\ =\ \sum_{j=1}^m \JJ^\prime_{\mu_j^{(1)}}(\eta^{(j)})=0.
$$

We begin with a precise statement of the algorithm used to construct $\tilde{\eta}_n$.

\begin{enumerate}
\item
Choose $R_j>1$ large enough so that
$$\|\eta^{(j)}\|_{H^3(|(x,z)|>R_j)} < \frac{\mu}{2^j}.$$
\item
Write $S_1=0$ and choose $S_j>S_{j-1}+2R_j+2R_{j-1}$ for $j=2,\ldots,m$.
\item
The sequence $\{\tilde{\eta}_n\}$ is defined by
$$\tilde{\eta}_n = \sum_{j=1}^m \tau_{S_j+(j-1)n}\eta^{(j)}, \qquad n \in {\mathbb N}.$$
\end{enumerate}

It is confirmed in Corollary \ref{Special sequence size},
Corollary \ref{Special PS property} and Proposition \ref{Special sequence minimises}
below that $\{\tilde{\eta}_n\}$ has the properties advertised in Theorem 4.1.

\begin{proposition} \label{Uniform convergence over sequences}
There exists a constant $C>0$ such that
$$\left\|\sum_{j=1}^m \tau_{S_j} \eta^{(j)}\right\|_3^2 \leq 3C^2D\mu.$$
for all choices of $\{S_j\}_{j=1}^m$. Moreover, in the case $m=\infty$ the series converges uniformly over all such sequences. 
\end{proposition}
{\bf Proof.} Defining
$$\eta^{(j,1)}(x,z):=\eta^{(j)}(x,z)\chi\left(\frac{|(x,z)|}{R_j}\right),$$
observe that the supports of the functions $\tau_{S_j}\eta^{(j,1)}$
are disjoint ($\mathrm{supp}\,\tau_{S_j}\eta^{(j,1)} \subseteq B_{2R_j}(-S_j,0)$), and that
$\eta^{(j,1)}$ and $\eta^{(j,2)}:=\eta^{(j)}-\eta^{(j,1)}$
satisfy
$$\|\eta^{(j,1)}\|_3 \leq C \|\eta^{(j)}\|_3, \qquad \|\eta^{(j,2)}\|_3 \leq C \frac{\mu}{2^j}$$
uniformly over $j$. It follows that
\begin{eqnarray*}
\left\|\sum_{j=i}^m \tau_{S_j}\eta^{(j)}\right\|_3^2 & = &
\left\|\sum_{j=i}^m (\tau_{S_j}\eta^{(j,1)}+\tau_{S_j}\eta^{(j,2)})\right\|_3^2 \\
& \leq & \left\|\sum_{j=i}^m \tau_{S_j}\eta^{(j,1)}\right\|_3^2 + 2 \left\|\sum_{j=i}^m \tau_{S_j}\eta^{(j,1)}\right\|_3
\sum_{j=i}^m \|\tau_{S_j}\eta^{(j,2)}\|_3 + \left(\sum_{j=i}^m \|\tau_{S_j}\eta^{(j,2)}\|_3 \right)^{\!\!2} \\
& = & \sum_{j=i}^m \|\tau_{S_j}\eta^{(j,1)}\|_3^2 + 2 \left(\sum_{j=i}^m \|\tau_{S_j}\eta^{(j,1)}\|_3^2\right)^{\!\!\frac{1}{2}}
\sum_{j=i}^m \|\tau_{S_j}\eta^{(j,2)}\|_3 + \left(\sum_{j=i}^m \|\tau_{S_j}\eta^{(j,2)}\|_3 \right)^{\!\!2} \\
& = & \sum_{j=i}^m \|\eta^{(j,1)}\|_3^2 + 2 \left(\sum_{j=i}^m \|\eta^{(j,1)}\|_3^2\right)^{\!\!\frac{1}{2}}
\sum_{j=i}^m \|\eta^{(j,2)}\|_3 + \left(\sum_{j=i}^m \|\eta^{(j,2)}\|_3 \right)^{\!\!2} \\
& \leq & C^2\left(\sum_{j=i}^m \|\eta^{(j)}\|_3^2 + 2\mu \left(\sum_{j=i}^m \|\eta^{(j)}\|_3^2\right)^{\!\!\frac{1}{2}}
\sum_{j=i}^m 2^{-j} + \mu^2\left(\sum_{j=i}^m 2^{-j} \right)^{\!\!2} \right);
\end{eqnarray*}
in the case $m=\infty$ the series on the left-hand side therefore converges uniformly over all sequences $\{S_j\}_{j=1}^\infty$.
Choosing $i=1$, we find that
$$
\left\|\sum_{j=1}^m \tau_{S_j}\eta^{(j)}\right\|_3^2
\ \leq\ C^2\left(\sum_{j=1}^m \|\eta^{(j)}\|_3^2 + 2\mu \left(\sum_{j=1}^m \|\eta^{(j)}\|_3^2\right)^{\!\!\frac{1}{2}}
+ \mu^2\right)\ \leq\ 2C^2D\mu + O(\mu^\frac{3}{2}).\eqno{\Box}
$$

\begin{corollary} \label{Special sequence size}
The sequence $\{\tilde{\eta}_n\}$ satisfies $\|\tilde{\eta}_n\|_3^2 \leq 3C^2D\mu$.
\end{corollary}

\begin{lemma}
The functional $\LL$ satisfies
$$\lim_{n \rightarrow \infty} \left[ \LL(\tilde{\eta}_n) - \sum_{i=1}^m \LL(\eta^{(i)})\right]=0, \qquad
\lim_{n \rightarrow \infty} \left\| \LL^\prime(\tilde{\eta}_n) - \sum_{i=1}^m \LL^\prime(\eta^{(i)})\right\|_1=0.$$
These limits also hold for $\KK$.
\end{lemma}
{\bf Proof.} Choose $\varepsilon>0$ and take $N$ large enough so that
$$\left\| \sum_{j=N+1}^\infty \tau_{S_j+(j-1)n}\eta^{(j)} \right\|_3 < \varepsilon,
\qquad \sum_{j=N+1}^\infty \LL(\eta^{(j)}) < \varepsilon
$$
if $m=\infty$ (see Proposition \ref{Uniform convergence over sequences}); write $N=m$ if $m<\infty$. Select $R>0$ large enough so that
$$\|\tau_{S_j}\eta^{(j)}\|_{H^3(|(x,z)|>R)} < \varepsilon, \qquad j=1,\ldots, N$$
and define
$$\zeta^{(j)}(x,z)=(\tau_{S_j}\eta^{(j)})(x,z)\chi\left(\frac{|(x,z)|}{R}\right), \qquad j=1,\ldots,N,$$
so that $\supp \zeta^{(j)} \subset B_{2R}(0)$ and $\|\zeta^{(j)}-\tau_{S_j}\eta^{(j)}\|_3 =O(\varepsilon)$
for $j=1,\ldots, N$; it follows that
$$\supp \tau_{(j-1)n}\zeta^{(j)} \subset {\mathbb R}^2\sm B_{n-2R}(0),
\qquad j=2,\ldots,N.$$

Observe that
\begin{eqnarray}
\lefteqn{\LL\left(\tau_{S_i}\eta^{(i)} + \sum_{j=i+1}^N \tau_{S_j+(j-i)n}\eta^{(j)}\right)
-\LL(\tau_{S_i}\eta^{(i)}) - \LL\left(\sum_{j=i+1}^N \tau_{S_j+(j-i)n}\eta^{(j)}\right)} \quad\nonumber \\
& & = \underbrace{\LL\left(\zeta^{(i)} + \sum_{j=i+1}^N \tau_{(j-i)n}\zeta^{(j)}\right)-\LL(\zeta^{(i)})
- \LL\left(\sum_{j=i+1}^N \tau_{(j-i)n}\zeta^{(j)}\right)}_{\displaystyle = o(1)}+ O(\varepsilon)\quad
\label{First step in L separation}
\end{eqnarray}
for $i=1,\ldots, N-1$ (see Theorem \ref{Key splitting theorem} in Appendix D), so that
\begin{eqnarray*}
\lim_{n \rightarrow \infty}
\LL\left(\tau_{S_i}\eta^{(i)} + \sum_{j=i+1}^N \tau_{S_j+(j-i)n}\eta^{(j)}\right)
& = & \LL(\tau_{S_i}\eta^{(i)}) + \lim_{n \rightarrow \infty} \LL\left(\sum_{j=i+1}^N \tau_{S_j+(j-i)n}\eta^{(j)}\right) \\
& = & \LL(\eta^{(i)})
+ \lim_{n \rightarrow \infty}\LL\left(\sum_{j=i+1}^N \tau_{S_j+(j-i)n}\eta^{(j)}\right)
\end{eqnarray*}
for $i=1,\ldots, N-1$; altogether one finds that
$$
\lim_{n \rightarrow \infty}\LL\left(\sum_{j=1}^N \tau_{S_j+(j-1)n}\eta^{(j)} \right)
= \sum_{j=1}^N \LL(\eta^{(j)}).
$$
This calculation establishes the result for $\LL$ in the case $m<\infty$; in the case $m=\infty$ we note that
\begin{eqnarray*}
\lefteqn{\lim_{n \rightarrow \infty}\LL\left(\sum_{j=1}^\infty \tau_{S_j+(j-1)n}\eta^{(j)}\right)
- \sum_{j=1}^\infty \LL(\eta^{(j)})} \quad \\
& & = \lim_{n \rightarrow \infty}\LL\left(\sum_{j=1}^N \tau_{S_j+(j-1)n}\eta^{(j)}\right)
- \sum_{j=1}^N \LL(\eta^{(j)}) + O(\varepsilon).
\end{eqnarray*}

The results for the other operators are obtained in a similar fashion (the $o(1)$ terms
in equation \eqn{First step in L separation} are identically zero for $\KK$ and $\KK^\prime$.)\qed

\begin{corollary} \label{Special PS property}
The sequence $\{\tilde{\eta}_n\}$ has the properties that
$$\lim_{n \rightarrow \infty} \JJ_\mu(\tilde{\eta}_n) = c_{\rho,\mu}, \qquad
\lim_{n \rightarrow \infty} \|\JJ_\mu^\prime(\tilde{\eta}_n)\|_1=0.$$
\end{corollary}
{\bf Proof.} Observe that
\begin{eqnarray*}
\lim_{n \rightarrow \infty} \JJ_\mu(\tilde{\eta}_n)
& = & \lim_{n \rightarrow \infty} \left[ \KK(\tilde{\eta}_n) + \frac{\mu^2}{\LL(\tilde{\eta}_n)} \right] \\
& = & \lim_{n \rightarrow \infty} \left[ \KK(\tilde{\eta}_n) + \frac{\mu^2}{\LL(\tilde{\eta}_n)^2}\LL(\tilde{\eta}_n) \right] \\
& = & \sum_{j=1}^m \KK(\eta^{(j)})
+ \lim_{n \rightarrow \infty}\left(\frac{\mu^2}{\LL(\tilde{\eta}_n)^2}\sum_{j=1}^m \LL(\eta^{(j)})\right) \\
& = & \sum_{j=1}^m \KK(\eta^{(j)})
+ \frac{\mu^2}{\displaystyle \lim_{n \rightarrow \infty} \LL(\tilde{\eta}_n)^2}\sum_{j=1}^m \LL(\eta^{(j)}) \\
& = & \sum_{j=1}^m \KK(\eta^{(j)}) + \sum_{j=1}^m \frac{(\mu_j^{(1)})^2}{\LL(\eta^{(j)})^2} \LL(\eta^{(j)}) \\
& = & \sum_{j=1}^m \KK(\eta^{(j)}) + \sum_{j=1}^m \frac{(\mu_j^{(1)})^2}{\LL(\eta^{(j)})} \\
& = & \sum_{j=1}^m \JJ_{\mu_j^{(1)}}(\eta^{(j)}) \\
& = & c_{\rho,\mu},
\end{eqnarray*}
where we have used the fact that
$$\frac{\mu}{\displaystyle \lim_{n \rightarrow \infty} \LL(\tilde{\eta}_n)} = \frac{\mu_j^{(1)}}{\LL(\eta^{(j)})}.$$
We similarly find that
\begin{eqnarray*}
\lim_{n \rightarrow \infty} \JJ_\mu^\prime(\tilde{\eta}_n)
& = & \lim_{n \rightarrow \infty} \left[ \KK^\prime(\tilde{\eta}_n) - \frac{\mu^2}{\LL(\tilde{\eta}_n)^2}
\LL^\prime(\tilde{\eta}_n) \right] \\
& = & \sum_{j=1}^m \KK^\prime(\eta^{(j)}) - \frac{\mu^2}{\displaystyle \lim_{n \rightarrow \infty}
\LL(\tilde{\eta}_n)^2} \sum_{j=1}^m \LL^\prime(\eta^{(j)}) \\
& = & 0
\end{eqnarray*}
because
$$\KK^\prime(\eta^{(j)})\ =\ \frac{(\mu_j^{(1)})^2}{\LL(\eta^{(j)})^2} \LL^\prime(\eta^{(j)})
\ =\ \frac{\mu^2}{\displaystyle \lim_{n \rightarrow \infty}
\LL(\tilde{\eta}_n)^2} \LL^\prime(\eta^{(j)});$$
the limits in these equations are taken in $H^1({\mathbb R}^2)$.\qed

\begin{proposition} \label{Special sequence minimises}
The sequence $\{\tilde{\eta}_n\}$ is a minimising sequence for $\JJ_\mu$ over $U \sm \{0\}$.
\end{proposition}
{\bf Proof.} Let us first note that $\{\tilde{\eta}_n\}$ is a minimising sequence for $\JJ_\mu$ over
$\tilde{U}\sm\{0\}$ since the existence of a minimising sequence $\{v_n\}$ for $\JJ_\mu$ over $\tilde{U}\sm\{0\}$
with $\lim_{n \rightarrow \infty} \JJ_\mu(v_n) < \lim_{n \rightarrow \infty} \JJ_\mu(\tilde{\eta}_n)$ would
lead to the contradiction
$$\lim_{n \rightarrow \infty} \JJ_{\rho,\mu}(v_n) = \lim_{n \rightarrow \infty} \JJ_\mu(v_n)
 < \lim_{n \rightarrow \infty} \JJ_\mu(\tilde{\eta}_n) = \lim_{n \rightarrow \infty} \JJ_{\rho,\mu}(\tilde{\eta}_n)
=c_{\rho,\mu}.$$
It follows from this fact and the estimate $\|\tilde{\eta}_n\|_3^2 \leq 3C^2D\mu$ that
$$\inf \{\JJ_\mu(\eta): \|\eta\|_3 \in (0,\tilde{M})\} = \inf \{\JJ_\mu(\eta): \|\eta\|_3 \in (0,\sqrt{3C^2D\mu})\}$$
for all $\tilde{M} \in (\sqrt{3C^2D\mu},M)$. The right-hand side of this equation does not depend
upon $\tilde{M}$; letting $\tilde{M} \rightarrow M$ on the left-hand side, one therefore
finds that
\begin{eqnarray*}
\inf \{\JJ_\mu(\eta): \|\eta\|_3 \in (0,M)\} & = & \inf \{\JJ_\mu(\eta): \|\eta\|_3 \in (0,\sqrt{3C^2D\mu})\} \\
& = & \lim_{n \rightarrow \infty} \JJ_\mu(\tilde{\eta}_n).
\end{eqnarray*}
\qed

\section{Strict sub-additivity} \label{SSA}
The goal of this section is to establish that the quantity
$$c_\mu = \inf_{\eta \in U\sm\{0\}} \JJ_\mu(\eta)$$
is a \emph{strictly sub-homogeneous} function of $\mu$, that is
$$
c_{a\mu} < a c_\mu, \qquad a>1.
$$
Its strict sub-homogeneity implies that $c_\mu$ also has the \emph{strict sub-additivity} property that
\begin{equation}
c_{\mu_1+\mu_2} < c_{\mu_1} + c_{\mu_2}, \qquad \mu_1, \mu_2>0 \label{Strict SH}
\end{equation}
(see Buffoni \cite[p.\ 48]{Buffoni04a}); inequality \eqn{Strict SH} plays a crucial role
in the variational theory presented in Section \ref{Stability} below.
The strict sub-homogeneity of $c_\mu$ follows from the fact that the function
\begin{equation}
a \mapsto a^{-\frac{5}{2}}\MM_{a^2\mu}(a\tilde{\eta}_n), \qquad a \in [1,2], \label{Decreasing function of a}
\end{equation}
where $\{\tilde{\eta}_n\}$ is the minimising sequence for $\JJ_\mu$ over $U\sm\{0\}$ constructed
in Section \ref{Special minimising sequence} above, is decreasing and strictly negative (see Lemma \ref{SH step 2}). The proof of the latter property,
which is given in Proposition \ref{SH step 1}, relies upon two observations, namely that
the `leading-order' term in $\MM_\mu(\tilde{\eta}_n)$ is
$-(\mu/\LL_2(\tilde{\eta}_n))^2\LL_3(\tilde{\eta}_n)$, and that the function \eqn{Decreasing function of a}
clearly has the required property if $\MM_\mu(\tilde{\eta}_n)$ is replaced by  $-(\mu/\LL_2(\tilde{\eta}_n))^2\LL_3(\tilde{\eta}_n)$
and $\LL_3(\tilde{\eta}_n)$ is positive. The idea behind the proof of Proposition \ref{SH step 1} is therefore to approximate
$\MM_\mu(\tilde{\eta}_n)$ with $-(\mu/\LL_2(\tilde{\eta}_n))^2\LL_3(\tilde{\eta}_n)$.

Writing
\begin{equation}
\MM_\mu(\eta) = - \frac{\mu^2\LL_3(\eta)}{\LL_2(\eta)^2} - \frac{\mu^2}{\LL_2(\eta)^2}(\LL_\mathrm{nl}(\eta)-\LL_3(\eta))
+\frac{\mu^2 \LL_\mathrm{nl}(\eta)^2}{\LL_2(\eta)^2\LL(\eta)}+\KK_\mathrm{nl}(\eta)
\label{lot in MM}
\end{equation}
and recalling that $\MM_\mu(\tilde{\eta}_n) \leq -c\mu^3$, one finds
that $-(\mu/\LL_2(\tilde{\eta}_n))^2\LL_3(\tilde{\eta}_n)$ is indeed the leading-order term in $\MM_\mu(\tilde{\eta}_n)$
(and is negative) provided that all other terms in the above formula
are $o(\mu^3)$;  straightforward estimates of the kind
$$\KK_j(\tilde{\eta}_n),\ \LL_j(\tilde{\eta}_n) = O(\|\tilde{\eta}_n\|_3^j)= O(\mu^{j/2})$$
however do not suffice for this purpose.
According to the calculations presented in Appendix B the function
$$\eta^\star_\mu(x,z)=\mu^2\Psi(\mu x, \mu^2 z), \qquad \Psi \in C_0^\infty({\mathbb R}^2),$$
whose length scales are those of the KP-I equation (cf.\ equation \eqn{AS explicit formula}),
has the property that $\MM_\mu(\eta^\star_\mu)=-(\mu/\LL_2(\eta^\star_\mu))^2\LL_3(\eta^\star_\mu)+o(\mu^3)$
(cf.\ equations \eqn{Also for SSH 1}--\eqn{Also for SSH 3}).
Although $\|\tilde{\eta}_n\|_3$ and $\|\eta_\mu^\star\|_3$ are both $O(\mu^{\frac{1}{2}})$, the function
$\eta^\star_\mu$ has the advantage that the $L^2({\mathbb R}^2)$-norms of its derivatives are higher
order with respect to $\mu$ (e.g.\ the quantity $\|\partial_x \eta^\star_\mu\|_0$ is not merely
$O(\mu^{\frac{1}{2}})$ but $O(\mu^{\frac{3}{2}})$). This fact allows one to obtain better estimates for
$\KK_j(\eta^\star_\mu)$ and $\LL_j(\eta^\star_\mu)$ (see below). Motivated by the expectation
that a minimiser, and hence a minimising sequence, should have the KP-I length scales, our
our strategy is therefore to show that $\tilde{\eta}_n$ is $O(\mu^{\frac{1}{2}})$ with respect to a norm on
$H^3({\mathbb R}^2)$ with weighted derivatives. To this end we consider the norm
$$
\nn \eta \nn_\alpha^2 := \int_{{\mathbb R}^2} \left(1+\mu^{-6\alpha}|k|^6+\mu^{-4\alpha}\frac{k_2^4}{|k|^4}\right)|\hat{\eta}|^2\dk
$$
and choose $\alpha>0$ as large as possible so that $\nn \tilde{\eta}_n \nn_\alpha$
is $O(\mu^{\frac{1}{2}})$.

We begin by observing that the norm of certain Fourier-multiplier operators
with respect to $\nn\cdot\nn_\alpha$ is proportional to a power of $\mu$
and using this fact to obtain some basic estimates.

\begin{proposition} \label{Powers of mu in multiplier estimates}
The estimates
$$\|\FF^{-1}[f_0 \hat{\eta}]\|_\infty \leq c\mu^{\frac{3\alpha}{2}}\nn \eta \nn_\alpha, \qquad \|\FF^{-1}[f_1 \hat{\eta}]\|_\infty \leq c\mu^{\frac{5\alpha}{2}}\nn \eta \nn_\alpha$$
and
$$\|\FF^{-1}[f_2 \hat{\eta}]\|_{L^4({\mathbb R}^2)}\leq c\mu^{\frac{11\alpha}{4}}\nn \eta \nn_\alpha$$
hold for all continuous functions $f_j: {\mathbb R}^2 \rightarrow {\mathbb R}$ which satisfy the estimates
$$|f_j(k)| \leq c|k|^j, \qquad j=0,1,2$$
for all $k \in {\mathbb R}^2$. 
\end{proposition}
{\bf Proof.} Using the Cauchy-Schwarz inequality, we find that
$$\|\FF^{-1}[f_j \hat{\eta}]\|_\infty^2 \leq c\||k|^{2j} \hat{\eta}\|_{L^1({\mathbb R}^2)}^2 \leq
c\left(\int_{{\mathbb R}^2} \frac{|k|^{2j}}{1+\mu^{-6\alpha}|k|^6 + \mu^{-4\alpha}k_2^4/|k|^4} \dk\right)
\nn \eta\nn_\alpha^2,$$
and
\begin{eqnarray*}
\int_{{\mathbb R}^2} \frac{|k|^{2j}}{1+\mu^{-6\alpha}|k|^6 + \mu^{-4\alpha}k_2^4/|k|^4} \dk
& = &
2\int_0^\infty\!\!\!\int_{-\frac{\pi}{2}}^{\frac{\pi}{2}} \frac{r^{2j+1}}{1+\mu^{-6\alpha}r^6 + \mu^{-4\alpha}\sin^4\theta}\dtheta\dr \\
& \leq & 
c\int_0^\infty\!\!\!\int_{-\frac{\pi}{2}}^{\frac{\pi}{2}}\frac{r^{2j+1}}{1+\mu^{-6\alpha}r^6 + \mu^{-4\alpha}\theta^4}\dtheta\dr \\
& \leq & 
c\int_0^\infty\!\!\!\int_0^\infty\frac{r^{2j+1}}{1+\mu^{-6\alpha}r^6 + \mu^{-4\alpha}\theta^4}\dtheta\dr \\
& = & c\mu^{(2j+3)\alpha}
\int_0^\infty\!\!\!\int_0^\infty\frac{r^{2j+1}}{1+r^6 + \theta^4}\dtheta\dr.
\end{eqnarray*}
The remaining estimate follows from the calculation
\begin{eqnarray*}
\lefteqn{\|\FF^{-1}[f_2 \hat{\eta}]\|_{L^4({\mathbb R}^2)}} \\
& \leq & c\||k|^2\hat{\eta}\|_{L^{4/3}({\mathbb R}^2)} \\
& = & c\left(\int_{{\mathbb R}^2} \!\!\left(\frac{|k|^8}{(1+\mu^{-6\alpha}|k|^6 + \mu^{-4\alpha}k_2^4/|k|^4)^2}\right)^{\!\! \frac{1}{3}}\!\!
\left(1+\mu^{-6\alpha}|k|^6+\mu^{-4\alpha}\frac{k_2^4}{|k|^4}\right)^{\!\!\frac{2}{3}}|\hat{\eta}|^\frac{4}{3}\dk\right)^{\!\! \frac{3}{4}} \\
& \leq & c\left(\int_{{\mathbb R}^2} \frac{|k|^8}{(1+\mu^{-6\alpha}|k|^6 + \mu^{-4\alpha}k_2^4/|k|^4)^2} \dk\right)^{\!\!\frac{1}{4}}
\nn \eta\nn_\alpha, \\
& \leq & c\left(\mu^{11\alpha}\int_0^\infty\!\!\!\int_0^\infty\frac{r^9}{(1+r^6 + \theta^4)^2}\dtheta\dr\right)^{\!\!\frac{1}{4}}
\nn \eta\nn_\alpha,
\end{eqnarray*}
in which the Hausdorff-Young and H\"{o}lder inequalities have been used.\qed

\begin{corollary} \hspace{1cm} \label{Specific mu powers}
\begin{list}{(\roman{count})}{\usecounter{count}}
\item
The estimates
$$
\|\eta\|_\infty \leq c \mu^{\frac{3\alpha}{2}} \nn \eta \nn_\alpha,
\quad \|\eta\|_Z \leq c \mu^{\frac{3\alpha}{2}} \nn \eta \nn_\alpha,
\quad \|K^0\eta\|_\infty \leq c \mu^{\frac{3\alpha}{2}} \nn \eta \nn_\alpha
$$
and
$$
\|\eta_{xx} \|_{L^4({\mathbb R}^2)} \leq c\mu^{\frac{11\alpha}{4}} \nn \eta \nn_\alpha, \quad
\|\eta_x\|_\infty \leq c \mu^{\frac{5\alpha}{2}} \nn \eta \nn_\alpha$$
hold for all $\eta \in H^3({\mathbb R}^2)$ and remain valid when $K^0$ is replaced by $L^0$ or $M^0$,
$\eta_x$ is replaced by $\eta_z$ and $\eta_{xx}$ is replaced by $\eta_{xz}$ or $\eta_{zz}$.
\item
The estimates
$$\|u K^0 \eta_x \|_0 \leq c\mu^{\frac{5\alpha}{2}}\|u\|_1 \nn \eta \nn_\alpha$$
$$\|u K^0 (\eta_{xx}) \|_0 \leq c\mu^{\frac{11\alpha}{4}}\|u\|_2 \nn \eta \nn_\alpha$$
hold for all $\eta$, $u \in H^3({\mathbb R}^2)$ and remain valid when $K^0$ is replaced by $L^0$ or $M^0$,
$\eta_x$ is replaced by $\eta_z$
and $\eta_{xx}$ is replaced by $\eta_{xz}$ or $\eta_{zz}$.
\end{list}
\end{corollary}
{\bf Proof.} (i) The first, fourth and fifth estimates are a direct consequence of Proposition \ref{Powers of mu in multiplier estimates},
while the second and third follow from the calculations
$$
\|\eta\|_Z\ =\ \|\eta\|_{1,\infty}+ \|\eta_{xx}\|_1 + \|\eta_{xz}\|_1+ \|\eta_{zz}\|_1\ \leq\ c(\mu^{\frac{3\alpha}{2}}\nn \eta \nn_\alpha+\mu^{2\alpha}\nn\eta\nn_\alpha)
$$
and
$$
|\FF[K^0(\eta)| \leq \hspace{-5mm}\underbrace{\frac{k_1^2}{|k|^2}}_{\displaystyle := g_1(k)}\hspace{-5mm} |\hat{\eta}|
+ \underbrace{\frac{k_1^2}{|k|^2}(|k|\coth|k|-1)}_{\displaystyle := g_2(k)} |\hat{\eta}|,
$$
where $g_1(k)=O(1)$, $g_2(k)=O(|k|)$, so that
$$
\|K^0\eta\|_0 \leq \|g_1 \hat{\eta}\|_0 + \|g_2 \hat{\eta}\|_0 \leq c\mu^{\frac{3\alpha}{2}}\nn \eta \nn_\alpha.
$$

(ii) Observe that
\begin{eqnarray*}
|\FF[K^0(\eta_x)]| & \leq & \underbrace{\frac{k_1^2}{|k|^2} |k_1|}_{\displaystyle := g_3(k)}\hspace{-1mm} |\hat{\eta}|
+ \underbrace{\frac{k_1^2}{|k|^2}(|k|\coth|k|-1)|k_1|}_{\displaystyle := g_4(k)} |\hat{\eta}|, \\
|\FF[K^0(\eta_{xx})]| & \leq & \underbrace{\frac{k_1^2}{|k|^2} k_1^2}_{\displaystyle := g_5(k)}\hspace{-2mm} |\hat{\eta}|
+ \underbrace{\frac{k_1^2}{|k|^2}(|k|\coth|k|-1)k_1^2}_{\displaystyle := g_6(k)} |\hat{\eta}|,
\end{eqnarray*}
where $g_3(k)=O(|k|)$, $g_4(k)$, $g_5(k)=O(|k|^2)$ and $g_6(k)=O(|k|^3)$, so that
\begin{eqnarray*}
\|u K^0\eta_x\|_0 & \leq & \|u\|_0\|\FF^{-1}[g_3 \hat{\eta}]\|_\infty
+ \|u\|_{L^4({\mathbb R}^2)} \|\FF^{-1}[g_4 \hat{\eta}]\|_{L^4({\mathbb R}^2)} \\
& \leq & c (\mu^{\frac{5\alpha}{2}}\|u\|_0\nn \eta \nn_\alpha + \mu^{\frac{11\alpha}{4}}\|u\|_1 \nn \eta \nn_\alpha),\\
\\
\|u K^0\eta_{xx}\|_0 & \leq & \|u\|_{L^4({\mathbb R}^2)} \|\FF^{-1}[g_5 \hat{\eta}]\|_{L^4({\mathbb R}^2)} + \|u\|_\infty \|g_6\hat{\eta}\|_0\\
& \leq & c (\mu^{\frac{11\alpha}{4}}\|u\|_1 \nn \eta \nn_\alpha + \mu^{3\alpha}\|u\|_2\nn \eta \nn_\alpha).
\end{eqnarray*}

The remaining estimates are obtained in a similar fashion.\qed

The next step is to show that any function $\tilde{\eta} \in U\sm\{0\}$ which satisfies
\begin{equation}
\|\tilde{\eta}\|_3^2 \leq c \mu, \quad \JJ_\mu(\tilde{\eta})<2\mu, \quad \|\JJ_\mu^\prime(\tilde{\eta})\|_1\leq c \mu^N
\label{Properties of tildeeta}
\end{equation}
for a sufficiently large natural number $N \in {\mathbb N}$ has the requisite property that 
$\nn \tilde{\eta} \nn_\alpha^2 = O(\mu)$ for $\alpha < 1$. Notice that $\JJ_\mu(\tilde{\eta})<2\mu$ implies
$\LL(\tilde{\eta})>\mu/2$ and hence $\LL_2(\tilde{\eta}) > c \mu$; the following result gives another useful
inequality for $\tilde{\eta}$.

\begin{proposition} \label{A/B estimate on speed}
The function $\tilde{\eta}$ satisfies the inequality
$$
\RR_1(\tilde{\eta}) + \tilde{\MM}_\mu(\tilde{\eta}) \leq \frac{\mu}{\LL(\tilde{\eta})} -1 \leq \RR_2(\tilde{\eta}) + \tilde{\MM}_\mu(\tilde{\eta}),
$$
where
\begin{eqnarray*}
\RR_1(\tilde{\eta}) & = & -\frac{\langle \JJ_\mu^\prime(\tilde{\eta}),\tilde{\eta}\rangle_0}{4\mu} + \frac{\langle \MM_\mu^\prime(\tilde{\eta}),\tilde{\eta}\rangle_0}{4\mu}, \nonumber \\
\RR_2(\tilde{\eta}) & = & -\frac{\langle \JJ_\mu^\prime(\tilde{\eta}),\tilde{\eta}\rangle_0}{4\mu} + \frac{\langle \MM_\mu^\prime(\tilde{\eta}),\tilde{\eta}\rangle_0}{4\mu}-\frac{\MM_\mu(\tilde{\eta})}{2\mu},
\end{eqnarray*}
and
$$
\tilde{\MM}_\mu(\tilde{\eta}) = \frac{\mu}{\LL(\tilde{\eta})}-\frac{\mu}{\LL_2(\tilde{\eta})}.
$$
\end{proposition}
{\bf Proof.} Taking the scalar product of the equation
$$\JJ_\mu^\prime(\tilde{\eta})=\KK_2^\prime(\tilde{\eta}) - \left(\frac{\mu}{\LL_2(\tilde{\eta})}\right)^{\!\!2}\LL_2^\prime(\tilde{\eta})+\MM_\mu^\prime(\tilde{\eta})$$
with $\tilde{\eta}$ yields the identity
$$
\frac{\mu}{\LL_2(\tilde{\eta})} = - \frac{\langle \JJ_\mu^\prime(\tilde{\eta}),\tilde{\eta}\rangle_0}{4\mu} + \frac{1}{2\mu}\left[\KK_2(\tilde{\eta}) + \frac{\mu^2}{\LL_2(\tilde{\eta})}\right] + \frac{\langle \MM_\mu^\prime(\tilde{\eta}),\tilde{\eta}\rangle_0}{4\mu}.
$$
The assertion follows by estimating the quantity in square brackets from above and below by means of the
inequalities
$$2\mu \leq \KK_2(\tilde{\eta}) + \frac{\mu^2}{\LL_2(\tilde{\eta})} = \JJ_\mu(\tilde{\eta}) - \MM_\mu(\tilde{\eta}) \leq 2\mu - \MM_\mu(\tilde{\eta})$$
(see \eqn{inf less than 2mu} and \eqn{Lower bound on quadratic part}).\qed

\begin{corollary} \label{Speed estimate}
The function $\tilde{\eta}$ satsfies the inequality
$$\left|\frac{\mu}{\LL(\tilde{\eta})}-1\right| \leq c(\mu^{N-\frac{1}{2}}+\mu^{\frac{3\alpha}{2}-\frac{1}{2}}\nn\tilde{\eta}\nn_\alpha^2).$$
\end{corollary}
{\bf Proof.} Remark \ref{Estimates for KKnl} implies that
$\KK_\mathrm{nl}(\tilde{\eta})$ and  $\langle\KK_\mathrm{nl}^\prime(\tilde{\eta}),\tilde{\eta}\rangle_0$
are $O(\|\tilde{\eta}\|_Z\|\tilde{\eta}\|_3^3)$,  while Proposition \ref{Z estimate for LL} shows that
$\LL_\mathrm{nl}(\tilde{\eta}) =O(\|\tilde{\eta}\|_Z \|\tilde{\eta}\|_3^2)$
and
$$
|\langle\LL_\mathrm{nl}^\prime(\tilde{\eta}),\tilde{\eta}\rangle_0|
\ \leq\ 3|\LL_3(\tilde{\eta})|+\|\LL_\mathrm{nl}^\prime(\tilde{\eta})-\LL_3^\prime(\tilde{\eta})\|_0\|\tilde{\eta}\|_0
\ \leq\ c\|\tilde{\eta}\|_Z \|\tilde{\eta}\|_3^2.
$$
These four quantities are therefore all $O(\mu^{\frac{3\alpha}{2}+\frac{1}{2}}\nn \tilde{\eta} \nn_\alpha^2)$
because $\|\tilde{\eta}\|_Z \leq c\mu^{\frac{3\alpha}{2}}\nn \tilde{\eta} \nn_\alpha$ (Corollary \ref{Specific mu powers}(i)) and
$\|\tilde{\eta}\|_3^2 \leq \|\tilde{\eta}\|_3 \nn \tilde{\eta} \nn \leq c \mu^{\frac{1}{2}}\nn \tilde{\eta} \nn_\alpha$.
Writing
$$
\MM_\mu(\tilde{\eta}) = \KK_\mathrm{nl}(\tilde{\eta}) - \frac{\mu^2\LL_\mathrm{nl}(\tilde{\eta})}{\LL(\tilde{\eta})\LL_2(\tilde{\eta})}, \qquad
\tilde{\MM}_\mu(\tilde{\eta}) = - \frac{\mu\LL_\mathrm{nl}(\tilde{\eta})}{\LL(\tilde{\eta})\LL_2(\tilde{\eta})}
$$
and
$$\langle \MM_\mu^\prime(\tilde{\eta}),\tilde{\eta} \rangle_0 = \langle \KK_\mathrm{nl}^\prime(\tilde{\eta}),\tilde{\eta} \rangle_0
-\frac{\mu^2\langle \LL_\mathrm{nl}^\prime(\tilde{\eta}),\tilde{\eta}\rangle_0}{\LL(\tilde{\eta})\LL_2(\tilde{\eta})}
+\frac{2\mu^2\LL_\mathrm{nl}(\tilde{\eta})}{\LL(\tilde{\eta})\LL_2(\tilde{\eta})} + \frac{2\mu^2\LL_\mathrm{nl}(\tilde{\eta})}{\LL(\tilde{\eta})^2}
+ \frac{\mu^2\LL_\mathrm{nl}(\tilde{\eta}) \langle \LL_\mathrm{nl}^\prime(\tilde{\eta}),\tilde{\eta}\rangle_0}{\LL(\tilde{\eta})^2\LL_2(\tilde{\eta})},$$
one finds that
\begin{equation}
|\RR_1(\tilde{\eta})|,\ |\RR_2(\tilde{\eta})|,\ |\tilde{\MM}_\mu(\tilde{\eta})| \leq c\mu^{\frac{3\alpha}{2}-\frac{1}{2}}\nn\tilde{\eta}\nn_\alpha^2. \label{Estimates on RR and tildeMM}
\end{equation}
The assertion follows from
\eqn{Estimates on RR and tildeMM}, the estimate $|\langle \JJ_\mu(\tilde{\eta}),\tilde{\eta} \rangle_0| \leq \|\JJ_\mu(\tilde{\eta})\|_0\|\tilde{\eta}\|_3 
\leq c\mu^{N+\frac{1}{2}}$ and Proposition \ref{A/B estimate on speed}.\qed

We proceed by applying the operators
$$-\langle (\cdot)_x,\tilde{\eta}_{xxx} \rangle_0 -\langle (\cdot)_z,\tilde{\eta}_{zzz} \rangle_0, \qquad \left\langle \FF[\cdot], \frac{k_2^2}{|k|^2}\hat{\tilde{\eta}}\right\rangle_{\!\!0}$$
to the identity
$$\JJ_\mu^\prime(\tilde{\eta})=\KK_2^\prime(\tilde{\eta}) - \LL_2^\prime(\tilde{\eta})
+ \KK_\mathrm{nl}^\prime(\tilde{\eta})-\left(\frac{\mu}{\LL(\tilde{\eta})}-1\right)\!\!\!\left(\frac{\mu}{\LL(\tilde{\eta})}+1\right)\LL_2^\prime(\tilde{\eta}) - \left(\frac{\mu}{\LL(\tilde{\eta})}\right)^{\!\!2}\LL_\mathrm{nl}^\prime(\tilde{\eta})$$
and estimating the resulting equations, namely
\begin{eqnarray}
\lefteqn{\int_{{\mathbb R}^2} \frac{k_2^2}{|k|^2}\left(1+\beta|k|^2 - \frac{k_1^2}{|k|^2} |k| \coth |k|\right)|\hat{\tilde{\eta}}|^2 \dk} \nonumber \\[1mm]
& & = \left\langle \FF[\JJ_\mu^\prime(\tilde{\eta})],\frac{k_2^2}{|k|^2}\hat{\tilde{\eta}}\right\rangle_{\!\!0}
-\left\langle \FF[\KK_\mathrm{nl}^\prime(\tilde{\eta})],\frac{k_2^2}{|k|^2}\hat{\tilde{\eta}}\right\rangle_{\!\!0}
+ \left(\frac{\mu}{\LL(\tilde{\eta})}\right)^{\!\!2}\left\langle \FF[\LL_\mathrm{nl}^\prime(\tilde{\eta})],\frac{k_2^2}{|k|^2}\hat{\tilde{\eta}}\right\rangle_{\!\!0}
\nonumber \\
& & \qquad\quad \mbox{} + \left(\frac{\mu}{\LL(\tilde{\eta})}-1\right)\!\!\!\left(\frac{\mu}{\LL(\tilde{\eta})}+1\right)\!\!\int_{{\mathbb R}^2} \frac{k_1^2}{|k|^2}|k|\coth|k|\frac{k_2^2}{|k|^2}|\hat{\tilde{\eta}}|^2\dk \label{Scaled norm 1}
\end{eqnarray}
and
\begin{eqnarray}
\lefteqn{\int_{{\mathbb R}^2} (k_1^4+k_2^4)\left(1+\beta|k|^2 - \frac{k_1^2}{|k|^2} |k| \coth |k|\right)|\hat{\tilde{\eta}}|^2 \dk} \nonumber \\[1mm]
& & = -\langle (\JJ_\mu^\prime(\tilde{\eta}))_x,\tilde{\eta}_{xxx}\rangle_0 -\langle (\JJ_\mu^\prime(\tilde{\eta}))_z,\tilde{\eta}_{zzz}\rangle_0 
+\langle (\KK_\mathrm{nl}^\prime(\tilde{\eta}))_x,\tilde{\eta}_{xxx}\rangle_0 + \langle (\KK_\mathrm{nl}^\prime(\tilde{\eta}))_z,\tilde{\eta}_{zzz}\rangle_0 \nonumber \\
& & \qquad\quad\mbox{}- \left(\frac{\mu}{\LL(\tilde{\eta})}\right)^{\!\!2}\langle (\LL_\mathrm{nl}^\prime(\tilde{\eta}))_x,\tilde{\eta}_{xxx}\rangle_0
- \left(\frac{\mu}{\LL(\tilde{\eta})}\right)^{\!\!2}\langle (\LL_\mathrm{nl}^\prime(\tilde{\eta}))_z,\tilde{\eta}_{zzz}\rangle_0\nonumber \\
& & \qquad\quad \mbox{} + \left(\frac{\mu}{\LL(\tilde{\eta})}-1\right)\!\!\!\left(\frac{\mu}{\LL(\tilde{\eta})}+1\right)\!\!\int_{{\mathbb R}^2} \frac{k_1^2}{|k|^2}|k|\coth|k|(k_1^4+k_2^4)|\hat{\tilde{\eta}}|^2\dk \label{Scaled norm 2}
\end{eqnarray}

\begin{lemma} \label{RHS of scaled norm inequalities}
Suppose that $N \geq \max(6\alpha + 1/2,1+3\alpha/2)$.
The right-hand sides of \eqn{Scaled norm 1} and \eqn{Scaled norm 2} are
respectively $O(\mu^{4\alpha}(\mu+\mu^{-\frac{\alpha}{2}+\frac{1}{2}}\nn \tilde{\eta} \nn_\alpha^2
+\mu^{-\frac{\alpha}{2}-\frac{1}{2}}\nn \tilde{\eta} \nn_\alpha^4))$ and $O(\mu^{6\alpha}(\mu+\mu^{-\frac{\alpha}{2}+\frac{1}{2}}\nn \tilde{\eta} \nn_\alpha^2
+\mu^{-\frac{\alpha}{2}-\frac{1}{2}}\nn \tilde{\eta} \nn_\alpha^4))$.
\end{lemma}
{\bf Proof.} We examine each term on the right-hand sides of \eqn{Scaled norm 1} and \eqn{Scaled norm 2} separately.
\begin{itemize}
\item
Clearly
$$\left|\left\langle \FF[\JJ_\mu^\prime(\tilde{\eta})],\frac{k_2^2}{|k|^2}\hat{\tilde{\eta}}\right\rangle_{\!\!0}\right| \leq \|\JJ^\prime_\mu(\tilde{\eta})\|_0 \|\tilde{\eta}\|_3 \leq c\mu^{4\alpha}\mu^{N+\frac{1}{2}-4\alpha}$$
and
$$|\langle(\JJ_\mu^\prime(\tilde{\eta}))_x,\tilde{\eta}_{xxx}\rangle_0|
\leq \|(\JJ_\mu^\prime(\tilde{\eta}))_x\|_0\|\tilde{\eta}\|_3 \leq c \mu^{6\alpha}\mu^{N+\frac{1}{2}-6\alpha}
\nn \tilde{\eta} \nn_\alpha^2,$$
$$|\langle(\JJ_\mu^\prime(\tilde{\eta}))_z,\tilde{\eta}_{zzz}\rangle_0|
\leq \|(\JJ_\mu^\prime(\tilde{\eta}))_z\|_0\|\tilde{\eta}\|_3 \leq c \mu^{6\alpha}\mu^{N+\frac{1}{2}-6\alpha}
\nn \tilde{\eta} \nn_\alpha^2$$
for $\alpha \leq 1$.
\item
Define
\begin{eqnarray*}
h_1(a,b) & = & \frac{\beta a(a^2+b^2)}{(1+\sqrt{1+a^2+b^2})\sqrt{1+a^2+b^2}}, \\
h_2(a,b) & = & \frac{\beta b(a^2+b^2)}{(1+\sqrt{1+a^2+b^2})\sqrt{1+a^2+b^2}},
\end{eqnarray*}
so that
$$\KK_\mathrm{nl}(\tilde{\eta}) = (h_1(\tilde{\eta}_x,\tilde{\eta}_z))_x + (h_2(\tilde{\eta}_x,\tilde{\eta}_z))_z,$$
and observe that $h_1, h_2: {\mathbb R}^2 \rightarrow {\mathbb R}$ are analytic at the origin, where they
have a third-order zero.

$\quad$It follows from the equation
$$(h_1(\tilde{\eta}_x,\tilde{\eta}_z))_x =  \partial_1 h_1(\tilde{\eta}_x,\tilde{\eta}_z)\tilde{\eta}_{xx} + \partial_2 h_1(\tilde{\eta}_x,\tilde{\eta}_z)\tilde{\eta}_{xz}$$
that
$$\|(h_1(\tilde{\eta}_x,\tilde{\eta}_z))_x\|_0  \leq c(\|\tilde{\eta}_x\|_\infty + \|\tilde{\eta}_z\|_\infty)\|\tilde{\eta}\|_3^2 \leq c\mu^{\frac{5\alpha}{2}+1}\nn \tilde{\eta} \nn_\alpha$$
(see Corollary \ref{Specific mu powers}(i)) and the same estimate clearly holds for $(h_2(\tilde{\eta}_x,\tilde{\eta}_z))_z$; we conclude that
$$\|\KK_\mathrm{nl}^\prime(\tilde{\eta})\|_0 \leq c\mu^{\frac{5\alpha}{2}+1}\nn \tilde{\eta} \nn_\alpha$$
and hence that
$$\left|
\left\langle \FF[\KK_\mathrm{nl}^\prime(\tilde{\eta})],\frac{k_2^2}{|k|^2}\hat{\tilde{\eta}}\right\rangle_{\!\!0}
\right|
\leq \|\KK_\mathrm{nl}^\prime(\tilde{\eta})\|_0 \left\|\frac{k_2^2}{|k|^2}\hat{\tilde{\eta}}\right\|_0 \leq c\mu^{4\alpha}\mu^{\frac{\alpha}{2}+1}\nn \tilde{\eta} \nn_\alpha^2.$$
Similarly, the equation
\begin{eqnarray*}
\lefteqn{(h_1(\tilde{\eta}_x,\tilde{\eta}_z))_{xx}
= \partial_1h_1(\tilde{\eta}_x,\tilde{\eta}_z)\tilde{\eta}_{xxx}+\partial_1^2h_1(\tilde{\eta}_x,\tilde{\eta}_z)\tilde{\eta}_{xx}^2} \hspace{1in}\\
& & \mbox{} + 2\partial^2_{1,2}h_1(\tilde{\eta}_x,\tilde{\eta}_z)\tilde{\eta}_{xx}\tilde{\eta}_{xz}
+ \partial_2^2h_1(\tilde{\eta}_x,\tilde{\eta}_z)\tilde{\eta}_{xz}^2 + \partial_2h_1(\tilde{\eta}_x,\tilde{\eta}_z)\tilde{\eta}_{xxz}
\end{eqnarray*}
implies that
\begin{eqnarray*}
\|(h_1(\tilde{\eta}_x,\tilde{\eta}_z))_{xx}\|_0 & \leq & c(\|\tilde{\eta}\|_3^2\|\tilde{\eta}_{xxx}\|_0 + (\|\tilde{\eta}_x\|_\infty+\|\tilde{\eta}_z\|_\infty)\|\tilde{\eta}\|_3^2
+\|\tilde{\eta}\|_3^2\|\tilde{\eta}_{xxz}\|_0) \\
& \leq & c (\mu^{3\alpha+1}\nn \tilde{\eta} \nn_{\alpha}+ \mu^{\frac{5\alpha}{2}+1}\nn \tilde{\eta} \nn_{\alpha}) \\
& \leq & c \mu^{\frac{5\alpha}{2}+1}\nn \tilde{\eta} \nn_\alpha,
\end{eqnarray*}
and analogous calculations yield
$$\|(h_2(\tilde{\eta}_x,\tilde{\eta}_z))_{xz}\|_0,\ \|(h_2(\tilde{\eta}_x,\tilde{\eta}_z))_{xz}\|_0,\ \|(h_2(\tilde{\eta}_x,\tilde{\eta}_z))_{xz}\|_0 \leq c
\mu^{\frac{5\alpha}{2}+1}\nn \tilde{\eta} \nn_\alpha.$$
Altogether one finds that
$$\|(\KK_\mathrm{nl}^\prime(\tilde{\eta}))_x\|_0,\ \|(\KK_\mathrm{nl}^\prime(\tilde{\eta}))_z\|_0 \leq c
\mu^{\frac{5\alpha}{2}+1}\nn \tilde{\eta} \nn_\alpha$$
and hence that
$$|\langle(\KK_\mathrm{nl}^\prime(\tilde{\eta}))_x,\tilde{\eta}_{xxx}\rangle_0|
\leq \|(\KK_\mathrm{nl}^\prime(\tilde{\eta}))_x\|_0\|\tilde{\eta}_{xxx}\|_0 \leq c \mu^{6\alpha}\mu^{-\frac{\alpha}{2}+1}
\nn \tilde{\eta} \nn_\alpha^2,$$
$$|\langle(\KK_\mathrm{nl}^\prime(\tilde{\eta}))_z,\tilde{\eta}_{zzz}\rangle_0|
\leq \|(\KK_\mathrm{nl}^\prime(\tilde{\eta}))_z\|_0\|\tilde{\eta}_{zzz}\|_0 \leq c \mu^{6\alpha}\mu^{-\frac{\alpha}{2}+1}
\nn \tilde{\eta} \nn_\alpha^2.$$

\item
Estimating the expression for $\LL_3^\prime(\tilde{\eta})$ given in Lemma \ref{Formulae for L3prime and L4prime} using
Corollary \ref{Specific mu powers}(i), we find that
\begin{eqnarray*}
\|\LL_3^\prime(\tilde{\eta})\|_0 & \leq & c(\|\tilde{\eta}_x\|_\infty\|\tilde{\eta}_x\|_0 + \|K^0\tilde{\eta}\|_\infty\|K^0\tilde{\eta}\|_0 + \|L^0\tilde{\eta}\|_\infty\|L^0\tilde{\eta}\|_0 + \|\tilde{\eta}\|_Z \|\tilde{\eta}\|_3) \\
& \leq & c(\|\tilde{\eta}\|_Z + \|K^0\tilde{\eta}\|_\infty + \|L^0\tilde{\eta}\|_\infty)\|\tilde{\eta}\|_3 \\
& \leq & c\mu^{\frac{3\alpha}{2}+\frac{1}{2}}\nn \tilde{\eta} \nn_\alpha.
\end{eqnarray*}
Because
$$\|\LL_\mathrm{nl}^\prime(\tilde{\eta})-\LL_3^\prime(\tilde{\eta})\|_0 \leq c \|\tilde{\eta}\|_Z \|\tilde{\eta}\|_3^2 \leq c\mu^{\frac{3\alpha}{2}+1}\nn \tilde{\eta} \nn_\alpha$$
(see Proposition \ref{Z estimate for LL}), one concludes that
$$\|\LL_\mathrm{nl}^\prime(\tilde{\eta})\|_0 \leq c\mu^{\frac{3\alpha}{2}+\frac{1}{2}}\nn \tilde{\eta} \nn_\alpha$$
and hence that
$$\left|
\left\langle \FF[\LL_\mathrm{nl}^\prime(\tilde{\eta})],\frac{k_2^2}{|k|^2}\hat{\tilde{\eta}}\right\rangle_{\!\!0}
\right|
\leq \|\LL_\mathrm{nl}^\prime(\tilde{\eta})\|_0 \left\|\frac{k_2^2}{|k|^2}\hat{\tilde{\eta}}\right\|_0 \leq c\mu^{4\alpha}\mu^{-\frac{\alpha}{2}+\frac{1}{2}}\nn \tilde{\eta} \nn_\alpha^2.$$

$\quad$It follows from Corollary \ref{Explicit formulae for LL3} that
\begin{eqnarray*}
\lefteqn{(\LL^\prime_3(\tilde{\eta}))_x = - \frac{3}{2}\tilde{\eta}_x\tilde{\eta}_{xx} - \tilde{\eta}\tilde{\eta}_{xxx} - (K^0\tilde{\eta})(K^0\tilde{\eta}_x) - (L^0 \tilde{\eta})(L^0\tilde{\eta}_x)} \\
& & \qquad\qquad\quad\mbox{}
-K^0(\tilde{\eta}_x K^0 \tilde{\eta}) - L^0(\tilde{\eta}_xL^0\tilde{\eta}) - K^0(\tilde{\eta} K^0\tilde{\eta}_x) - L^0(\tilde{\eta} L^0 \tilde{\eta}_x).
\end{eqnarray*}
Using this formula, the estimates
\begin{eqnarray*}
\|K^0(\tilde{\eta}_x K^0 \tilde{\eta})\|_0 & \leq & c\|\tilde{\eta}_x K^0 \tilde{\eta}\|_1 \\
& \leq & c (\|\tilde{\eta}_x K^0 \tilde{\eta}\|_0 + \|\tilde{\eta}_{xx} K^0 \tilde{\eta}\|_0 + \|\tilde{\eta}_x K^0 \tilde{\eta}_x\|_0) \\
& \leq & c (\|\tilde{\eta}_x\|_\infty \|K^0 \tilde{\eta}\|_0 + \|\tilde{\eta}_{xx}\|_{L^4({\mathbb R}^2)} \|K^0\tilde{\eta}\|_{L^4({\mathbb R}^2)}+ \|\tilde{\eta}_x K^0 \tilde{\eta}_x\|_0) \\
& \leq & c\mu^{\frac{5\alpha}{2}}\|\tilde{\eta}\|_3\nn \tilde{\eta} \nn_\alpha, \\
\\
\|K^0(\tilde{\eta} K^0 \tilde{\eta}_x)\|_0 & \leq & c\|\tilde{\eta} K^0 \tilde{\eta}_x\|_1 \\
& \leq & c (\|\tilde{\eta} K^0 \tilde{\eta}_x\|_0 + \|\tilde{\eta}_x K^0 \tilde{\eta}_x\|_0 + \|\tilde{\eta} K^0 \tilde{\eta}_{xx}\|_0) \\
& \leq & c\mu^{\frac{5\alpha}{2}}\|\tilde{\eta}\|_3\nn \tilde{\eta} \nn_\alpha
\end{eqnarray*}
(see Corollary \ref{Specific mu powers})
and the corresponding calculations for $L^0(\tilde{\eta}_x L^0\tilde{\eta})$ and $L^0(\tilde{\eta} L^0 \tilde{\eta}_x)$, one finds that
\begin{eqnarray*}
\|\LL^\prime_3(\tilde{\eta}))_x\|_0 & \leq & c (\|\tilde{\eta}_x\|_\infty \|\tilde{\eta}_{xx}\|_0 + \|\tilde{\eta}\|_\infty \|\tilde{\eta}_{xxx}\|_0 \\
& & \qquad\qquad\mbox{}
+\mu^{\frac{5\alpha}{2}}(\|K^0\tilde{\eta}\|_1 + \|L^0\tilde{\eta}\|_1)\nn \tilde{\eta} \nn _\alpha + \mu^{\frac{5\alpha}{2}}\|\tilde{\eta}\|_3\nn \tilde{\eta} \nn_\alpha) \\
& \leq & c\mu^{\frac{5\alpha}{2}}\|\tilde{\eta}\|_3\nn \tilde{\eta} \nn_\alpha \\
& \leq & c\mu^{\frac{5\alpha}{2}+\frac{1}{2}}\nn \tilde{\eta} \nn_\alpha,
\end{eqnarray*}
in which Corollary \ref{Specific mu powers} has again been used.

$\quad$Treating the right-hand side of the formula
\begin{eqnarray*}
\lefteqn{(\LL^\prime_4(\tilde{\eta}))_x} \\
& & = K^0\tilde{\eta}_x\Big(\tilde{\eta}\tilde{\eta}_{xx}+K^0(\tilde{\eta} K^0\tilde{\eta}) + L^0(\tilde{\eta} L^0 \tilde{\eta})\Big)
+ K^0\tilde{\eta}\Big(\tilde{\eta}\tilde{\eta}_{xx}+K^0(\tilde{\eta} K^0\tilde{\eta}) + L^0(\tilde{\eta} L^0 \tilde{\eta})\Big)_{\!\!x} \\
& & \mbox{} + L^0\tilde{\eta}_x\Big(\tilde{\eta}\tilde{\eta}_{xz}+L^0(\tilde{\eta} K^0\tilde{\eta}) + M^0(\tilde{\eta} L^0 \tilde{\eta})\Big)
+ L^0\tilde{\eta}\Big(\tilde{\eta}\tilde{\eta}_{xz}+L^0(\tilde{\eta} K^0\tilde{\eta}) + M^0(\tilde{\eta} L^0 \tilde{\eta})\Big)_{\!\!x} \\
& & \mbox{} + (K^2(\tilde{\eta})\tilde{\eta})_x
\end{eqnarray*}
in the same manner and using the additional estimate
$$\|(K^2(\tilde{\eta})\tilde{\eta})_x\|_0 \leq c\|\tilde{\eta}\|_Z^2\|\tilde{\eta}\|_3 \leq c \mu^{3\alpha}\nn \tilde{\eta} \nn_\alpha^3,$$
(see Proposition \ref{Z estimate for LL}), one finds that
$$
\|(\LL^\prime_4(\tilde{\eta}))_x\|_0 \leq c (\mu^{\frac{5\alpha}{2}+\frac{1}{2}}\nn \tilde{\eta} \nn_\alpha + \mu^{3\alpha}\nn \tilde{\eta} \nn_\alpha^3).$$
Finally
$$\|(\LL_\mathrm{nl}^\prime(\tilde{\eta}) - \LL_3^\prime(\tilde{\eta}) - \LL_4^\prime(\tilde{\eta}))_x\|_0
\leq \|\tilde{\eta}\|_Z^2 \|\tilde{\eta}\|_3^2 \leq c\mu^{3\alpha}\nn \tilde{\eta} \nn_\alpha^3$$
(see Proposition \ref{Z estimate for LL}), and altogether these calculations yield
$$\|\LL_\mathrm{nl}^\prime(\tilde{\eta}))_x\|_0 \leq c (\mu^{\frac{5\alpha}{2}+\frac{1}{2}}\nn \tilde{\eta} \nn_\alpha + \mu^{3\alpha}\nn \tilde{\eta} \nn_\alpha^3),$$
so that
$$|\langle(\LL_\mathrm{nl}^\prime(\tilde{\eta}))_x,\tilde{\eta}_{xxx}\rangle_0|
\leq \|(\LL_\mathrm{nl}^\prime(\tilde{\eta}))_x\|_0\|\tilde{\eta}_{xxx}\|_0 \leq c \mu^{6\alpha}(\mu^{-\frac{\alpha}{2}+\frac{1}{2}}\nn \tilde{\eta} \nn_\alpha^2+ \nn \tilde{\eta} \nn^4).$$

$\quad$The estimate
$$|\langle(\LL_\mathrm{nl}^\prime(\tilde{\eta}))_z,\tilde{\eta}_{zzz}\rangle_0|
\leq \|(\LL_\mathrm{nl}^\prime(\tilde{\eta}))_z\|_0\|\tilde{\eta}_{zzz}\|_0 \leq c \mu^{6\alpha}(\mu^{-\frac{\alpha}{2}+\frac{1}{2}}\nn \tilde{\eta} \nn_\alpha^2+ \nn \tilde{\eta} \nn^4)$$
is obtained in the same fashion.

\item
Combing the estimates
\begin{eqnarray*}
\int_{{\mathbb R}^2} \frac{k_1^2}{|k|^2}|k|\coth|k|\frac{k_2^2}{|k|^2}|\hat{\tilde{\eta}}|^2\dk
& \leq & \int_{{\mathbb R}^2} \left(1+\frac{|k|^2}{3}\right)\frac{k_2^2}{|k|^2} |\hat{\tilde{\eta}}|^2\dk \leq c \mu^{4\alpha} \\
& \leq & \int_{{\mathbb R}^2} \left(\frac{k_2^2}{|k|^2}+\frac{|k|^2}{3}\right) |\hat{\tilde{\eta}}|^2\dk \\
& \leq & c \mu^{2\alpha} \nn \tilde{\eta} \nn_\alpha^2
\end{eqnarray*}
and
$$\int_{{\mathbb R}^2} \frac{k_1^2}{|k|^2}|k|\coth|k|(k_1^4+k_2^4)|\hat{\tilde{\eta}}|^2\dk
\leq \int_{{\mathbb R}^2} \left(1+\frac{|k|^2}{3}\right)|k|^4 |\hat{\tilde{\eta}}|^2\dk \leq c \mu^{4\alpha} \nn \tilde{\eta} \nn_\alpha^2$$
with Corollary \ref{Speed estimate}, one finds that
\begin{eqnarray*}
& & \hspace{1cm} \left(\frac{\mu}{\LL(\tilde{\eta})}-1\right)\!\!\!\left(\frac{\mu}{\LL(\tilde{\eta})}+1\right)\!\!\int_{{\mathbb R}^2} \frac{k_1^2}{|k|^2}|k|\coth|k|\frac{k_2^2}{|k|^2}|\hat{\tilde{\eta}}|^2\dk\hspace{1cm} \\[1mm]
& & \hspace{2.5cm}\leq c\mu^{4\alpha}(\mu^{N-\frac{1}{2}-2\alpha}\nn \tilde{\eta} \nn_\alpha^2 + \mu^{-\frac{1}{2}-\frac{\alpha}{2}}
\nn \tilde{\eta} \nn_\alpha^4)\hspace{6cm}\hphantom{\Box}
\end{eqnarray*}
and
\begin{eqnarray*}
& & \hspace{1cm}\left(\frac{\mu}{\LL(\tilde{\eta})}-1\right)\!\!\!\left(\frac{\mu}{\LL(\tilde{\eta})}+1\right)\!\!\int_{{\mathbb R}^2} \frac{k_1^2}{|k|^2}(k_1^4+k_2^4)|k|\coth|k||\hat{\tilde{\eta}}|^2\dk\hspace{1cm} \\[1mm]
& & \hspace{2.5cm}\leq c\mu^{6\alpha}(\mu^{N-\frac{1}{2}-2\alpha}\nn \tilde{\eta} \nn_\alpha^2 + \mu^{-\frac{1}{2}-\frac{\alpha}{2}}
\nn \tilde{\eta} \nn_\alpha^4).\hspace{5cm}\Box
\end{eqnarray*}
\end{itemize}

Estimating the left-hand sides of \eqn{Scaled norm 1}, \eqn{Scaled norm 2} from below and using
Lemma \ref{RHS of scaled norm inequalities} to estimate their right-hand sides from above, one finds that
$$\mu^{-4\alpha}\int_{{\mathbb R}^2} \frac{k_2^4}{|k|^4} |\hat{\tilde{\eta}}|^2 \dk \leq c(\mu+\mu^{-\frac{\alpha}{2}+\frac{1}{2}}\nn \tilde{\eta} \nn_\alpha^2
+\mu^{-\frac{\alpha}{2}-\frac{1}{2}}\nn \tilde{\eta} \nn_\alpha^4),$$
$$\mu^{-6\alpha}\int_{{\mathbb R}^2} |k|^6 |\hat{\tilde{\eta}}|^2 \dk
\leq c(\mu+\mu^{-\frac{\alpha}{2}+\frac{1}{2}}\nn \tilde{\eta} \nn_\alpha^2
+\mu^{-\frac{\alpha}{2}-\frac{1}{2}}\nn \tilde{\eta} \nn_\alpha^4),$$
and adding these inequalities to $\|\tilde{\eta}\|_0^2 \leq c\mu$ yields
\begin{equation}
\nn \tilde{\eta} \nn_\alpha^2 \leq c(\mu+\mu^{-\frac{\alpha}{2}+\frac{1}{2}}\nn \tilde{\eta} \nn_\alpha^2
+\mu^{-\frac{\alpha}{2}-\frac{1}{2}}\nn \tilde{\eta} \nn_\alpha^4),
\label{Preliminary scaled norm estimate}
\end{equation}
from which we deduce our final estimate for $\nn \tilde{\eta} \nn_\alpha$.

\begin{theorem} \label{Property of special minimising sequence}
The inequality $\nn \tilde{\eta} \nn_\alpha^2 \leq c \mu$ holds for each $\alpha<1$.
\end{theorem}
{\bf Proof.} Define $Q=\{\alpha \in [0,1): \nn \tilde{\eta} \nn_\alpha^2 =O(\mu)\}$. The inequality $\nn \tilde{\eta} \nn_{\alpha_1}^2
\leq \nn \tilde{\eta} \nn_{\alpha_2}^2$ for $\alpha_1 \leq \alpha_2$ shows that $[0,\alpha] \subset Q$
whenever $\alpha \in Q$; furthermore $0 \subseteq Q$ because $\nn \tilde{\eta} \nn_0^2 \leq c\| \tilde{\eta} \|_3^2 = O(\mu)$.
Suppose that $\alpha^\star:=\sup Q$ is strictly less than unity and choose $\varepsilon>0$ so that
$\alpha^\star+49\varepsilon<1$. Writing \eqn{Preliminary scaled norm estimate} in the form
$$
\frac{\nn \tilde{\eta} \nn_\alpha^2}{\mu} \leq c\left(1+\mu^{\frac{1}{2}-\frac{\alpha}{2}}\left(\frac{\nn \tilde{\eta} \nn_\alpha^2}{\mu}\right)
+\mu^{\frac{1}{2}-\frac{\alpha}{2}}\left(\frac{\nn \tilde{\eta} \nn_\alpha^2}{\mu}\right)^{\!\!2}\right)
$$
with $\alpha=\alpha^\star+\varepsilon$ and estimating
$$
\nn \tilde{\eta} \nn_{\alpha^\star+\varepsilon}^2 \leq \mu^{-12\varepsilon}\nn \tilde{\eta} \nn_{\alpha^\star-\varepsilon}^2,
$$
one finds that
$$
\frac{\nn \tilde{\eta} \nn_{\alpha^\star+\varepsilon}^2}{\mu}\ \leq\ c\Bigg(1+\mu^{\frac{1}{2}-\frac{\alpha^\star}{2}-\frac{25\varepsilon}{2}}
\underbrace{\left(\frac{\nn \tilde{\eta} \nn_{\alpha^\star-\varepsilon}^2}{\mu}\right)}_{\displaystyle = O(1)}
+\mu^{\frac{1}{2}-\frac{\alpha^\star}{2}-\frac{49\varepsilon}{2}}
\underbrace{\left(\frac{\nn \tilde{\eta} \nn_{\alpha^\star-\varepsilon}^2}{\mu}\right)^{\!\!2}}_{\displaystyle = O(1)}\Bigg)
\ \leq\ c,
$$
which leads to the contradiction that $\alpha^\star+\varepsilon \in Q$. It follows that $\alpha^\star=1$ and
$\nn \tilde{\eta} \nn_\alpha^2 = O(\mu)$ for each $\alpha<1$.\qed

It remains to confirm the discussed property of the mapping \eqn{Decreasing function of a}
and to deduce the strict sub-additivity inequality \eqn{Strict SH}. The following preliminary
result is used in the proof.

\begin{proposition} \label{SH step 0}
The quantities $\MM_\mu(\tilde{\eta})$ and $\tilde{\MM}_\mu(\tilde{\eta})$ satisfy the estimates
\begin{eqnarray*}
\MM_\mu(\tilde{\eta}) & = & -\left(\frac{\mu}{\LL_2(\tilde{\eta})}\right)^{\!\!2}\LL_3(\tilde{\eta})+O(\|\tilde{\eta}\|_Z^2\|\tilde{\eta}\|_3^2), \\
\langle \MM_\mu^\prime(\tilde{\eta}),\tilde{\eta}\rangle_0
& = & \left(\frac{\mu}{\LL_2(\tilde{\eta})}\right)^{\!\!2}\LL_3(\tilde{\eta})+O(\|\eta\|_Z^2\|\eta\|_3^2), \\
\tilde{\MM}_\mu(\tilde{\eta}) & = & -\mu^{-1}\left(\frac{\mu}{\LL_2(\tilde{\eta})}\right)^{\!\!2}\LL_3(\tilde{\eta})+O(\mu^{-1}\|\tilde{\eta}\|_Z^2\|\tilde{\eta}\|_3^2).
\end{eqnarray*}
\end{proposition}
{\bf Proof.} Remark \ref{Estimates for KKnl} and Proposition \ref{Z estimate for LL} imply that
$\KK_\mathrm{nl}(\tilde{\eta})$, $\langle\KK_\mathrm{nl}^\prime(\tilde{\eta}),\tilde{\eta}\rangle_0$,
$\LL_\mathrm{nl}(\tilde{\eta})-\LL_3(\tilde{\eta})$ and
$\langle\LL_\mathrm{nl}^\prime(\tilde{\eta})-\LL_3^\prime(\tilde{\eta}),\tilde{\eta}\rangle_0$
are all $O(\|\tilde{\eta}\|_Z^2\|\tilde{\eta}\|_3^2)$,  while
Proposition \ref{Z estimate for LL} shows that
$\LL_\mathrm{nl}(\tilde{\eta}) =O(\|\tilde{\eta}\|_Z \|\tilde{\eta}\|_3^2)$
and
$$
|\langle\LL_\mathrm{nl}^\prime(\tilde{\eta}),\tilde{\eta}\rangle_0|
\ \leq\ 3|\LL_3(\tilde{\eta})|+|\langle\LL_\mathrm{nl}^\prime(\tilde{\eta})-\LL_3^\prime(\tilde{\eta}), \tilde{\eta}\rangle_0|
\ \leq\ c\|\tilde{\eta}\|_Z \|\tilde{\eta}\|_3^2.
$$
The assertions follow by estimating the right-hand sides of the identities \eqn{lot in MM},
\begin{eqnarray*}
\lefteqn{\langle \MM_\mu^\prime(\eta),\eta \rangle_0 }\\
& & = \langle \KK_\mathrm{nl}^\prime(\eta),\eta \rangle + \frac{\mu^2\LL_3(\eta)}{\LL_2(\eta)^2} + \frac{4\mu^2}{\LL_2(\eta)^2}(\LL_\mathrm{nl}(\eta)-\LL_3(\eta))
-\frac{\mu^2}{\LL_2(\eta)^2}\langle\LL_\mathrm{nl}^\prime(\eta) -\LL_3^\prime(\eta),\eta \rangle_0 \\
& & \qquad \mbox{}-\frac{4\mu^2\LL_\mathrm{nl}(\eta)^2}{\LL_2(\eta)^2\LL(\eta)}
-\frac{2\mu^2\LL_\mathrm{nl}(\eta)^2}{\LL_2(\eta)\LL(\eta)^2}
+ \frac{\mu^2\LL_\mathrm{nl}(\eta) \langle \LL_\mathrm{nl}^\prime(\eta),\eta\rangle_0}{\LL_2(\eta)^2\LL(\eta)}
+ \frac{\mu^2\LL_\mathrm{nl}(\eta) \langle \LL_\mathrm{nl}^\prime(\eta),\eta\rangle_0}{\LL_2(\eta)\LL(\eta)^2}
\end{eqnarray*}
and
$$\tilde{\MM}_\mu(\eta) = - \frac{\mu\LL_3(\eta)}{\LL_2(\eta)^2} - \frac{\mu}{\LL_2(\eta)^2}(\LL_\mathrm{nl}(\eta)-\LL_3(\eta))
+\frac{\mu \LL_\mathrm{nl}(\eta)^2}{\LL_2(\eta)^2\LL(\eta)}$$
using these rules.\qed

\begin{proposition} \label{SH step 1}
The function
$$a \mapsto a^{-\frac{5}{2}}\MM_{a^2\mu}(a\tilde{\eta}_n), \qquad a \in [1,2]$$
is decreasing and strictly negative.
\end{proposition}
{\bf Proof.} Because $\{\tilde{\eta}_n\}$ is a minimising sequence for $\JJ_\mu$ over $U \sm \{0\}$ with $\lim_{n \rightarrow \infty}
\|\JJ_\mu^\prime(\tilde{\eta}_n)\|_1 =0$ it has the properties
\eqn{Properties of tildeeta} (see Remark \ref{General properties without rho}) and we may assume that
$\|\JJ_\mu^\prime(\tilde{\eta}_n)\|_1 \leq \mu^{6\alpha+\frac{1}{2}}$, so that $\nn \tilde{\eta}_n \nn_\alpha^2 = O(\mu)$.
It follows that
\begin{eqnarray*}
\frac{\mathrm{d}}{\mathrm{d}a}\left(a^{-\frac{5}{2}}\MM_{a^2\mu}(a\tilde{\eta}_n)\right)
& = & a^{-\frac{7}{2}}\left(
{\textstyle -\frac{5}{2}}\MM_{a^2\mu}(a\tilde{\eta}_n)+\langle \tilde{\MM}^\prime_{a^2\mu}(a\tilde{\eta}_n),a\tilde{\eta}_n\rangle_0
+ 4a^2\mu\tilde{\MM}_{a^2\mu}(a\tilde{\eta}_n)\right) \nonumber \\
& = & \frac{1}{2}a^{-\frac{1}{2}}\Bigg[-\!\!\left(\frac{\mu}{\LL_2(\tilde{\eta}_n)}\right)^{\!\!2}\LL_3(\tilde{\eta}_n)+a\hspace{-7mm}\underbrace{O(\|\tilde{\eta}_n\|_Z^2\|\tilde{\eta}_n\|_3^2)}_{\displaystyle\begin{array}{c}=O(\mu^{3\alpha}\nn \tilde{\eta}_n \nn_\alpha^2 \|\tilde{\eta}_n\|_3^2)\\[1mm]
=O(\mu^{3\alpha+2})\\[1mm] = O(\mu^4)\end{array}}\hspace{-8mm}\Bigg], \\
& = & \frac{1}{2}a^{-\frac{1}{2}}\Big[ \MM_\mu(\tilde{\eta}_n) +\underbrace{O(\|\tilde{\eta}_n\|_Z^2\|\tilde{\eta}_n\|_3^2)}_{\displaystyle=O(\mu^4)} \Big] \\[1mm]
& <  & 0
\end{eqnarray*}
in which Proposition \ref{SH step 0} and the fact that $\MM_\mu(\tilde{\eta}_n) \leq -c\mu^3$ (see Remark \ref{General properties without rho}) have been used. We conclude that
$$a^{-\frac{5}{2}}\MM_{a^2\mu}(a\tilde{\eta}_n)\ \leq\ \MM_\mu(\tilde{\eta}_n)\ <\ 0,
\qquad a \in [1,2].\eqno{\Box}$$

\begin{lemma} \label{SH step 2}
The strict sub-homogeneity property
$$c_{a\mu} < ac_\mu$$
holds for each $a>1$.
\end{lemma}
{\bf Proof.} It suffices to establish this result for $a \in (1,4]$ (see Buffoni \cite[p.\ 56]{Buffoni04a}).

Replacing $a$ by $a^{\frac{1}{2}}$, we find from the above proposition that
$$\MM_{a\mu}(a^{\frac{1}{2}}\tilde{\eta}_n)\ \leq a^{\frac{5}{4}}\MM_\mu(\tilde{\eta}_n), \qquad a \in (1,4]$$
and therefore that
\begin{eqnarray*}
c_{a\mu} & \leq & \KK(a^{\frac{1}{2}}\tilde{\eta}_n) + \frac{a^2\mu^2}{\LL(a^{\frac{1}{2}}\tilde{\eta}_n)}  \\
& = & \KK_2(a^{\frac{1}{2}}\tilde{\eta}_n) + \frac{a^2\mu^2}{\LL_2(a^{\frac{1}{2}}\tilde{\eta}_n)}
+\MM_{a\mu}(a^{\frac{1}{2}}\tilde{\eta}_n) \\
& \leq & a\left(\KK_2(\tilde{\eta}_n)+\frac{\mu^2}{\LL_2(\tilde{\eta}_n)}\right) + a^{\frac{5}{4}}\MM_\mu(\tilde{\eta}_n) \\
& = & a\left(\KK_2(\tilde{\eta}_n)+\frac{\mu^2}{\LL_2(\tilde{\eta}_n)}+\MM_\mu(\tilde{\eta}_n)\right)
+ (a^{\frac{5}{4}}-a)\MM_\mu(\tilde{\eta}_n) \\
& \leq & \JJ_\mu(\tilde{\eta}_n) - c(a^{\frac{5}{4}}-a)\mu^3
\end{eqnarray*}
for $a \in (1,4]$. In the limit $n \rightarrow \infty$ the above inequality yields
$$c_{a\mu}\ \leq\ ac_\mu - c(a^{\frac{5}{4}}-a)\mu^3\ < ac_\mu.\eqno{\Box}$$

\section{Conditional energetic stability} \label{Stability}

The following theorem, which is proved using the results of Section \ref{MS},
and \ref{SSA}, is our final result concerning the set of minimisers
of $\JJ_\mu$ over $U\sm\{0\}$. 

\begin{theorem} \label{Key minimisation theorem}
\hspace{1in}
\begin{list}{(\roman{count})}{\usecounter{count}}
\item
The set $C_\mu$ of minimisers of $\JJ_\mu$ over $U \sm \{0\}$ is non-empty.
\item
Suppose that $\{\eta_n\}$ is a minimising sequence for $\JJ_\mu$ on $U\sm\{0\}$ which satisfies
\begin{equation}
\sup_{n\in{\mathbb N}} \|\eta_n\|_3 < M. \label{Sup less than M}
\end{equation}
There exists a sequence $\{(x_n,z_n)\} \subset {\mathbb R}^2$ with the property that
a subsequence of\linebreak $\{\eta_n(x_n+\cdot,z_n+\cdot)\}$ converges
in $H^r({\mathbb R}^2)$, $0 \leq r < 3$ to a function $\eta \in C_\mu$.
\end{list}
\end{theorem}
{\bf Proof.} It suffices to prove part (ii), since an application of this result to the sequence
$\{\tilde{\eta}\}$ constructed in Section \ref{Special minimising sequence} above yields part (i).

In order to establish part (ii) we choose $\tilde{M} \in (\sup \|\eta_n\|_3, M)$, so that
$\{\eta_n\}$ is also a minimising sequence for the functional $\JJ_{\rho,\mu}$
discussed in Section \ref{MS} (the existence of a minimising sequence $\{v_n\}$
for $\JJ_{\rho,\mu}$ with $\lim_{n \rightarrow \infty} \JJ_{\rho,\mu}(v_n) <
\lim_{n \rightarrow \infty} \JJ_{\rho,\mu}(\eta_n)$ would lead to the contradiction
$$\lim_{n \rightarrow \infty} \JJ_\mu(v_n)\ \leq\ \lim_{n \rightarrow \infty} \JJ_{\rho,\mu}(v_n)
\ <\  \lim_{n \rightarrow \infty} \JJ_{\rho,\mu}(\eta_n)
\ =\ \lim_{n \rightarrow \infty} \JJ_\mu(\eta_n)\ =\ c_\mu.\ )$$
We may therefore study $\{\eta_n\}$ using
the theory given there, noting that  the sequence $\{u_n\}$
defined in equation \eqn{Definition of un} does not have the `dichotomy' property:
the existence of sequences $\{\eta_n^{(1)}\}$, $\{\eta_n^{(2)}\}$ with the features
listed in Lemma \ref{Splitting properties 1} and Corollary \ref{J splitting} is
incompatible with the  strict sub-additivity property \eqn{Strict SH} of $c_\mu$.
Recall that the positive numbers $\mu^{(1)}$, $\mu^{(2)}$ defined in
Lemma \ref{J splitting} sum to $\mu$;  this fact leads to the contradiction
\begin{eqnarray*}
c_\mu & < & c_{\mu^{(1)}} + c_{\mu^{(2)}} \\
& \leq & \lim_{n \rightarrow \infty}\JJ_{\mu^{(1)}}(\eta_n^{(1)}) + \lim_{n \rightarrow \infty}\JJ_{\mu^{(2)}}(\eta_n^{(2)}) \\
& = & \lim_{n \rightarrow \infty}\JJ_\mu(\eta_n) \\
& = & c_\mu,
\end{eqnarray*}
where Corollary \ref{J splitting} has been used.
We conclude that $\{u_n\}$ has the `concentration' property and hence
$\eta_n(\cdot+x_n,\cdot+z_n) \rightarrow \eta$ in
$H^r({\mathbb R}^2)$ for every $r \in [0,3)$ (see the proof of Lemma \ref{Concentration}), whereby
$\JJ_\mu(\eta)=\lim_{n \rightarrow \infty} \JJ_\mu(\eta_n(\cdot+x_n,\cdot+z_n)) = c_\mu$, so that $\eta$
is a minimiser of $\JJ_\mu$ over $U\sm\{0\}$.\qed

The next step is to relate the above result to our original problem
finding minimisers of $\EE(\eta,\Phi)$ subject to the constraint $\II(\eta,\Phi) = 2\mu$,
where $\EE$ and $\II$ are defined in equations \eqn{Definition of E} and
\eqn{Definition of I}.

\begin{theorem} \label{Result for constrained minimisation} \hspace{1cm}
\begin{list}{(\roman{count})}{\usecounter{count}}
\item
The set $D_\mu$ of minimisers of $\EE$ on the set
$$S_\mu=\{(\eta,\Phi) \in U \times H_\star^{1/2}({\mathbb R}^2): \II(\eta,\Phi)=2\mu\}$$
is non-empty.
\item
Suppose that $\{(\eta_n,\Phi_n)\} \subset S_\mu$ is a minimising sequence for $\EE$ with
the property that
$$
\sup_{k\in{\mathbb N}} \|\eta_n\|_3 < M.
$$
There exists a sequence $\{(x_n,z_n)\} \subset {\mathbb R}^2$ with the property that
a subsequence of\linebreak $\{\eta_n(x_n+\cdot,z_n+\cdot),\Phi_n(x_n+\cdot,z_n+\cdot)\}$
converges  in
$H^r({\mathbb R}^2) \times H_\star^{1/2}({\mathbb R}^2)$, $0 \leq r < 3$ to a function in $D_\mu$.
\end{list}
\end{theorem}
{\bf Proof.} (i) We consider the minimisation problem in two steps.

\emph{1. Fix $\eta \in U\sm\{0\}$
and minimise $\EE(\eta,\cdot)$ over $T_\mu:=\{\Phi \in H_\star^{1/2}({\mathbb R}^2): \II(\eta,\Phi)=2\mu\}$.}
Let $\{\Phi_n\} \subset H_\star^{1/2}({\mathbb R}^2)$ be a minimising sequence for $\EE(\eta,\cdot)$
on $T_\mu$. The sequence is clearly bounded (because $\|\Phi_n\|_{\star,1/2} \rightarrow \infty$ as $n \rightarrow
\infty$ would imply that $\EE(\eta,\Phi_n) \rightarrow \infty$ as $n \rightarrow \infty$ and contradict the
fact that $\{\Phi_n\}$ is a minimising sequence for $\EE(\eta,\cdot)$ on $T_\mu$) and hence admits a
weakly convergent subsequence (still denoted by $\{\Phi_n\}$).
Notice that $\EE(\eta,\cdot)$ is weakly lower semicontinuous on $H_\star^{1/2}({\mathbb R}^2)$
($\Phi \mapsto \left(\int_{{\mathbb R}^2} \Phi G(\eta) \Phi \dx\dz\right)^{\!\frac{1}{2}}$ is equivalent to its usual norm)
and $\II(\eta,\cdot)$ is weakly continuous on $H_\star^{1/2}({\mathbb R}^2)$; a familiar argument
shows that the sequence $\{\Phi_n\}$ converges to a minimiser $\Phi_\eta$ of $\EE(\eta,\cdot)$ on
$T_\mu$.

\emph{2. Minimise $\EE(\eta,\Phi_\eta)$ over $U\sm\{0\}$.}
Because $\Phi_\eta$ minimises $\EE(\eta,\cdot)$ over $T_\mu$ there exists a
Lagrange multiplier $\lambda_\eta$ such that $$G(\eta)\Phi_\eta = \lambda_\eta \eta_x,$$
and a straightforward calculation shows that
$\Phi_\eta = \lambda_\eta G(\eta)^{-1}\eta_x$, $\lambda_\eta = \mu/\LL(\eta)$
(which also confirms the uniqueness of $\Phi_\eta$). According to Theorem \ref{Key minimisation theorem}(i)
the set $C_\mu$ of minimisers of $\JJ_\mu(\eta):=\EE(\eta,\Phi_\eta)$ over $U \sm \{0\}$ is not empty; it follows
that $D_\mu$ is also not empty.

(ii) Let $\{(\eta_n,\Phi_n)\} \subset U \times H_\star^{1/2}({\mathbb R}^2)$ be a minimising sequence for
$\EE$ over $S_\mu$ with $\sup \|\eta_n\|_3 < M$. The inequality
$$\EE(\eta_n,\Phi_{\eta_n}) \leq \EE(\eta_n,\Phi_n)$$
implies that $\{(\eta_n,\Phi_{\eta_n})\} \subset U \times H_\star^{1/2}({\mathbb R}^2)$ is also
a minimising sequence; it follows that $\{\eta_n\} \subset U\sm\{0\}$ is a minimising sequence for
$\JJ_\mu$ which therefore converges (up to translations and
subsequences) in $H^r({\mathbb R}^2)$, $0 \leq r < 3$
to a minimiser $\eta$ of $\JJ_\mu$ over $U\sm\{0\}$ (see Theorem \ref{Key minimisation theorem}(ii)).

The relations
$$\Phi_{\eta_n} = \frac{\mu G^{-1}(\eta_n)\eta_{nx}}{\LL(\eta_n)}, \qquad
\Phi_\eta = \frac{\mu G^{-1}(\eta)\eta_{x}}{\LL(\eta)}$$
show that $\Phi_{\eta_n} \rightarrow \Phi_\eta$ in $H_\star^{1/2}({\mathbb R}^2)$, and using
this result and the calculation
\begin{eqnarray*}
\lefteqn{c \|\Phi_n - \Phi_{\eta_n}\|_{\star,1/2}^2} \qquad \\
& \leq & \int_{{\mathbb R}^2} (\Phi_n-\Phi_{\eta_n}) G(\eta_n)(\Phi_n-\Phi_{\eta_n}) \dx\dz \\
& = & \int_{{\mathbb R}^2} \Phi_n G(\eta_n) \Phi_n \dx\dz + \int_{{\mathbb R}^2} \Phi_{\eta_n} G(\eta_n) \Phi_{\eta_n} \dx\dz - 2 \int_{{\mathbb R}^2} \Phi_n G(\eta_n) \Phi_{\eta_n} \dx\dz \\
& = & \int_{{\mathbb R}^2} \Phi_n G(\eta_n) \Phi_n \dx\dz + \int_{{\mathbb R}^2} \Phi_{\eta_n} G(\eta_n) \Phi_{\eta_n} \dx\dz -4\lambda_{\eta_n}\mu \\
& = & \int_{{\mathbb R}^2} \Phi_n G(\eta_n) \Phi_n \dx\dz - \int_{{\mathbb R}^2} \Phi_{\eta_n} G(\eta_n) \Phi_{\eta_n} \dx\dz \\
& = & 2 \EE(\eta_n,\Phi_n) - 2\EE(\eta_n, \Phi_{\eta_n}) \\
& \rightarrow & 0
\end{eqnarray*}
as $n \rightarrow \infty$, one finds that $\Phi_n \rightarrow \Phi_\eta$ in $H_\star^{1/2}({\mathbb R}^2)$
as $n \rightarrow \infty$.\qed

It is also possible to obtain a bound on the speed of the waves described by functions in $D_\mu$.

\begin{lemma} \label{ND speed of wave}
The fully localised solitary wave corresponding to $(\tilde{\eta},\tilde{\Phi})\in D_\mu$ is subcritical, that is its
dimensionless speed is less than unity.
\end{lemma}
{\bf Proof.} The dimensionless speed of the wave is the Lagrange multiplier $\lambda$ in
the equations
$$\mathrm{d}\EE[\tilde{\eta},\tilde{\Phi}] = \lambda\mathrm{d} \II[\tilde{\eta},\tilde{\Phi}]$$
satisfied by by the constrained minimiser $(\tilde{\eta},\tilde{\Phi})$. Examining the second component of
this equation, one finds that
$$\lambda = \frac{\mu}{\LL(\tilde{\eta})}$$
(see the calculation in the proof of Theorem \ref{Result for constrained minimisation}(i)). 

Because $\tilde{\eta}$ minimises $\JJ_\mu$ over $U \sm \{0\}$ it has the properties \eqn{Properties of tildeeta} with
$\JJ_\mu^\prime(\tilde{\eta})=0$. It follows from Proposition \ref{A/B estimate on speed} that
$$
\frac{\mu}{\LL(\tilde{\eta})} -1 \leq
\frac{\langle \MM_\mu^\prime(\tilde{\eta}),\tilde{\eta}\rangle_0}{4\mu}
-\frac{\MM_\mu(\tilde{\eta})}{2\mu} + \tilde{\MM}_\mu(\tilde{\eta}),
$$
and using Proposition \ref{SH step 0}, one finds that
\begin{eqnarray*}
\frac{\langle \MM_\mu^\prime(\tilde{\eta}),\tilde{\eta}\rangle_0}{4\mu}-\frac{\MM_\mu(\tilde{\eta})}{2\mu} + \tilde{\MM}_\mu(\tilde{\eta})
& = &
-\frac{1}{4\mu}\left(\frac{\mu}{\LL_2(\tilde{\eta})}\right)^{\!\!2}\LL_3(\tilde{\eta})+\underbrace{O(\mu^{-1}\|\tilde{\eta}\|_Z^2\|\tilde{\eta}\|_3^2)}_{\displaystyle\begin{array}{c}=O(\mu^{3\alpha}\nn \tilde{\eta} \nn_\alpha^2)\\[1mm]
=O(\mu^{3\alpha+1})\\[1mm] = O(\mu^3)\end{array}} \\
& = & \frac{1}{4\mu}\Big[ \MM_\mu(\eta) +\underbrace{O(\|\tilde{\eta}\|_Z^2\|\tilde{\eta}\|_3^2)}_{\displaystyle=O(\mu^4)} \Big] + O(\mu^3) \\
& \leq & -c\mu^2 + O(\mu^3)
\end{eqnarray*}
because $\MM_\mu(\tilde{\eta}) \leq -c\mu^3$ (see Remark \ref{General properties without rho}).\qed

Our stability result (Theorem \ref{CES} below)
is obtained from Theorem \ref{Result for constrained minimisation}
under the following assumption concerning the
well-posedness of the hydrodynamic problem with small initial data.\\
\\
{\bf (Well-posedness assumption)} There exists a subset  $S$
of $U \times H_\star^{1/2}({\mathbb R}^2)$ with the following properties.
\begin{list}{(\roman{count})}{\usecounter{count}}
\item
The closure of $S \sm D_\mu$ in $L^2({\mathbb R}^2)$ has a non-empty intersection with
$D_\mu$.
\item
For each $(\eta_0,\Phi_0) \in S$
there exists $T>0$ and a continuous function $t \mapsto (\eta(t), \Phi(t)) \in 
U \times H_\star^{1/2}({\mathbb R}^2)$, $t \in [0, T]$ such that $(\eta(0),\Phi(0)) = (\eta_0,\Phi_0)$,
$$\EE(\eta(t),\Phi(t)) = \EE(\eta_0,\Phi_0),\ \II(\eta(t),\Phi(t))=\II(\eta_0,\Phi_0), \qquad t \in [0,T]$$
and 
$$\sup_{t \in [0,T]} \|\eta(t)\|_3 < M.$$
\end{list}

\begin{theorem} \label{CES}
Choose $r \in [0,3)$. For each $\varepsilon>0$  there exists $\delta>0$ such that
$$(\eta_0,\Phi_0) \in S,\ \dist((\eta_0,\Phi_0), D_\mu) < \delta \quad \Rightarrow \quad
\dist((\eta(t),\Phi(t)), D_\mu)<\varepsilon,$$
for $t\in[0,T]$,
where `$\dist$' denotes the distance in $H^r({\mathbb R}^2) \times H_\star^{1/2}({\mathbb R}^2)$.
\end{theorem}
{\bf Proof.} This result is proved by contradiction. Suppose the assertion is false: there exists a
real number $\varepsilon>0$ and sequences $\{(\eta_{0,n},\Phi_{0,n})\} \subset
U \times H_\star^{1/2}({\mathbb R}^2)$, $\{T_n\} \subset (0,\infty)$,
and $\{(\eta_n(\cdot),\Phi_n(\cdot))\} \in C([0,T_n], U \times H_\star^{1/2}({\mathbb R}^2))$ such that
\begin{eqnarray}
(\eta_n(0),\Phi_n(0)) & = & (\eta_{0,n},\Phi_{0,n}), \nonumber \\
\EE(\eta_n(t),\Phi_n(t)) & = & \EE(\eta_{0,n},\Phi_{0,n}), \qquad t \in [0,T_n], \nonumber \\
\II(\eta_n(t),\Phi_n(t)) & = & \II(\eta_{0,n},\Phi_{0,n}), \qquad t \in [0,T_n], \nonumber \\
\dist ((\eta_{0,n},\Phi_{0,n}),D_\mu) & < & \frac{1}{n} \label{Property 4}
\end{eqnarray}
together with a sequence $\{t_n\} \subset {\mathbb R}$ with $t_n \in [0,T_n]$ and
\begin{equation}
\dist((\eta_n(t_n),\Phi_n(t_n)),D_\mu) \geq \varepsilon. \label{Property 5}
\end{equation}

Inequality \eqn{Property 4} asserts the existence of
a sequence $\{(\tilde{\eta}_n, \tilde{\Phi}_{\tilde{\eta}_n})\}$
in $D_\mu$ such that
\begin{equation}
\|\eta_{0,n}-\tilde{\eta}_n\|_r \rightarrow 0, \label{Compare two sequences 1}
\end{equation}
\begin{equation}
\|\Phi_{0,n}-\tilde{\Phi}_{\tilde{\eta}_n}\|_{\star,1/2} \rightarrow 0
 \label{Compare two sequences 2}
\end{equation}
as $n \rightarrow \infty$.
Observe that $\tilde{\eta}_n$ minimises $\JJ_\mu$ over $U \sm \{0\}$ and satisfies
$\sup \|\tilde{\eta}_n\|_3 < M$; the sequence $\{\tilde{\eta}_n\}$ therefore converges
(up to subsequences and translations) in $H^r({\mathbb R}^2)$ to a
minimiser $\eta$ of $\JJ_\mu$ over $U \sm \{0\}$, and \eqn{Compare two sequences 1}
shows the same is true of $\{\eta_{0,n}\}$. Furthermore, the relations
$$\Phi_{\tilde{\eta}_n} = \frac{\mu G^{-1}(\tilde{\eta}_n)\tilde{\eta}_{nx}}{\LL(\tilde{\eta}_n)}, \qquad
\Phi_\eta = \frac{\mu G^{-1}(\eta)\eta_{x}}{\LL(\eta)}$$
show that $\Phi_{\tilde{\eta}_n} \rightarrow \Phi_\eta$ in $H_\star^{1/2}({\mathbb R}^2)$, and combining
this fact with \eqn{Compare two sequences 2}, we find that $\Phi_{0,n} \rightarrow \Phi_\eta$ in
$H_\star^{1/2}({\mathbb R}^2)$ as $n \rightarrow \infty$. Altogether, these arguments show that
$$\EE(\eta_{0,n},\Phi_{0,n}) \rightarrow \EE(\eta,\Phi_\eta),$$
$$2\mu_n:=\int_{{\mathbb R}^2} \partial_x \eta_{0,n} \Phi_{0,n}\dx\dz \rightarrow
\int_{{\mathbb R}^2} \eta_x \Phi_\eta \dx\dz :=2\mu$$
as $n \rightarrow \infty$.

Define $\hat{\eta}_n = \eta_n(t_n)$, $\hat{\Phi}_n = (\mu/\mu_n)\Phi_n(t_n)$, so that
\begin{eqnarray*}
\lim_{n \rightarrow \infty} \EE(\hat{\eta}_n, \hat{\Phi}_n) & = & \lim_{n \rightarrow \infty} \EE(\eta_n(t_n),\Phi_n(t_n)) \\
& = & \lim_{n \rightarrow \infty} \EE(\eta_{0,n},\Phi_{0,n}) \\
& = & \EE(\eta,\Phi_\eta), \\[2mm]
\II(\hat{\eta}_n, \hat{\Phi}_n)
& = & \frac{\mu}{\mu_n} \int_{{\mathbb R}^2} \partial_x \eta_n(t_n)\Phi_n(t_n)\dx\dz \\
& = &  \frac{\mu}{\mu_n} \int_{{\mathbb R}^2} \partial_x \eta_{0,n} \Phi_{0,n} \dx\dz \\
& = & 2\mu;
\end{eqnarray*}
it follows that $\{(\hat{\eta}_n, \hat{\Phi}_n)\}$ is a minimising sequence for $\EE$ over
$S_\mu$ with
$\sup\|\hat{\eta}_n\|_3 < M$, and Theorem \ref{Result for constrained minimisation} therefore asserts that
$\dist ((\hat{\eta}_n,\hat{\Phi}_n), D_\mu) \rightarrow 0$ as $n \rightarrow \infty$.
On the other hand, it follows from \eqn{Property 5} and the limit
$\|(\hat{\eta}_n,\hat{\Phi}_n)-(\eta_n(t_n),\Phi_n(t_n))\|_r \rightarrow 0$ as $n \rightarrow \infty$ that
$$\dist ((\hat{\eta}_n,\hat{\Phi}_n), D_\mu)) \geq \frac{\varepsilon}{2}.\eqno{\Box}$$

\appendix
\section*{Appendix A: Proof of Proposition \ref{Regularity proposition 2}}
\setcounter{section}{1}

This proposition is proved using the observation that
the unique weak solution $u \in H_\star^1(\Sigma)$ to the boundary-value problem
\eqn{BC for un 1}--\eqn{BC for un 3}, where $F_1$, $F_2$, $F_3 \in L^2(\Sigma)$
and we have dropped the superscript $n$ for notational simplicity,
is given by the explicit formula
$$u={\mathcal F}^{-1}\left[\int_0^1 \left\{G(y,\tilde{y})(\i k_1 \hat{F}_1 + \i k_2 \hat{F}_2) - G_{\tilde{y}}(y,\tilde{y})\hat{F}_3 \right\} \dtildey \right],$$
in which
$$G(y,\tilde{y}) = \left\{\begin{array}{cccc}
\displaystyle -\frac{\cosh|k|y \cosh |k|(1-\tilde{y})}{|k|\sinh|k|}, & & & 0 \leq y \leq \tilde{y} \leq 1, \\
\\
\displaystyle -\frac{\cosh|k|\tilde{y} \cosh |k|(1-y)}{|k|\sinh|k|}, & & & 0 \leq \tilde{y} \leq y \leq 1
\end{array}\right.$$

We proceed by writing the formula for $u$ in the alternative form
$$u={\mathcal F}^{-1}\left[\int_0^1 \left\{G(y,\tilde{y})(\i k_1 \hat{F}_1 + \i k_2 \hat{F}_2) + H_y(y,\tilde{y})\hat{F}_3\right\} \dtildey\right],$$
where
$$H(y,\tilde{y}) = \left\{\begin{array}{cccc}
\displaystyle -\frac{\sinh|k|y \sinh |k|(1-\tilde{y})}{|k|\sinh|k|}, & & & 0 \leq y \leq \tilde{y} \leq 1, \\
\\
\displaystyle -\frac{\sinh|k|\tilde{y} \sinh |k|(1-y)}{|k|\sinh|k|}, & & & 0 \leq \tilde{y} \leq y \leq 1;
\end{array}\right.$$
elementary calculations yield the following estimates for $G$ and $H$.

\begin{proposition}
The estimates
$$|G| \leq \frac{c}{|k|^2}, \qquad |H| \leq c, \qquad \left|\left\{\begin{array}{c}\!\!\!\partial_y\!\!\! \\ \!\!\!\partial_{\tilde{y}}\!\!\!\end{array}\right\}\!\left\{
\begin{array}{c} \!\!\!G\!\!\! \\ \!\!\!H\!\!\! \end{array}\right\}\right| \leq c$$
for $|k| \leq \delta$ and
$$
\left|\left\{
\begin{array}{c} \!\!\!G\!\!\! \\ \!\!\!H\!\!\! \end{array}\right\}\right| \leq \frac{c}{|k|}\e^{-|k||\tilde{y}-y|}, \qquad
\left|\left\{\begin{array}{c}\!\!\!\partial_y\!\!\! \\ \!\!\!\partial_{\tilde{y}}\!\!\!\end{array}\right\}\!\left\{
\begin{array}{c} \!\!\!G\!\!\! \\ \!\!\!H\!\!\! \end{array}\right\}\right| \leq c\, \e^{-|k||\tilde{y}-y|}
$$ 
for $|k| \geq \delta$ hold for each sufficiently small positive number $\delta$.
\end{proposition}

\begin{lemma}
Suppose that $F_1$, $F_2$, $F_3 \in H^r(\Sigma)$, $r \geq 0$. The function $u$
satisfies
$$\|u_x\|_r,\ \|u_y\|_r,\ \|u_z\|_r \leq c(\|F_1\|_r + \|F_2\|_r + \|F_3\|_r).$$
\end{lemma}
{\bf Proof.} The representation
$$
u_x = {\mathcal F}^{-1}\left[\int_0^1 \left\{\i k_1G(y,\tilde{y})(\i k_1 \hat{F}_1 + \i k_2 \hat{F}_2) + \i k_1H_y(y,\tilde{y})\hat{F}_3 \right\}\dtildey\right]
$$
yields the formulae
\begin{eqnarray*}
{\mathcal F}[\partial_x^n u_x] & = & \int_0^1 \left\{\i^{n+1} k_1^{n+1}G(y,\tilde{y})(\i k_1 \hat{F}_1 + \i k_2 \hat{F}_2) + \i^{n+1} k_1^{n+1}H_y(y,\tilde{y})\hat{F}_3 \dtildey\right\}, \\
{\mathcal F}[\partial_z^n u_x] & = & \int_0^1 \left\{\i^{n+1} k_1k_2^nG(y,\tilde{y})(\i k_1 \hat{F}_1 + \i k_2 \hat{F}_2) + \i^{n+1} k_1k_2^nH_y(y,\tilde{y})\hat{F}_3 \dtildey\right\}, \qquad n \in {\mathbb N}_0
\end{eqnarray*}
and
\begin{eqnarray*}
\lefteqn{{\mathcal F}[\partial_y^{2n} u_x] =
-\int_0^1|k|^{2n} G(y,\tilde{y})(k_1^2 \hat{F}_1 + k_1k_2 \hat{F}_2)\dtildey
-\sum_{j=1}^n |k|^{2j-2} \partial_y^{2n-2j}(k_1^2  \hat{F}_1 + k_1k_2 \hat{F}_2)} \qquad\\
& & \mbox{} +\int_0^1|k|^{2n} H_y(y,\tilde{y})\i k_1 \hat{F}_3\dtildey
+\sum_{j=1}^n |k|^{2j-2}\partial_y^{2n-2j+1}\i k_1 \hat{F}_3, \hspace{1.3in} \\
\lefteqn{{\mathcal F}[\partial_y^{2n+1} u_x] =
-\int_0^1|k|^{2n} G_y(y,\tilde{y})(k_1^2 \hat{F}_1 + k_1k_2 \hat{F}_2)\dtildey
-\sum_{j=1}^n |k|^{2j-2} \partial_y^{2n-2j+1}(k_1^2  \hat{F}_1 + k_1k_2 \hat{F}_2)} \qquad\\
& & \mbox{} +\int_0^1|k|^{2n+2} H(y,\tilde{y})\i k_1 \hat{F}_3\dtildey
+\sum_{j=1}^{n+1} |k|^{2j-2}\partial_y^{2n-2j+2}\i k_1 \hat{F}_3, \qquad n \in {\mathbb N}_0,
\end{eqnarray*}
which are established by mathematical induction.

Observe that
\begin{eqnarray}
\|\partial_x^n u_x\|_0 & \leq &
\left\|\int_0^1 G(y,\tilde{y})k_1^{n+2} \hat{F}_1 \dtildey \right\|_0 \nonumber \\
& & \qquad\mbox{}
+\left\|\int_0^1 G(y,\tilde{y})k_1^{n+1}k_2 \hat{F}_2 \dtildey \right\|_0
+\left\|\int_0^1 H_y(y,\tilde{y})k_1^{n+1} \hat{F}_3 \dtildey \right\|_0.
\label{Representation of ux}
\end{eqnarray}
Writing
$$\left\|\int_0^1 G(y,\tilde{y})k_1^{n+2} \hat{F}_1 \dtildey \right\|_0^2 = I_1 + I_2,$$
where
\begin{eqnarray*}
I_1 & = & \int_0^1 \int_{|k| < \delta} \left| \int_0^1 G(y,\tilde{y})k_1^{n+2}\hat{F}_1 \dtildey \right|^2 \dk \dy \\
& \leq & c \int_0^1 \int_{|k| < \delta} |\hat{F}_1|^2 \dk\dy \\
& \leq & c \|F_1\|_0^2, \\
\\
I_2 & = & \int_0^1 \int_{|k| > \delta} \left| \int_0^1 G(y,\tilde{y})k_1^{n+2}\hat{F}_1 \dtildey \right|^2 \dk \dy \\
& \leq & c \int_0^1 \int_{|k| > \delta} \left| \int_0^1 |k|^{n+1}\e^{-|k||y-\tilde{y}|}|\hat{F}_1|\dtildey\right|^2 \dk \dy \\
& \leq & c \int_{|k| > \delta} |k|^{2n+2} \int_0^1 \left[\int_0^1\e^{-|k||y-\tilde{y}|}|\hat{F}_1|\dtildey\right]^2 \dy\dk \\
& \leq & c \int_{|k| > \delta} |k|^{2n+2} \int_0^1\Bigg[ \underbrace{\int_0^1 \e^{-|k||y-\tilde{y}|} \dtildey}_{\displaystyle=O(|k|^{-1})} \int_0^1\e^{-|k||y-\tilde{y}|}|\hat{F}_1|^2\dtildey\Bigg] \dy \dk\\
& \leq & c \int_{|k| > \delta} |k|^{2n+1} \int_0^1 \underbrace{\int_0^1\e^{-|k||y-\tilde{y}|}\dy}_{\displaystyle=O(|k|^{-1})}\, |\hat{F}_1|^2\dtildey \dk\\
& \leq & c \int_{|k| > \delta} |k|^{2n} \int_0^1 |\hat{F}_1|^2 \dy\dk \\
& \leq & c \|F_1\|_n^2
\end{eqnarray*}
and estimating the second and third terms on the right-hand side of \eqn{Representation of ux}
in a similar fashion, we find that
$$\|\partial_x^n u_x\|_0^2 \leq c (\|F_1\|_n^2 + \|F_2\|_n^2 + \|F_3\|_n^2), \qquad n \in {\mathbb N}_0.$$
Analogous calculations yield the estimates
$$\|\partial_z^n u_x\|_0^2 \leq c (\|F_1\|_n^2 + \|F_2\|_n^2 + \|F_3\|_n^2) \qquad n \in {\mathbb N}_0$$
and
\begin{eqnarray*}
\|\partial_y^{2n} u_x\|_0^2 & \leq & c (\|F_1\|_{2n}^2 + \|F_2\|_{2n}^2 + \|F_3\|_{2n}^2), \\
\|\partial_y^{2n+1} u_x\|_0^2 & \leq & c (\|F_1\|_{2n+1}^2 + \|F_2\|_{2n+1}^2 + \|F_3\|_{2n+1}^2),
\qquad n \in {\mathbb N}_0.
\end{eqnarray*}

Altogether we have established that
$$\|u_x\|_m \leq c (\|F_1\|_m + \|F_2\|_m + \|F_3\|_m), \qquad m \in {\mathbb N}_0$$
and a similar argument yields the corresponding estimates for $u_y$ and $u_z$. The
advertised result follows by interpolation.\qed

\begin{remark}
Suppose that $r \geq 2$. A straightforward calculation shows that $u$ is a strong solution of
the boundary-value problem \eqn{BC for un 1}--\eqn{BC for un 3}:
equation \eqn{BC for un 1} holds in $H^{r-2}(\Sigma)$ while
equations \eqn{BC for un 2}, \eqn{BC for un 3} hold in $H^{r-3/2}({\mathbb R}^2)$.
\end{remark}

\section*{Appendix B: Proof of Lemma \ref{Test function 1}}
\setcounter{section}{3}

This lemma asserts the existence of $\eta_\mu^\star \in C_0^\infty({\mathbb R}^2)$
with $\JJ_\mu(\eta_\mu^\star)<2\mu-c\mu^3$.
In constructing the `test function' $\eta_\mu^\star$ one is guided by the principle that it should approximate a
minimiser of $\JJ_\mu$. Numerical experiments suggest that the KP-I equation correctly
models fully localised solitary waves in its formal region of validity
(Parau, Vanden-Broeck \& Cooker \cite{ParauVandenBroeckCooker05a});
we therefore expect a minimiser of $\JJ_\mu$ to have the KP length scales (cf.\ equation
\eqn{AS explicit formula}) and use a test function of the form
$$\eta^\star(x,z)=\gamma^2\Psi(\gamma x, \gamma^2 z),\qquad 0<\gamma \ll 1$$
with an appropriate choice of $\Psi \in C_0^\infty([-{\textstyle\frac{1}{2}},{\textstyle\frac{1}{2}}]^2)$ and $\gamma=\gamma(\mu)$. Furthermore, in the following analysis we are confronted with the task of estimating
$$\int_{{\mathbb R}^2} \frac{k_2^2}{|k|^2}|\hat{\eta}^\star|^2 \dk=\gamma \int_{{\mathbb R}^2} \frac{\gamma^4k_2^2}{\gamma^2 k_1^2+\gamma^4 k_2^2}|\hat{\Psi}|^2\dk,$$
which is generally $O(\gamma)$; in the special case
$$\Psi(x,z): = \psi_x(x,z),$$
where $\psi$ also belongs to $C_0^\infty([-{\textstyle\frac{1}{2}},{\textstyle\frac{1}{2}}]^2)$, it
is however $O(\gamma^3)$ because
$$
\int_{{\mathbb R}^2} \frac{\gamma^4k_2^2}{\gamma^2 k_1^2+\gamma^4 k_2^2}|\hat{\Psi}|^2\dk
= \int_{{\mathbb R}^2} \frac{\gamma^2k_1^2k_2^2}{k_1^2+\gamma^2 k_2^2}|\hat{\psi}|^2\dk
= \gamma^2 \int_{{\mathbb R}^2} \frac{k_1^2}{k_1^2+\gamma^2 k_2^2}|\hat{\psi}_z|^2\dk.
$$
It is therefore advantageous to work with this special class of test functions. (Note that the explicit
fully localised solitary-wave solution \eqn{Explicit KP formula} to the KP-I equation
may be written as
$$u(x,z)=\frac{\partial}{\partial x}\left(\frac{-8x}{3+x^2+z^2}\right)$$
and therefore also has this form.)

We begin by computing the leading-order terms in the asymptotic expansions
of $\KK(\eta^\star)$ and $\LL(\eta^\star)$ in powers of $\gamma$. Observe that
$$
\KK_2(\eta^\star)
\ =\ \int_{{\mathbb R}^2}\left\{\frac{1}{2}(\eta^\star)^2 + \frac{\beta}{2}(\eta^\star_x)^2\right\}\dx\dz
\ =\  \frac{\gamma}{2}\int_{{\mathbb R}^2} \Psi^2 \dx\dz
+ \frac{\gamma^3\beta}{2} \int_{{\mathbb R}^2} \Psi_x^2 \dx\dz
$$
and
\begin{eqnarray*}
\LL_2(\eta^\star)
& = & \frac{1}{2}\int_{{\mathbb R}^2} \eta^\star K_0 \eta^\star \dx\dz \\
& = & \frac{1}{2}\int_{{\mathbb R}^2} \frac{k_1^2}{|k|^2}|\hat{\eta}^\star|^2 \dk
+\frac{1}{6}\int_{{\mathbb R}^2} k_1^2 |\hat{\eta}^\star|^2 \dk
+ \frac{1}{2}\int_{{\mathbb R}^2}
\frac{k_1^2}{|k|^2} (|k|\coth|k|\!-\!1\!-\!{\textstyle\frac{1}{3}}|k|^2)|\hat{\eta}^\star|^2 \dk,
\end{eqnarray*}
in which
\begin{eqnarray*}
\int_{{\mathbb R}^2} \frac{k_1^2}{|k|^2}|\hat{\eta}^\star|^2 \dk
& = & \int_{{\mathbb R}^2} |\hat{\eta}^\star|^2\dk - \int_{{\mathbb R}^2} \frac{k_2^2}{|k|^2}
|\hat{\eta}^\star|^2\dk \\
& = & \gamma \int_{{\mathbb R}^2} \Psi^2 \dx\dz
- \gamma^3 \int_{{\mathbb R}^2} \frac{k_1^2}{k_1^2+\gamma^2 k_2^2}|\hat{\psi}_z|^2\dk,
\end{eqnarray*}
and
$$\int_{{\mathbb R}^2}
\frac{k_1^2}{|k|^2} (|k|\coth|k|-1-{\textstyle\frac{1}{3}}|k|^2)|\hat{\eta}^\star|^2 \dk
\ \leq\ c\int_{{\mathbb R}^2} k_1^2 |k|^2 |\hat{\eta}^\star|^2 \dk\ =\ O(\gamma^5),$$
so that
$$
\LL_2(\eta^\star)
= \frac{\gamma}{2} \int_{{\mathbb R}^2} \Psi^2 \dx\dz
+\frac{\gamma^3}{6}\int_{{\mathbb R}^2} \Psi_x^2\dx\dz
- \frac{\gamma^3}{2} \int_{{\mathbb R}^2} \frac{k_1^2}{k_1^2+\gamma^2 k_2^2}|\hat{\psi}_z|^2\dk+O(\gamma^5).
$$

Recall that
\begin{equation}
\KK_\mathrm{nl}(\eta^\star)\ =\ O((\|\eta^\star_x\|_\infty + \|\eta^\star_z\|_\infty)^2 \|\eta^\star\|_3^2)\ =\ O(\gamma^7)
\label{Also for SSH 1}
\end{equation}
(Remark \ref{Estimates for KKnl}) and
\begin{equation}
\LL_\mathrm{nl}(\eta^\star) - \LL_3(\eta^\star)\ =\  O(\|\eta^\star\|_{1,\infty}^2 \|\eta^\star\|_3^2)\ =\ O(\gamma^5),
\label{Also for SSH 2}
\end{equation}
where
$$
\LL_3(\eta^\star) = \frac{1}{2}\int_{{\mathbb R}^2} \Big\{ (\eta^\star_x)^2\eta^\star
-\eta^\star (K^0\eta^\star)^2 -\eta^\star (L^0\eta^\star)^2\Big\}\dx\dz
$$
(Proposition \ref{Z estimate for LL} and Corollary \ref{Explicit formulae for LL3}).
We calculate
$$\int_{{\mathbb R}^2} (\eta^\star_x)^2 \eta^\star \dx\dz\ =\ O(\|\eta^\star\|_\infty\|\eta^\star_x\|_0^2)\ =\ O(\gamma^5)$$
and
\begin{eqnarray*}
\int_{{\mathbb R}^2} \eta^\star(L^0\eta^\star)^2\dx\dz
& \leq & \|\eta^\star\|_\infty \int_{{\mathbb R}^2} \frac{k_1^2k_2^2}{|k|^4} |k|^2 \coth^2|k| |\hat{\eta}^\star|^2 \dk \\
& = & \|\eta^\star\|_\infty \int_{{\mathbb R}^2} \frac{k_1^2k_2^2}{|k|^4} |\hat{\eta}^\star|^2\dk +
O(\|\eta^\star\|_\infty \|\eta^\star_z\|_0^2) \\
& = & \gamma^3 \|\eta^\star\|_\infty \int_{{\mathbb R}^2} \frac{k_1^4}{(k_1^2+\gamma^2k_2^2)^2}
|\hat{\psi}_z|^2\dk_1\dk_1 + O(\|\eta^\star\|_\infty \|\eta^\star_z\|_0^2) \\
& = & O(\gamma^5),
\end{eqnarray*}
because $|k|^2\coth^2|k|-1= O(|k|^2)$; furthermore
\begin{eqnarray*}
\lefteqn{\int_{{\mathbb R}^2}  \eta^\star(K^0\eta^\star)^2 \dx\dz} \quad\\
& = & \int_{{\mathbb R}^2} \eta^\star\left( \FF^{-1}\left[\frac{k_1^2}{|k|^2}\hat{\eta}^\star\right]\right)^{\!\!2}\dx\dz\\
& & \hspace{1cm}\mbox{}+2 \int_{{\mathbb R}^2} \eta^\star\,\FF^{-1}\left[\frac{k_1^2}{|k|^2}\eta^\star\right]\FF^{-1}\left[
\frac{k_1^2}{|k|^2}(|k|\coth|k|-1)\hat{\eta}^\star\right]\dx\dz \\
& & \hspace{1cm}\mbox{}
+ \int_{{\mathbb R}^2} \eta^\star\left( \FF^{-1}\left[\frac{k_1^2}{|k|^2}(|k|\coth|k|-1)\hat{\eta}^\star\right]
\right)^{\!\!2}\dx\dz \\
& = & \int_{{\mathbb R}^2} \eta^\star\left( \FF^{-1}\left[\frac{k_1^2}{|k|^2}\hat{\eta}^\star\right]\right)^{\!\!2}\dx\dz
+O(\|\eta^\star\|_\infty\|\eta^\star\|_0\|\eta^\star_x\|_0)+O(\|\eta^\star\|_\infty\|\eta^\star_x\|_0^2) \\
& = & \int_{{\mathbb R}^2} (\eta^\star)^3\dx\dz - 2 \int_{{\mathbb R}^2} (\eta^\star)^2 \FF^{-1}
\left[\frac{k_2^2}{|k|^2}\hat{\eta}^\star\right]\dx\dz
+ \int_{{\mathbb R}^2} \eta^\star\left(\FF^{-1}
\left[\frac{k_2^2}{|k|^2}\hat{\eta}^\star\right]\right)^{\!\!2}\dx\dz + O(\gamma^4) \\
& = & \int_{{\mathbb R}^2} (\eta^\star)^3\dx\dz
+O\left(\|\eta^\star\|_\infty \left(
 \|\eta^\star\|_0 \left\|\FF^{-1}\left[\frac{k_2^2}{|k|^2}\hat{\eta}^\star\right]\right\|_0
+\left\|\FF^{-1}\left[\frac{k_2^2}{|k|^2}\hat{\eta}^\star\right]\right\|_0^2\right)\right)+O(\gamma^4) \\
& = & \int_{{\mathbb R}^2} (\eta^\star)^3\dx\dz + O(\gamma^4)
\end{eqnarray*}
because $|k|\coth|k|-1= O(|k|)$ and
$$\int_{{\mathbb R}^2} \frac{k_2^4}{|k|^4} |\hat{\eta}^\star|^2 \dk
\ =\ \gamma^3 \int_{{\mathbb R}^2} \frac{\gamma^2 k_1^2k_2^2}{(k_1^2+\gamma^2k_2^2)^2}
|\hat{\psi}_z|^2\dk\ =\ O(\gamma^3).$$
Altogether these computations show that
\begin{equation}
\LL_\mathrm{nl}(\eta^\star) = -\frac{\gamma^3}{2}\int_{{\mathbb R}^2} \Psi^3 \dx\dz + O(\gamma^4).
\label{Also for SSH 3}
\end{equation}

Combining the above results shows that
\begin{equation}
\KK(\eta^\star) = \frac{\gamma}{2}\int_{{\mathbb R}^2} \Psi^2 \dx\dz
+ \frac{\gamma^3\beta}{2} \int_{{\mathbb R}^2} \Psi_x^2 \dx\dz + O(\gamma^5)
\label{KK(eta star)}
\end{equation}
and
\begin{eqnarray}
\lefteqn{\LL(\eta^\star) = \frac{\gamma}{2}\int_{{\mathbb R}^2} \Psi^2\dx\dz
+ \frac{\gamma^3}{6}\int_{{\mathbb R}^2} \Psi_x^2 \dx\dz} \nonumber \\
& & \hspace{1.5cm}\mbox{}
- \frac{\gamma^3}{2}\int_{{\mathbb R}^2} \frac{k_1^2}{k_1^2+\gamma^2 k_2^2} |\hat{\psi}_z|^2\dk
- \frac{\gamma^3}{2}\int_{{\mathbb R}^2} \Psi^3 \dx\dz + O(\gamma^4).
\label{LL(eta star)}
\end{eqnarray}
Let $\gamma$ be a solution of the equation $\mu=\LL(\eta^\star)$, so that
$\gamma = 2\mu/\|\Psi\|_0^2 + o(\mu)$; using equations \eqn{KK(eta star)} and
\eqn{LL(eta star)}, one finds that
\begin{eqnarray*}
\lefteqn{\JJ(\eta^\star_\mu)-2\mu} \\
& & \!\!= \KK(\eta^\star_\mu)-\LL(\eta^\star_\mu) \\
& & \!\!= \frac{\gamma^3}{2}\int_{{\mathbb R}^2} \Big((\beta-{\textstyle\frac{1}{3}})\Psi_x^2 + \Psi^3\Big)\dx\dz
+\frac{\gamma^3}{2}\int_{{\mathbb R}^2} \frac{k_1^2}{k_1^2+\gamma^2 k_2^2} |\hat{\psi}_z|^2\dk
+ O(\gamma^4) \\
& & \!\!= \frac{\gamma^3}{2}\int_{{\mathbb R}^2} \Big((\beta-{\textstyle\frac{1}{3}})\psi_{xx}^2 + \psi_x^3\Big)\dx\dz
+\frac{\gamma^3}{2}\int_{{\mathbb R}^2} \frac{k_1^2}{k_1^2+\gamma^2 k_2^2} |\hat{\psi}_z|^2\dk
+ O(\gamma^4).
\end{eqnarray*}
Finally, let us choose $\tilde{\psi} \in C_0^\infty([-{\textstyle\frac{1}{2}},{\textstyle\frac{1}{2}}]^2)$ such that
$$\int_{{\mathbb R}^2} \tilde{\psi}_x^3 \dx\dz < 0$$
and set $\psi=A\tilde{\psi}$; it follows that
\begin{eqnarray*}
\lefteqn{\JJ_\mu(\eta^\star_\mu)-2\mu} \\
& & \!\!= \frac{\gamma^3}{2}\left[
A^2\int_{{\mathbb R}^2} (\beta-{\textstyle\frac{1}{3}})\tilde{\psi}_{xx}^2 \dx\dz
+A^2\int_{{\mathbb R}^2} \frac{k_1^2}{k_1^2+\gamma^2 k_2^2} |\hat{\tilde{\psi}}_z|^2\dk
+A^3\int_{{\mathbb R}^2}\tilde{\psi}_x^3\dx\dz\right]
+ O(\gamma^4)\\
& & \!\!< 0
\end{eqnarray*}
for sufficiently large values of $A$.

\section*{Appendix C: Proof of Lemma \ref{Splitting properties 1}}
\setcounter{section}{2}
Note that the proof of this lemma given below uses results concerning
the nonlocal operator $\LL$ which are presented in the following appendix.

(i) These results are deduced from the observation that
\begin{eqnarray*}
\lefteqn{\|\eta_n\|_{H^2(M_n^{(1)} < |(x,z)| < M_n^{(2)})}^2} \qquad\\
 & = & \mathop{\int}_{B_{\!M_n^{(2)}}\!(0)} \!\!\!\!u_n(x,z)\dx\dz
- \mathop{\int}_{B_{\!M_n^{(1)}}\!(0)} \!\!\!\!u_n(x,z)\dx\dz \\
& \rightarrow & \kappa-\kappa \\
& = & 0
\end{eqnarray*}
as $n \rightarrow \infty$.

In particular, we find that
\begin{eqnarray*}
\|\eta_n^{(1)}\|_2^2 & = & \|\eta_n^{(1)}\|_{H^2(|(x,z)| < 2M_n^{(1)})}^2\\
& = & \|\eta_n^{(1)}\|_{H^2(|(x,z)| < M_n^{(1)})}^2 + \|\eta_n^{(1)}\|_{H^2(M_n^{(1)} < |(x,z)| < 2M_n^{(1)})}^2 \\
& \rightarrow & \kappa
\end{eqnarray*}
as $n \rightarrow \infty$ since $\|\eta_n^{(1)}\|_{H^2(|(x,z)| < M_n^{(1)})}^2 \rightarrow \kappa$,
$$\|\eta_n^{(1)}\|_{H^2(M_n^{(1)} < |(x,z)| < 2M_n^{(1)})}\ \leq\ \|\eta_n^{(1)}\|_{H^2(M_n^{(1)} < |(x,z)| < M_n^{(2)})}\ \leq
\ c\|\eta_n\|_{H^2(M_n^{(1)} < |(x,z)| < M_n^{(2)})}\ \rightarrow\ 0$$
as $n \rightarrow \infty$, and a similar argument yields the second limit. Finally, observe that
\begin{eqnarray*}
\lefteqn{\|\eta_n-\eta_n^{(1)}-\eta_n^{(2)}\|_2^2} \qquad\\
& = & \|\eta_n-\eta_n^{(1)}-\eta_n^{(2)}\|_{H^2(M_n^{(1)} < |(x,z)| < M_n^{(2)})}^2 \\
& \leq & \|\eta_n\|_{H^2(M_n^{(1)} < |(x,z)| < M_n^{(2)})}^2
+\|\eta_n^{(1)}\|_{H^2(M_n^{(1)} < |(x,z)| < M_n^{(2)})}^2
+\|\eta_n^{(2)}\|_{H^2(M_n^{(1)} < |(x,z)| < M_n^{(2)})}^2 \\
& \leq & c\|\eta_n\|_{H^2(M_n^{(1)} < |(x,z)| < M_n^{(2)})}^2 \\
& \rightarrow & 0
\end{eqnarray*}
as $n \rightarrow \infty$.

(ii) The equation
$$\eta_n^{(1)}+\eta_n^{(2)} = \eta_n \tilde{\chi},$$
where
$$\tilde{\chi}(x,z) := \chi\left(\frac{|(x,z)|}{M_n^{(1)}}\right)\!\!+\!1\!-\chi\left(\frac{|(x,z)|}{M_n^{(2)}}\right)$$
satisfies
$$|\partial_x^i\partial_z^j \tilde{\chi}(x,z)| \leq \frac{1}{(M_n^{(1)})^{i+j}}, \qquad
(x,z) \in {\mathbb R}^2,$$
and the estimate
$$\|\eta_n \tilde{\chi}\|_3 \leq \|\eta_n\|_3\|\tilde{\chi}\|_\infty + c \|\eta_n\|_3\max_{i+j=1,2,3} \|\partial_x^i\partial_z^j \chi\|_\infty$$
imply that
\begin{equation}
\|\eta_n^{(1)}+\eta_n^{(2)}\|_3^2\ \leq\ \|\eta_n\|_3^2 + o(1); \label{Sum in Y}
\end{equation}
replacing $\{\eta_n\}$ by a subsequence if necessary,
one finds from the above result and the bound $\sup \|\eta_n\|_3 < M$ that
$\sup \|\eta_n^{(1)}+\eta_n^{(2)}\|_3 < M$.
The estimates on the suprema of $\|\eta_n^{(1)}\|_3$ and $\|\eta_n^{(2)}\|_3$ follow from the inequalities
$$\|\eta_n^{(1)}\|_3 \leq \|\eta_n^{(1)}+\eta_n^{(2)}\|_3, \qquad \|\eta_n^{(2)}\|_3 \leq \|\eta_n^{(1)}+\eta_n^{(2)}\|_3.$$

(iii) Clearly
$$\KK(\eta_n^{(1)}+\eta_n^{(2)}) - \KK(\eta_n^{(1)}) - \KK(\eta_n^{(2)}) \rightarrow 0$$
as $n \rightarrow \infty$; in fact the expression appearing in this limit vanishes identically since
$\{\eta_n^{(1)}\}$ and $\{\eta_n^{(2)}\}$ have disjoint supports and $\KK$ is a local operator.
Because the derivative of $\KK$ is bounded on $U$, we find that
$$|\KK(\eta_n)-\KK(\eta_n^{(1)}+\eta_n^{(2)})|\ \leq\ c\|\eta_n-\eta_n^{(1)}-\eta_n^{(2)}\|_3\ \rightarrow\ 0$$
and therefore that
$$
\KK(\eta_n)-\KK(\eta_n^{(1)})-\KK(\eta_n^{(2)})\ =\ 
\underbrace{\KK(\eta_n) - \KK(\eta_n^{(1)}+\eta_n^{(2)})}_{\displaystyle =o(1)}
+ \underbrace{\KK(\eta_n^{(1)}+\eta_n^{(2)}) - \KK(\eta_n^{(1)}) - \KK(\eta_n^{(2)})}_{\displaystyle =o(1)}$$
as $n \rightarrow \infty$.

Turning to the result for $\LL$, note that according to the above argument
it suffices to establish that
\begin{equation}
\lim_{n \rightarrow \infty}\Big(\LL(\eta_n^{(1)}+\eta_n^{(2)}) - \LL(\eta_n^{(1)}) - \LL(\eta_n^{(2)})\Big)=0.
\label{Intermediate L estimate}
\end{equation}
This limit is in turn obtained by approximating $\{\eta_n^{(1)}\}$ by a sequence $\{\eta_n^{(3)}\}$
of functions in $U$ with uniform compact support (and for which $\sup \|\eta_n^{(3)}+\eta_n^{(2)}\|_3 < M$);
in Appendix D (Theorem \ref{Key splitting theorem}) it is shown that
$$\lim_{n \rightarrow \infty}\Big(\LL(\eta_n^{(3)}+\eta_n^{(2)}) - \LL(\eta_n^{(3)}) - \LL(\eta_n^{(2)})\Big)=0.$$

Choose $\varepsilon>0$. The `dichotomy' property of the sequence $\{u_n\}$
asserts the existence of a positive real number $R$ such that
$$\|\eta_n\|_{H^2(|(x,z)|<R)}^2=\mathop{\int}_{B_R(0)} u_n(x,z) \dx\dz
\ \geq\ \kappa - {\textstyle\frac{1}{2}}\varepsilon^2.$$
Taking $n$ large enough so that $M_n^{(1)}>2R$, we find that
$$\|\eta_n^{(1)}\|_2^2 - \|\eta_n^{(1)}\|_{H^2(|(x,z)|>R)}^2\ =\ \|\eta_n^{(1)}\|_{H^2(|(x,z)|<R)}^2
=\|\eta_n\|_{H^2(|(x,z)|<R)}^2\ \geq\ \kappa-{\textstyle\frac{1}{2}}\varepsilon^2,$$
whereby
$$\|\eta_n^{(1)}\|_{H^2(|(x,z)|>R)}^2\ \leq\ \|\eta_n^{(1)}\|_2^2
-(\kappa-{\textstyle\frac{1}{2}}\varepsilon^2)\ <\ \varepsilon^2$$
for sufficiently large $n$, since $\|\eta_n^{(1)}\|_2^2 \rightarrow \kappa$ as $m \rightarrow \infty$.
Define $\eta_n^{(3)}=\eta_n^{(1)}\chi_R$ ($=\eta_n\chi_R$), where
$$\chi_R(x,z) = \chi\left(\frac{|(x,z)|}{R}\right),$$
so that $\supp \eta_n^{(3)} \subset B_{2R}(0,0)$
and
$$
\|\eta_n^{(1)}-\eta_n^{(3)}\|_2\ =\ \|\eta_n^{(1)}-\eta_n^{(3)}\|_{H^2(|(x,z)|>R)}\ \leq\ c\left\|1-\chi_R\right\|_{2,\infty}
\|\eta_n^{(1)}\|_{H^2(|(x,z)|>R)}\ =\ O(\varepsilon),
$$
$$
\|\eta_n^{(1)}-\eta_n^{(3)}\|_3\ \leq\ c\left\|1-\chi_R\right\|_{3,\infty}
\|\eta_n^{(1)}\|_3\ =\ O(1);$$
it follows by interpolation that
$$\|\eta_n^{(1)}-\eta_n^{(3)}\|_{2+t}\ \leq\ \|\eta_n^{(1)}-\eta_n^{(3)}\|_2^{1-t} \|\eta_n^{(1)}-\eta_n^{(3)}\|_3^t
\ \leq\ c\varepsilon^{1-t}.$$

Proceeding as in part (ii), we find that
$$\|\partial_x^i\partial_z^j(\eta_n^{(3)}+\eta_n^{(2)})\|_0^2\ \leq\ \|\partial_x^i\partial_z^j\eta_n\|_0^2 + O(R^{-2}) \|\eta_n\|_2^2\ \leq\ \|\partial_x^i\partial_z^j\eta_n\|_0^2 + O(R^{-2})$$
and furthermore
\begin{eqnarray*}
\|\eta_n^{(3)} + \eta_n^{(2)}\|_2^2
& = & \|\eta_n^{(1)} + \eta_n^{(2)}\|_2^2 + O(\varepsilon^{-2}) \\
& \leq & \|\eta_n\|_2^2 + o(1) + O(\varepsilon^{-2}),
\end{eqnarray*}
so that $\sup\|\eta_n^{(3)} + \eta_n^{(2)}\|_3$ and hence $\sup\|\eta_n^{(3)}\|_3$ are strictly less than $M$
for sufficiently small values of $\varepsilon$ (where $R$ is replaced by a larger number if necessary).
The estimate \eqn{Intermediate L estimate} is now obtained by choosing $t \in (\frac{1}{2},1)$ and noting that
\begin{eqnarray*}
\lefteqn{\LL(\eta_n^{(1)}+\eta_n^{(2)}) - \LL(\eta_n^{(1)}) - \LL(\eta_n^{(2)})} \qquad\\
& & = \LL(\eta_n^{(3)}+\eta_n^{(2)}) - \LL(\eta_n^{(3)}) - \LL(\eta_n^{(2)}) + O(\|\eta_n^{(1)}-\eta_n^{(3)}\|_{2+t}) \\
& & = o(1) + O(\varepsilon^{1-t})
\end{eqnarray*}
as $m \rightarrow \infty$,
in which Theorem \ref{Key splitting theorem} has been used.

The estimates for $\KK^\prime$ and $\LL^\prime$ are obtained in a similar fashion.

(iv) Suppose that $\lim_{n \rightarrow \infty} \LL(\eta_n^{(2)})=0$.
The calculation
\begin{eqnarray*}
\lim_{n \rightarrow \infty} \JJ_{\rho,\mu}(\eta_n^{(1)}) & = & 
\lim_{n \rightarrow \infty} \left\{\KK(\eta_n^{(1)})+\frac{\mu^2}{\LL(\eta_n^{(1)})}+\rho(\|\eta_n^{(1)}\|_3^2)\right\} \\
& = & \lim_{n \rightarrow \infty}
\left\{\KK(\eta_n^{(1)}) + \frac{\mu^2}{\LL(\eta_n)}+\rho(\|\eta_n^{(1)}\|_3^2)\right\} \\
& \leq & \lim_{n \rightarrow \infty}
\left\{\KK(\eta_n^{(1)})+\KK(\eta_n^{(2)})+ \frac{\mu^2}{\LL(\eta_n)}+\rho(\|\eta_n^{(1)}+\eta_n^{(2)}\|_3^2)\right\} \\
& \leq & \lim_{n \rightarrow \infty} \left\{\KK(\eta_n)+\frac{\mu^2}{\LL(\eta_n)} + \rho(\|\eta_n\|_3^2)\right\} \\
& = & c_{\rho,\mu},
\end{eqnarray*}
in which part (iii) and the facts that $\KK(\eta_n^{(2)})>0$ and
$$\|\eta_n^{(1)}\|_3 \leq \|\eta_n^{(1)}+\eta_n^{(2)}\|_3, \qquad
\lim_{n \rightarrow \infty} \|\eta_n^{(1)}+\eta_n^{(2)}\|_3
\leq \lim_{n \rightarrow \infty} \|\eta_n\|_3$$
have been used, shows
that $\JJ_{\rho,\mu}(\eta_n^{(1)}) \rightarrow c_{\rho,\mu}$ and $\KK(\eta_n^{(2)}) \rightarrow 0$ as $n
\rightarrow \infty$. It follows that
$$\|\eta_n^{(2)}\|_1^2\ \leq\ c\KK(\eta_n^{(2)})\ \rightarrow
\ 0$$
(see Proposition \ref{Quadratic estimates}) and hence $\|\eta_n^{(2)}\|_2 \rightarrow 0$ as $n \rightarrow \infty$,
which contradicts part (i). One similarly finds that
the assumption $\lim_{n \rightarrow \infty} \LL(\eta_n^{(1)})=0$ leads to the contradiction
$\|\eta_n^{(1)}\|_2 \rightarrow 0$ as $n \rightarrow \infty$.

\section*{Appendix D: Pseudo-local properties of the operator ${\mathcal L}$}
\setcounter{section}{4}
Consider two sequences $\{v_n^{(1)}\}$, $\{v_n^{(2)}\}$ with $\sup \|v_n^{(1)}+v_n^{(1)}\|_3 < M$ and 
$\supp v_n^{(1)} \subset B_{2R}(0)$, $\supp v_n^{(2)} \subset {\mathbb R}^2 \sm B_{N_n}(0)$,
where $R>0$ and $\{N_n\}$ is an increasing, unbounded sequence
of positive real numbers. Clearly
\begin{eqnarray*}
& & \KK(v_n^{(1)}+v_n^{(2)}) - \KK(v_n^{(1)})-\KK(v_n^{(2)}) \rightarrow0, \\
& & \KK^\prime(v_n^{(1)}+v_n^{(2)}) - \KK^\prime(v_n^{(1)})-\KK^\prime(v_n^{(2)})\rightarrow0, \\
& & \langle \KK^\prime(v_n^{(2)}),v_n^{(1)}\rangle_0 \rightarrow 0
\end{eqnarray*}
as $n \rightarrow \infty$; indeed the expressions appearing in these limits are identically zero for
for sufficiently large values of $m$ because the supports of $\{v_n^{(1)}\}$ and $\{v_n^{(2)}\}$
are disjoint and $\KK$, $\KK^\prime$ are local operators. In this appendix we show that
the result is also valid for the nonlocal operators $\LL$ and $\LL^\prime$.
Our strategy is to approximate the expressions appearing in the limits by a finite number of
terms involving integral operators which are estimated using the following result.

\begin{proposition} \label{Basic integral estimate}
The integral
$$\int_{N_1^{x,z}}\int_{N_2^{\tilde{x},\tilde{z}}} \left(\frac{1}{|x-\tilde{x}|^2 + |z-\tilde{z}|^2}\right)^{\!\!\tilde{m}}
\dtildex\dtildez\dx\dz,$$
in which
$$N_1^{x,z} = \{(x,z): |(x,z)| < 2R\}, \qquad N_2^{\tilde{x},\tilde{z}} = \{(\tilde{x},\tilde{z}): |(\tilde{x},\tilde{z})| > N_n\}$$
and  $\tilde{m} > 1$, converges to zero as $n \rightarrow \infty$.
\end{proposition}
{\bf Proof.} Define
$$N_3^{x^\prime,z^\prime} = \{(x^\prime,z^\prime): |(x^\prime,z^\prime)| > N_n-2R\}$$
and observe that
\begin{eqnarray*}
\lefteqn{\int_{N_1^{x,z}}\int_{N_2^{\tilde{x},\tilde{z}}} \left(\frac{1}{|x-\tilde{x}|^2 + |z-\tilde{z}|^2}\right)^{\!\!\tilde{m}}
\dtildex\dtildez\dx\dz} \qquad\qquad\\
& \leq & \int_{N_1^{x,z}}\int_{N_3^{x^\prime,z^\prime}} \left(\frac{1}{|x^\prime|^2 + |z^\prime|^2}\right)^{\!\!\tilde{m}}
\dx^\prime\dz^\prime \dx\dz \\
& = & \frac{\pi^2}{\tilde{m}-1}\frac{4R^2}{(N_n-2R)^{2\tilde{m}-2}} \\
& \rightarrow & 0
\end{eqnarray*}
as $n \rightarrow \infty$.\qed

In Theorems \ref{Integral operators 1}, \ref{Integral operators 2} below we list the integral
operators appearing in the definitions of $\LL(v_n^{(1)}+v_n^{(2)}) - \LL(v_n^{(1)})-\LL(v_n^{(2)})$,
$\LL^\prime(v_n^{(1)}+v_n^{(2)}) - \LL^\prime(v_n^{(1)})-\LL^\prime(v_n^{(2)})$ and
$\langle \LL^\prime(v_n^{(2)}),v_n^{(1)}\rangle_0$ and
obtain the relevant estimates using Proposition \ref{Basic integral estimate}.

\begin{theorem} \label{Integral operators 1}
Define integral operators
\begin{eqnarray*}
\GG_1(\cdot) & = & \FF^{-1}\left[ -\int_0^1 k_1^2 (\i k_1)^{m_1} (\i k_2)^{m_2} |k|^{2m_3}
G \FF[\cdot]\dtildey \right], \\
\GG_2(\cdot) & = & \FF^{-1}\left[ -\int_0^1 k_1k_2 (\i k_1)^{m_1} (\i k_2)^{m_2} |k|^{2m_3}
G \FF[\cdot]\dtildey \right], \\
\GG_3(\cdot) & = & \FF^{-1}\left[ -\int_0^1 k_2^2 (\i k_1)^{m_1} (\i k_2)^{m_2} |k|^{2m_3}
G \FF[\cdot]\dtildey \right], \\
\GG_4(\cdot) & = & \FF^{-1}\left[ \int_0^1 (\i k_1)^{m_1+1} (\i k_2)^{m_2} |k|^{2m_3}
G_y \FF[\cdot]\dtildey \right], \\
\GG_5(\cdot) & = & \FF^{-1}\left[ \int_0^1 (\i k_1)^{m_1+1} (\i k_2)^{m_2} |k|^{2m_3}
H_y \FF[\cdot]\dtildey \right], \\
\GG_6(\cdot) & = & \FF^{-1}\left[ \int_0^1 (\i k_1)^{m_1} (\i k_2)^{m_2+1} |k|^{2m_3}
G_y \FF[\cdot]\dtildey \right], \\
\GG_7(\cdot) & = & \FF^{-1}\left[ \int_0^1 (\i k_1)^{m_1} (\i k_2)^{m_2+1} |k|^{2m_3}
H_y \FF[\cdot]\dtildey \right], \\
\GG_8(\cdot) & = & \FF^{-1}\left[ \int_0^1 (\i k_1)^{m_1} (\i k_2)^{m_2} |k|^{2m_3+2}
H \FF[\cdot]\dtildey \right],
\end{eqnarray*}
where $m_1$, $m_2$ and $m_3$ are non-negative integers.
\begin{list}{(\roman{count})}{\usecounter{count}}
\item
The estimates
$$\|\GG_j(P_n)\|_{L^2(|(x,z)|<2R)} \rightarrow 0, \qquad j=1,\ldots,8$$
hold for every sequence $\{P_n\} \subset L^2(\Sigma)$ of functions with the properties that
$$\supp P_n \subset \{(x,y,z) \in \Sigma: (x,z) \in {\mathbb R}^2\sm  B_{N_n}(0,0)\},
\qquad \|P_n\|_0 \leq c.$$
\item
The estimates
$$\|\GG_j(Q_n)\|_{L^2(|(x,z)|>N_n)} \rightarrow 0, \qquad j=1,\ldots,8$$
hold for every sequence $\{Q_n\} \subset L^2(\Sigma)$ of functions with the properties that
$$\supp Q_n \subset \{(x,y,z) \in \Sigma: (x,z) \in B_{2R}(0,0)\},
\qquad \|Q_n\|_0 \leq c.$$
\end{list}
\end{theorem}
{\bf Proof.} Observe that
\begin{eqnarray*}
G(y,\tilde{y}) & = & -\frac{\e^{-|k||\tilde{y}-y|}}{2|k|(1-\e^{-2|k|})} - \frac{\e^{-|k|(\tilde{y}+y)}}{2|k|(1-\e^{-2|k|})}
-\frac{\e^{-|k|(2-\tilde{y}-y)}}{2|k|(1-\e^{-2|k|})} - \frac{\e^{-|k|(1-|\tilde{y}-y|)}}{2|k|(\e^{|k|}-\e^{-|k|})} ,\\
G_y(y,\tilde{y}) & = & -s\frac{\e^{-|k||\tilde{y}-y|}}{2(1-\e^{-2|k|})} + \frac{\e^{-|k|(\tilde{y}+y)}}{2(1-\e^{-2|k|})}
-\frac{\e^{-|k|(2-\tilde{y}-y)}}{2(1-\e^{-2|k|})}+s\frac{\e^{-|k|(1-|\tilde{y}-y|)}}{2(\e^{|k|}-\e^{-|k|})},
\end{eqnarray*}
where
$$s=\left\{\begin{array}{cl} 1, & \quad y < \tilde{y}, \\[2mm] -1, & \quad \tilde{y} <y, \end{array}\right.$$
and there are of course analogous formulae for $H$ and $H_y$. Write
$$\GG_1 = \GG_{1,\mathrm{a}} + \GG_{1,\mathrm{b}} + \GG_{1,\mathrm{c}} + \GG_{1,\mathrm{d}},$$
where
$$\GG_{1,\mathrm{a}}(\cdot) = \FF^{-1}\left[ -\int_0^1 k_1^2 (\i k_1)^{m_1} (\i k_2)^{m_2} |k|^{2m_3}
\frac{\e^{-|k||\tilde{y}-y|}}{2|k|(1-\e^{-2|k|})} \FF[\cdot]\dtildey \right]$$
and $\GG_{1,\mathrm{b}}$, $\GG_{1,\mathrm{c}}$ and $\GG_{1,\mathrm{d}}$ are defined by
replacing the final multiplier in this formula by respectively the second, third and fourth summand
in the formula for $G$; the other integral operators $\GG_j$, $j=2,\ldots,8$ admit similar decompositions.
We proceed by estimating $\| G_{1,\mathrm{a}}(P_n)\|_{L^2(|(x,z)|<2R)}$;
the same technique yields the corresponding
estimates for the other operators. Our strategy is to show that the quantity
$$I=\FF^{-1}\left[ -k_1^2 (\i k_1)^{m_1} (\i k_2)^{m_2} |k|^{2m_3}
\frac{\e^{-|k||\tilde{y}-y|}}{2|k|(1-\e^{-2|k|})}\right]$$
can be written as a finite sum $\sum_i K_i(x,z;y,\tilde{y}),$ where $K_i$ is $O(|(x,z)|^{-m_i})$ for
some $m_i > 1$, uniformly in $y \neq \tilde{y}$, since
\begin{eqnarray*}
\lefteqn{\left\|\FF^{-1}\left[ \int_0^1 \FF[K_i(x,z;y,\tilde{y})]\FF[P_n]\dtildey\right]\right\|_{L^2(|(x,z)|<2R)}} \quad\\
& & \!\!\!\!= \Bigg(\int_{N_1^{x,z}}\int_0^1\Bigg|\int_{N_2^{\tilde{x},\tilde{z}}} \int_0^1 K_i(x-\tilde{x},z-\tilde{z};y,\tilde{y}) P_n(\tilde{x},\tilde{y},\tilde{z})\dtildey\dtildex\dtildez \Bigg|^2 \dy\dx\dz\Bigg)^{\!\!\frac{1}{2}} \\
& & \!\!\!\!\leq \Bigg( \int_{N_1^{x,z}}\int_0^1\int_{N_2^{\tilde{x},\tilde{z}}} \int_0^1 |K_i(x-\tilde{x},z-\tilde{z};y,\tilde{y})|^2 \dtildey\dtildex\dtildez \dy\dx\dz\Bigg)^{\!\!\frac{1}{2}} \|P_n\|_0 \\
& & \!\!\!\! \leq c\Bigg(\int_{N_1^{x,z}}\int_{N_2^{\tilde{x},\tilde{z}}} \left(\frac{1}{|x-\tilde{x}|^2 + |z-\tilde{z}|^2}\right)^{\!\!m_i}
\dtildex\dtildez\dx\dz\Bigg)^{\!\!\frac{1}{2}}  \|P_n\|_0 \\
& & \!\!\!\! \rightarrow 0
\end{eqnarray*}
as $n \rightarrow \infty$, where we have used Proposition \ref{Basic integral estimate} and the fact that
$\|P_n\|_0 \leq c$.

Choose $\delta>0$, define
\begin{eqnarray*}
I_1 & = & \FF^{-1}\left[ -k_1^2 (\i k_1)^{m_1} (\i k_2)^{m_2} |k|^{2m_3}\chi_\delta(|k|)
\frac{\e^{-|k||\tilde{y}-y|}}{2|k|(1-\e^{-2|k|})}\right], \\
I_2 & = & \FF^{-1}\left[ -k_1^2 (\i k_1)^{m_1} (\i k_2)^{m_2} |k|^{2m_3}(1-\chi_\delta(|k|))
\frac{\e^{-|k||\tilde{y}-y|}}{2|k|(1-\e^{-2|k|})}\right],
\end{eqnarray*}
where
$$\chi_\delta(r) = \chi\left(\frac{r}{\delta}\right), \qquad r \in [0,\infty),$$
and observe that
$$I_1=
-\int_{{\mathbb R}^2} \left\{ \begin{array}{c} -\i\sin (k_1 x) \\ \cos(k_1 x) \end{array}\right\}
\e^{\i k_2 z}\left[
k_1^2 (\i k_1)^{m_1} (\i k_2)^{m_2} |k|^{2m_3-2}
\frac{|k|\chi_\delta(|k|)}{2(1-\e^{-2|k|})}\e^{-|k||\tilde{y}-y|}\right]\dk,$$
where we write $-\i\sin(k_1 x)$ for odd values of $m_1$ and $\cos(k_1 x)$ for even values of $m_1$. Taking the limit $\varepsilon \downarrow 0$ in the equation
\begin{eqnarray*}
\lefteqn{\int\limits_{{\mathbb R}\setminus [-\varepsilon,\varepsilon]}\!\!\int_{-\infty}^\infty
\left\{ \begin{array}{c} \i\sin (k_1 x) \\ -\cos(k_1 x) \end{array}\right\}
\e^{\i k_2 z}\left[
k_1^2 (\i k_1)^{m_1} (\i k_2)^{m_2} |k|^{2m_3-2}
\frac{|k|\chi_\delta(|k|)}{2(1-\e^{-2|k|})}\e^{-|k||\tilde{y}-y|}\right]\dk} \\
& = & -\frac{1}{x^2}\int\limits_{{\mathbb R}\setminus [-\varepsilon,\varepsilon]}\int_{-\infty}^\infty
\left\{ \begin{array}{c} \i\partial_{k_1}^2(\sin (k_1 x)) \\ \partial_{k_1}^2
(1-\cos(k_1 x)) \end{array}\right\}
\e^{\i k_2 z} \\
& & \hspace{1.5in} \times \left[
k_1^2 (\i k_1)^{m_1} (\i k_2)^{m_2} |k|^{2m_3-2}
\frac{|k|\chi_\delta(|k|)}{2(1-\e^{-2|k|})}\e^{-|k||\tilde{y}-y|}\right]\dk \\
& = & -\frac{1}{x^2}\int\limits_{{\mathbb R}\setminus [-\varepsilon,\varepsilon]}\int_{-\infty}^\infty
\left\{ \begin{array}{c} \i \sin (k_1 x) \\ 1-\cos(k_1 x) \end{array}\right\} \e^{\i k_2 z} \\
& & \hspace{1.5in} \times 
\partial_{k_1}^2\left[
k_1^2 (\i k_1)^{m_1} (\i k_2)^{m_2} |k|^{2m_3-2}
\frac{|k|\chi_\delta(|k|)}{2(1-\e^{-2|k|})}\e^{-|k||\tilde{y}-y|}\right]\dk,
\end{eqnarray*}
one finds that
$$I _1= -\frac{1}{x^2}\int_{{\mathbb R}^2}
\left\{ \begin{array}{c} \i \sin (k_1 x) \\ 1-\cos(k_1 x) \end{array}\right\}
\e^{\i k_2 z}\partial_{k_1}^2\left[
k_1^2 (\i k_1)^{m_1} (\i k_2)^{m_2} |k|^{2m_3-2}
\frac{|k|\chi_\delta(|k|)}{2(1-\e^{-2|k|})}\e^{-|k||\tilde{y}-y|}\right]\dk.$$
(It is necessary to exclude the set $k_2 \in [-\varepsilon,\varepsilon]$ to guarantee continuity
of all integrands and hence allow integration by parts in the above calculation.) It follows that
\begin{eqnarray*}
|I_1| & \leq & \frac{1}{|x|^{3/2}} \int_{{\mathbb R}^2}
\underbrace{\frac{1}{|k_1 x|^{1/2}}
\left\{ \begin{array}{c} |\sin (k_1 x)|^{1/2} \\ |1-\cos(k_1 x)|^{1/2} \end{array}\right\}}_{\displaystyle = O(1)}
\,\underbrace{\left\{ \begin{array}{c} |\sin (k_1 x)|^{1/2} \\ |1-\cos(k_1 x)|^{1/2} \end{array}\right\}}_{\displaystyle = O(1)} \\
& & \hspace{1in} \times
|k_1|^{1/2}\underbrace{\left|\partial_{k_1}^2\left[
k_1^2 (\i k_1)^{m_1} (\i k_2)^{m_2} |k|^{2m_3-2}
\frac{|k|\chi_\delta(|k|)}{2(1-\e^{-2|k|})}\e^{-|k||\tilde{y}-y|}\right]\right|}_{\displaystyle = O(|k|^{m_1+m_2+2m_3-2})}
\dk \\
& = & O(|x|^{-3/2}).
\end{eqnarray*}
A similar calculation yields the complementary estimate $|I_1| = O(|z|^{-3/2})$, and combining these results,
we conclude that
$$|I_1| = O(|(x,z)|^{-3/2}).$$

Turning to $I_2$, note that
\begin{eqnarray*}
I_2 & = & -\int_{{\mathbb R}^2}
\e^{\i(k_1 x + k_2 z)}\left[
k_1^2 (\i k_1)^{m_1} (\i k_2)^{m_2} |k|^{2m_3-1}\e^{-|k||\tilde{y}-y|}
\frac{1-\chi_\delta(|k|)}{2(1-\e^{-2|k|})}\right]\dk \\
& = & - \frac{1}{(\i x)^n} \int_{{\mathbb R}^2}
\partial_{k_1}^n(\e^{\i(k_1 x + k_2 z)})\left[
k_1^2 (\i k_1)^{m_1} (\i k_2)^{m_2} |k|^{2m_3-1}\e^{-|k||\tilde{y}-y|}
\frac{1-\chi_\delta(|k|)}{2(1-\e^{-2|k|})}\right]\dk \\
& = & - \left(\frac{\i}{x}\right)^n \int_{{\mathbb R}^2}
\e^{\i(k_1 x + k_2 z)}\partial_{k_1}^n\left[
k_1^2 (\i k_1)^{m_1} (\i k_2)^{m_2} |k|^{2m_3-1}\e^{-|k||\tilde{y}-y|}
\frac{1-\chi_\delta(|k|)}{2(1-\e^{-2|k|})}\right]\dk \\
& = & - \left(\frac{\i}{x}\right)^n \int_{{\mathbb R}^2}
\e^{\i(k_1 x + k_2 z)}\partial_{k_1}^n\left[
k_1^2 (\i k_1)^{m_1} (\i k_2)^{m_2} |k|^{2m_3-1}\e^{-|k||\tilde{y}-y|}\right]
\frac{1-\chi_\delta(|k|)}{2(1-\e^{-2|k|})}\dk \\[2mm]
& & \qquad\mbox{}+O(|x|^{-n}),
\end{eqnarray*}
where $n=m_1+m_2+2m_3+5$; here we use the fact that
\begin{eqnarray*}
\lefteqn{\left|\e^{\i(k_1 x + k_2 z)}\partial_{k_1}^i[
k_1^2 (\i k_1)^{m_1} (\i k_2)^{m_2} |k|^{2m_3-1}\e^{-|k||\tilde{y}-y|}]
\partial_{k_1}^j\left[\frac{1-\chi_\delta(|k|)}{2(1-\e^{-2|k|})}\right]\right|} \\
& & \leq \underbrace{c|k|^{m_1+m_2+2m_3+1-i} \left|\partial_{k_1}^j\left[\frac{1-\chi_\delta(|k|)}{2(1-\e^{-2|k|})}\right]\right|}_{\displaystyle \in L^1({\mathbb R}^2)}\hspace{1.5in}
\end{eqnarray*}
because every derivative of $(1-\chi_\delta(|k|))(1-\mathrm{e}^{-2|k|})^{-1}$ either has compact support or
decays exponentially quickly as $|k| \rightarrow \infty$.

The next step is to use the computations
\begin{eqnarray*}
\lefteqn{\partial_{k_1}^n \left[
k_1^2 (\i k_1)^{m_1} (\i k_2)^{m_2} |k|^{2m_3-1}\e^{-|k||\tilde{y}-y|}\right] } \qquad \quad\\
& & = \sum_{p=0}^n \left(\begin{array}{c} n \\ p \end{array}\right) \partial_{k_1}^p
[\e^{-|k||\tilde{y}-y|}]\partial_{k_1}^{n-p} \left[
k_1^2 (\i k_1)^{m_1} (\i k_2)^{m_2} |k|^{2m_3-1}\right],
\end{eqnarray*}
$$\partial_{k_1}^p[\e^{-|k||\tilde{y}-y|}] = \sum_{q=0}^p |\tilde{y}-y|^q \e^{-|k||\tilde{y}-y|} O(|k|^{q-p}),$$
$$\partial_{k_1}^{n-p} \left[ k_1^2 (\i k_1)^{m_1} (\i k_2)^{m_2} |k|^{2m_3-1}\right] = O(|k|^{p-4}),$$
to find that
\begin{eqnarray*}
\lefteqn{\left| \partial_{k_1}^n \left[k_1^2 (\i k_1)^{m_1} (\i k_2)^{m_2} |k|^{2m_3-1}\e^{-|k||\tilde{y}-y|}\right] \right|} \qquad \qquad\\
& = & \sum_{p=0}^n \sum_{q=0}^p |\tilde{y}-y|^q \e^{-|k||\tilde{y}-y|}O(|k|^{q-p})O(|k|^{p-4}) \\
& = & \sum_{q=0}^n |\tilde{y}-y|^q \e^{-|k||\tilde{y}-y|}O(|k|^{q-4})\\
& = & O(|k|^{-3}) + O(|\tilde{y}-y|^2)\e^{-|k||\tilde{y}-y|} + \sum_{q=5}^n |\tilde{y}-y|^{q}O(|k|^{q-4})\e^{-|k||\tilde{y}-y|}
\end{eqnarray*}
Altogether these calculations show that
\begin{eqnarray*}
I_2 & = & O(|x|^{-n})\left[1+\int_{{\mathbb R}^2} |k|^{-3}\dk
+ \sum_{q=4}^n \int_{{\mathbb R}^2} |\tilde{y}-y|^{q-2} |k|^{q-4}\e^{-|k||\tilde{y}-y|} \dk\right] \\
& = & O(|x|^{-n}) \left[1+\int_{{\mathbb R}^2} |k|^{-3}\dk
+ \sum_{q=4}^n \int_{{\mathbb R}^2} |s|^{q-4} \e^{-|s|} \ds_1\ds_2\right] \\
& = & O(|x|^{-n})
\end{eqnarray*}

A similar calculation yields the complementary estimate $|I_2| = O(|z|^{-n})$,
from which it follows that
$$|I_2| = O(|(x,z)|^{-n}).\eqno{\Box}$$

The following result is proved in the same way as Theorem \ref{Integral operators 1}.

\begin{theorem} \label{Integral operators 2}
Define integral operators
\begin{eqnarray*}
\GG_9(\cdot) & = & \FF^{-1}\left[ -k_1^2 \frac{\cosh |k|y}{|k|\sinh|k|} \FF[\cdot]\right], \\
\GG_{10}(\cdot) & = & \FF^{-1}\left[ -k_1k_2 \frac{\cosh |k|y}{|k|\sinh|k|} \FF[\cdot]\right], \\
\GG_{11}(\cdot) & = & \FF^{-1}\left[ \i k_1 \frac{\sinh |k|y}{\sinh|k|} \FF[\cdot]\right],
\end{eqnarray*}
\begin{list}{(\roman{count})}{\usecounter{count}}
\item
The estimates
$$\|\GG_j(p_n)\|_{H^2(|(x,z)|<2R)} \rightarrow 0, \qquad j=9,10,11$$
hold for every sequence $\{p_n\} \subset L^2({\mathbb R}^2)$ of functions with the properties that
$$\supp p_n \subset {\mathbb R}^2\sm  B_{N_n}(0,0),
\qquad \|p_n\|_0 \leq c.$$
\item
The estimates
$$\|\GG_j(q_n)\|_{H^2(|(x,z)|>N_n)} \rightarrow 0, \qquad j=9,10,11$$
hold for every sequence $\{q_n\} \subset L^2({\mathbb R}^2)$ of functions with the properties that
$$\supp q_n \subset B_{2R}(0,0),
\qquad \|q_n\|_0 \leq c.$$
\end{list}
\end{theorem}

The following lemma is the key step in the proof that 
$\LL(v_n^{(1)}+v_n^{(1)}) - \LL(v_n^{(2)})-\LL(v_n^{(2)}) \rightarrow 0$
and $\LL^\prime(v_n^{(1)}+v_n^{(1)}) - \LL^\prime(v_n^{(2)})-\LL^\prime(v_n^{(2)}) \rightarrow 0$
as $n \rightarrow \infty$.

\begin{lemma} \label{Key splitting lemma}
Let $\{v_n^{(1)}\}$ and $\{v_n^{(2)}\}$ be sequences in $U$ 
which have the properties that\linebreak
$\supp v_n^{(1)} \subset B_{2R}(0)$, $\supp v_n^{(2)} \subset {\mathbb R}^2 \sm B_{N_n}(0)$
and $\sup \|v_n^{(1)}+v_n^{(2)}\|_3 < M$. The estimates
\begin{eqnarray*}
& & \lim_{n \rightarrow \infty} \|u_x^j(v_n^{(1)}+v_n^{(2)})-u_x^j(v_n^{(1)})\|_{H^1(|(x,z)|<2R)} = 0, \\
& & \lim_{n \rightarrow \infty} \|u_y^j(v_n^{(1)}+v_n^{(2)})-u_y^j(v_n^{(1)})\|_{H^1(|(x,z)|<2R)} =0, \\
& & \lim_{n \rightarrow \infty} \|u_z^j(v_n^{(1)}+v_n^{(2)})-u_z^j(v_n^{(1)})\|_{H^1(|(x,z)|<2R)} =0
\end{eqnarray*}
and
\begin{eqnarray*}
& & \lim_{n \rightarrow \infty} \|u_x^j(v_n^{(1)}+v_n^{(2)})-u_x^j(v_n^{(2)})\|_{H^1(|(x,z)|>N_n)} =0, \\
& & \lim_{n \rightarrow \infty} \|u_y^j(v_n^{(1)}+v_n^{(2)})-u_y^j(v_n^{(2)})\|_{H^1(|(x,z)|>N_n)} =0, \\
& & \lim_{n \rightarrow \infty} \|u_z^j(v_n^{(1)}+v_n^{(2)})-u_z^j(v_n^{(2)})\|_{H^1(|(x,z)|>N_n)} =0
\end{eqnarray*}
hold for each $j\in {\mathbb N}_0$.
\end{lemma}
{\bf Proof.} This result is proved by mathematical induction.

Observe that
\begin{eqnarray*}
\|u_x^0(v_n^{(1)}+v_n^{(2)})-u_x^0(v_n^{(1)})\|_{H^1(|(x,z)|<2R)} 
& = & \|u_x^0(v_n^{(2)})\|_{H^1(|(x,z)|<2R)} \\
& = & \left\|\FF^{-1}\left[\frac{-k_1^2\cosh|k|y}{|k|\sinh |k|}\hat{v}_n^{(2)}\right]\right\|_{H^1(|(x,z)|<2R)} \\
& \rightarrow & 0, \\
\\
\|u_x^0(v_n^{(1)}+v_n^{(2)})-u_x^0(v_n^{(2)})\|_{H^1(|(x,z)|>N_n)} 
& = & \|u_x^0(v_n^{(1)})\|_{H^1(|(x,z)|>N_n)} \\
& = & \left\|\FF^{-1}\left[\frac{-k_1^2\cosh|k|y}{|k|\sinh |k|}\hat{v}_n^{(1)}\right]\right\|_{H^1(|(x,z)|>N_n)} \\
& \rightarrow & 0
\end{eqnarray*}
as $n \rightarrow \infty$ according to Theorem \ref{Integral operators 2}.

Suppose that the result holds for all $i \leq j$. It follows that
\begin{eqnarray*}
\lefteqn{\|v_n^{(1)}(u_x^i(v_n^{(1)}+v_n^{(2)})-u_x^i(v_n^{(1)}))\|_1} \quad\\
& & \leq c\|v_n^{(1)}\|_{1,\infty}
\|u_x^i(v_n^{(1)}+v_n^{(2)})-u_x^i(v_n^{(1)})\|_{H^1(|(x,z)|<2R)} \\
& & \leq c\|v_n^{(1)}\|_{5/2}
\|u_x^i(v_n^{(1)}+v_n^{(2)})-u_x^i(v_n^{(1)})\|_{H^1(|(x,z)|<2R)} \\
& & \rightarrow 0
\end{eqnarray*}
and
\begin{eqnarray*}
\lefteqn{\|(v_n^{(1)})_x(u_y^i(v_n^{(1)}+v_n^{(2)})-u_y^i(v_n^{(1)}))\|_1}\quad \\
& & \leq c\big(\|v_n^{(1)}\|_{1,\infty}
\|u_y^i(v_n^{(1)}+v_n^{(2)})-u_y^i(v_n^{(1)})\|_{H^1(|(x,z)|<2R)} \\
& & \qquad\mbox{}
+ \|(v_n^{(1)})_{xx}(u_y^i(v_n^{(1)}+v_n^{(2)})-u_y^i(v_n^{(1)}))\|_{L^2(|(x,z)|<2R)} \\
& & \qquad\mbox{}
+ \|(v_n^{(1)})_{xz}(u_y^i(v_n^{(1)}+v_n^{(2)})-u_y^i(v_n^{(1)}))\|_{L^2(|(x,z)|<2R)}\big) \\
& & \leq c\big(\|v_n^{(1)}\|_{1,\infty}
\|u_y^i(v_n^{(1)}+v_n^{(2)})-u_y^i(v_n^{(1)})\|_{H^1(|(x,z)|<2R)} \\
& & \qquad\mbox{}
+ \|v_n^{(1)}\|_{W^{2,4}({\mathbb R}^2)} \|u_y^i(v_n^{(1)}+v_n^{(2)})-u_y^i(v_n^{(1)})\|_{L^4(|(x,z)|<2R)}\big) \\
& & \leq c\|v_n^{(1)}\|_{5/2}
\|u_x^i(v_n^{(1)}+v_n^{(2)})-u_x^i(v_n^{(1)})\|_{H^1(|(x,z)|<2R)}\\
& & \rightarrow 0
\end{eqnarray*}
and similarly 
$$\left\|\left\{\begin{array}{c}\!\! v_n^{(2)} \\ \!\!(v_n^{(2)})_x \end{array}\!\!\right\}
(u_x^i(v_n^{(1)}+v_n^{(2)})-u_x^i(v_n^{(2)}))\right\|_1 \rightarrow 0$$
as $n \rightarrow \infty$, so that
\begin{eqnarray}
F_1^{j+1}(v_n^{(1)}+v_n^{(2)})
& = & -(v_n^{(1)}+v_n^{(2)})u_x^j(v_n^{(1)}+v_n^{(2)}) + y((v_n^{(1)})_x+(v_n^{(2)})_x)u_y^j(v_n^{(1)}+v_n^{(2)}) \nonumber \\
& = & -v_n^{(1)} u_x^j(v_n^{(1)}) - v_n^{(2)} u_x^j(v_n^{(2)}) + y (v_n^{(1)})_xu_y^j(v_n^{(1)})
+y(v_n^{(2)})_xu_y^j(v_n^{(2)}) + \underline{o}(1)\nonumber \\
& = & F_1^{j+1}(v_n^{(1)}) + F_1^{j+1}(v_n^{(2)}) + \underline{o}(1), \label{F result 1}
\end{eqnarray}
where the symbol $\underline{o}(1)$ denotes a quantity which converges to zero as $n \rightarrow
\infty$ in $H^1(\Sigma)$. A similar argument shows that
\begin{eqnarray}
F_2^{j+1}(v_n^{(1)}+v_n^{(2)}) & = & F_2^{j+1}(v_n^{(1)}) + F_2^{j+1}(v_n^{(2)}) + \underline{o}(1), 
\label{F result 2}\\
F_3^{j+1}(v_n^{(1)}+v_n^{(2)}) & = & F_3^{j+1}(v_n^{(1)}) + F_3^{j+1}(v_n^{(2)}) + \underline{o}(1),
\label{F result 3}
\end{eqnarray}
where we have used the result that
\begin{eqnarray*}
(v_n^{(1)}+v_n^{(2)})^\ell & = & (v_n^{(1)})^\ell + (v_n^{(2)})^\ell, \\
((v_n^{(1)})_x+(v_n^{(2)})_x)^\ell & = & (v_n^{(1)})_x^\ell + (v_n^{(2)})_x^\ell, \\
((v_n^{(1)})_z+(v_n^{(2)})_z)^\ell & = & (v_n^{(1)})_z^\ell + (v_n^{(2)})_z^\ell
\end{eqnarray*}
for each $\ell \in {\mathbb N}$ (since $v_n^{(1)}$ and $v_n^{(2)}$ have disjoint supports).

Notice that
\begin{eqnarray*}
\lefteqn{\|u_x^{j+1}(v_n^{(1)}+v_n^{(2)})-u_x^{j+1}(v_n^{(1)})\|_{L^2(|(x,z)|<2R)}} \\
& \leq & \left\|\FF\left[-\!\!\int_0^1 k_1^2 G \FF[F_1^{j+1}(v_n^{(1)}+v_n^{(2)})-F_1^{j+1}(v_n^{(1)})]
\dtildey \right]\right\|_{L^2(|(x,z)|<2R)}\\
& & \mbox{}+\left\|\FF\left[-\!\!\int_0^1 k_1k_2 G \FF[F_2^{j+1}(v_n^{(1)}+v_n^{(2)})-F_2^{j+1}(v_n^{(1)})]
\dtildey \right]\right\|_{L^2(|(x,z)|<2R)} \\
& & \mbox{}+\left\|\FF\left[\int_0^1 \i k_1 H_y \FF[F_3^{j+1}(v_n^{(1)}+v_n^{(2)})-F_3^{j+1}(v_n^{(1)})]
\dtildey \right]\right\|_{L^2(|(x,z)|<2R)} \\
& \leq & \left\|\FF\left[-\!\!\int_0^1 k_1^2 G \FF[F_1^{j+1}(v_n^{(2)})]\dtildey\right]\right\|_{L^2(|(x,z)<2R)}
+ \left\|\FF\left[-\!\!\int_0^1 k_1^2 G \FF[\underline{o}(1)]\dtildey\right]\right\|_{L^2(|(x,z)<2R)} \\
& & \mbox{}+\left\|\FF\left[-\!\!\int_0^1 k_1k_2 G \FF[F_2^{j+1}(v_n^{(2)})]\dtildey\right]\right\|_{L^2(|(x,z)<2R)}
+ \left\|\FF\left[-\!\!\int_0^1 k_1k_2 G \FF[\underline{o}(1)]\dtildey\right]\right\|_{L^2(|(x,z)<2R)} \\
& & \mbox{}+\left\|\FF\left[\int_0^1 \i k_1 H_y \FF[F_3^{j+1}(v_n^{(2)})]\dtildey\right]\right\|_{L^2(|(x,z)<2R)}
+ \left\|\FF\left[\int_0^1\i k_1 H_y \FF[\underline{o}(1)]\dtildey\right]\right\|_{L^2(|(x,z)<2R)} \\
& \rightarrow & 0
\end{eqnarray*}
as $n \rightarrow \infty$ because
$$\supp F_j^{j+1}(v_n^{(2)}) \subset \{(x,y,z) \in \Sigma: |(x,z)| >N_n\}, \qquad
j=1,2,3$$
and
$$\left\|\FF\left[\int_0^1\left\{\begin{array}{c} \!\!-k_1^2\\ \!\!-k_1k_2\end{array}\right\}
 G \FF[\underline{o}(1)]\dtildey\right]\right\|_{L^2(|(x,z)<2R)},\ 
 \left\|\FF\left[\int_0^1\i k_1 H_y \FF[\underline{o}(1)]\dtildey\right]\right\|_{L^2(|(x,z)<2R)}
 =o(1).$$
Similar calculations show that
$$\|u_x^{j+1}(v_n^{(1)}+v_n^{(2)})-u_x^{j+1}(v_n^{(2)})\|_{L^2(|(x,z)|>N_n)} \rightarrow 0$$
and yield the estimates
$$
\left\|\left\{\begin{array}{c}\!\partial_x\! \\ \!\partial_z\!\end{array}\right\}
(u_x^{j+1}(v_n^{(1)}+v_n^{(2)})-u_x^{j+1}(v_n^{(2)}))\right\|_{L^2(|(x,z)|<2R)} \rightarrow 0,
$$
$$
\left\|\left\{\begin{array}{c}\!\partial_x\! \\ \!\partial_z\!\end{array}\right\}
(u_x^{j+1}(v_n^{(1)}+v_n^{(2)})-u_x^{j+1}(v_n^{(2)}))\right\|_{L^2(|(x,z)|>N_n)} \rightarrow 0
$$
as $n \rightarrow \infty$.

Finally, observe that
\begin{eqnarray*}
\lefteqn{\|u_{xy}^{j+1}(v_n^{(1)}+v_n^{(2)})-u_{xy}^{j+1}(v_n^{(1)})\|_{L^2(|(x,z)|<2R)}} \quad \\
& \leq & \left\|\FF\left[-\!\!\int_0^1 k_1^2 G_y \FF[F_1^{j+1}(v_n^{(1)}+v_n^{(2)})-F_1^{j+1}(v_n^{(1)})]
\dtildey \right]\right\|_{L^2(|(x,z)|<2R)}\\
& & \mbox{}+\left\|\FF\left[-\!\!\int_0^1 k_1k_2 G_y \FF[F_2^{j+1}(v_n^{(1)}+v_n^{(2)})-F_2^{j+1}(v_n^{(1)})]
\dtildey \right]\right\|_{L^2(|(x,z)|<2R)} \\
& & \mbox{}+\left\|\FF\left[\int_0^1 \i k_1 |k|^2H \FF[F_3^{j+1}(v_n^{(1)}+v_n^{(2)})-F_3^{j+1}(v_n^{(1)})]
\dtildey \right]\right\|_{L^2(|(x,z)|<2R)} \\
& & \mbox{} + \|\partial_x(F_3^{j+1}(v_n^{(1)}+v_n^{(2)})-F_3^{j+1}(v_n^{(1)}))\|_{L^2(|(x,z)|<2R)}.
\end{eqnarray*}
The argument given above shows that the first three terms on the right-hand side of this inequality
are $o(1)$, while the fourth is equal to
$$\underbrace{\|F_{3x}^{j+1}(v_n^{(2)})\|_{L^2(|(x,z)|<2R)}}_{\textstyle = 0}+o(1).$$
A similar calculation shows that
$$\|u_{xy}^{j+1}(v_n^{(1)}+v_n^{(2)}) - u_{xy}^{j+1}(v_n^{(1)})\|_{L^2(|(x,z)|>N_n)}\rightarrow 0$$
as $n \rightarrow \infty$.

Altogether these calculations show that
$$\|u_x^{j+1}(v_n^{(1)}+v_n^{(2)})-u_x^{j+1}(v_n^{(1)})\|_{H^1(|(x,z)|<2R)} \rightarrow 0,$$
$$\|u_x^{j+1}(v_n^{(1)}+v_n^{(2)})-u_x^{j+1}(v_n^{(2)})\|_{H^1(|(x,z)|>N_n)} \rightarrow 0$$
as $n \rightarrow \infty$, and 
the corresponding results for $u_y^{j+1}$ and $u_z^{j+1}$ are obtained in a similar fashion.\qed

\begin{corollary}
Every sequence $\{v_n^{(2)}\}$ in $U$ with the property that
$\supp v_n^{(2)} \subset {\mathbb R}^2 \sm B_{N_n}(0)$ satisfies the estimates
\begin{eqnarray*}
& & \lim_{n \rightarrow \infty} \|u_x^j(v_n^{(2)})\|_{H^1(|(x,z)|<2R)} = 0, \\
& & \lim_{n \rightarrow \infty} \|u_y^j(v_n^{(2)})\|_{H^1(|(x,z)|<2R)} =0, \\
& & \lim_{n \rightarrow \infty} \|u_z^j(v_n^{(2)})\|_{H^1(|(x,z)|<2R)} =0
\end{eqnarray*}
for each $j \in {\mathbb N}_0$.
\end{corollary}
{\bf Proof.} This result follows directly from Lemma \ref{Key splitting lemma}
with $v_n{^{(1)}}=0$, $m \in {\mathbb N}$.\qed

We now have all the ingredients to prove the final result.

\begin{theorem} \label{Key splitting theorem}
The estimates
$$\lim_{n \rightarrow \infty} \Big(\LL(v_n^{(1)}+v_n^{(2)}) - \LL(v_n^{(1)})-\LL(v_n^{(2)}) \Big)=0,$$
$$\lim_{n \rightarrow \infty} \Big\|\LL^\prime(v_n^{(1)}+v_n^{(2)}) - \LL^\prime(v_n^{(1)})-\LL^\prime(v_n^{(2)}) \Big\|_1=0,$$
$$\lim_{i \rightarrow \infty} \langle \LL^\prime(v_n^{(2)}),v_n^{(1)} \rangle_0 =0$$
hold for all sequences $\{v_n^{(1)}\}$ and $\{v_n^{(2)}\}$ in $U$  which have 
the properties that
$\supp v_n^{(1)} \subset B_{2R}(0)$, $\supp v_n^{(2)} \subset {\mathbb R}^2 \sm B_{N_n}(0)$
and $\sup \|v_n^{(1)}+v_n^{(2)}\|_3<M$.
\end{theorem}
{\bf Proof.} Recall that
$$\LL(\eta) = \sum_{j=2}^\infty \LL_i(\eta), \qquad \LL_j(\eta) = -\frac{1}{2}\int_{{\mathbb R}^2}
u_x^{j-2}|_{y=1}\eta\dx\dz,$$
where the series converges uniformly in $U$.
Choose $\tilde{\varepsilon}>0$ and select $N \geq 2$ large enough so that
$$\left| \sum_{j=N+1}^\infty \LL_i(\eta) \right| < \tilde{\varepsilon}, \qquad \eta \in U.$$
Observe that
\begin{eqnarray*}
\lefteqn{\LL_j(\eta_n^{(1)}+\eta_n^{(2)}) - \LL_j(\eta_n^{(1)}) - \LL_j(\eta_n^{(2)})} \qquad\\
& = & -\frac{1}{2}\int_{{\mathbb R}^2} \eta_n^{(1)}(u_x^{j-2}(\eta_n^{(1)}+\eta_n^{(2)})-u_x^{j-2}(\eta_n^{(1)}))|_{y=1}\dx\dz \\
& & \qquad\mbox{}-\frac{1}{2}\int_{{\mathbb R}^2} \eta_n^{(2)}(u_x^{j-2}(\eta_n^{(1)}+\eta_n^{(2)})-u_x^{j-2}(\eta_n^{(2)}))|_{y=1}\dx\dz,
\end{eqnarray*}
whereby
\begin{eqnarray*}
\lefteqn{|\LL_j(\eta_n^{(1)}+\eta_n^{(2)})-\LL_j(\eta_n^{(1)})-\LL_j(\eta_n^{(2)})|} \qquad\\
& \leq & \uunderbrace{\|\eta_n^{(1)}\|_3}_{\textstyle O(1)}
\uunderbrace{\|u_x^{j-2}(\eta_n^{(1)}+\eta_n^{(2)})-u_x^{j-2}(\eta_n^{(1)})\|_{H^1(|(x,z)|<2R)}}_{\textstyle o(1)} \\
& & \mbox{}+\uunderbrace{\|\eta_n^{(2)}\|_3}_{\textstyle O(1)}
\uunderbrace{\|u_x^{j-2}(\eta_n^{(1)}+\eta_n^{(2)})-u_x^{j-2}(\eta_n^{(1)})\|_{H^1(|(x,z)|>N_n)}}_{\textstyle o(1)} \\
& = & o(1)
\end{eqnarray*}
as $n \rightarrow \infty$ for $m=2, \ldots, N$, in which Lemma \ref{Key splitting lemma} has been
used. It follows that
\begin{eqnarray*}
\lefteqn{|\LL(\eta_n^{(1)}+\eta_n^{(2)})-\LL(\eta_n^{(1)})-\LL(\eta_n^{(2)})|} \\
& = & \underbrace{\sum_{j=2}^N \{ \LL_j(\eta_n^{(1)}+\eta_n^{(2)})-\LL_j(\eta_n^{(1)})-\LL_j(\eta_n^{(2)}) \}}_{\textstyle o(1)}
+ \underbrace{\sum_{j=N+1}^\infty \{ \LL_j(\eta_n^{(1)}+\eta_n^{(2)})-\LL_j(\eta_n^{(1)})-\LL_j(\eta_n^{(2)}) \}}_{\textstyle O(\tilde{\varepsilon})}
\end{eqnarray*}
as $n \rightarrow \infty$.

The other estimates are obtained by applying the same argument to the formula for $\LL^\prime$
given in Theorem \ref{Formulae for the gradients}.\qed

\end{document}